\documentclass[11pt]{article}

\usepackage{amsthm}
\usepackage{amsmath, bm}
\usepackage{amssymb}
\usepackage{natbib}
\usepackage{geometry}
\usepackage{subfig}
\usepackage{graphicx}
\usepackage{tikz}
\usepackage{multirow}
\usepackage{listings}
\usepackage{xcolor}
\usepackage{authblk}
\usepackage{float}
\usepackage{url}
\usepackage{calc}
\usepackage{stmaryrd}

\usepackage{subfiles}
\usepackage{indentfirst}
\usepackage{fullpage}

\newtheorem{Def}{Definition}[section]
\newtheorem{Thm}{Theorem}[section]

\newtheorem{theorem}{Theorem}[section]
\newtheorem{lemma}[theorem]{Lemma}

\newcommand{\diag}[1]{{\rm diag}\LRp{#1}}
\newcommand{\td}[2]{\frac{{\rm d}#1}{{\rm d}{\rm #2}}}

\newcommand{\LRp}[1]{\left( #1 \right)}

\newcommand{\LRb}[1]{\left| #1 \right|}
\newcommand{\LRc}[1]{\left\{ #1 \right\}}

\newcommand{\avg}[1] {\ensuremath{\LRc{\!\{#1\}\!}}}

\newcommand{\fnt}[1]{\bm{\mathsf{ #1}}}

\title{Entropy stable discontinuous Galerkin methods for the shallow water equations with subcell positivity preservation}
\author[1]{Xinhui Wu}
\author[2]{Nathaniel Trask}
\author[1]{Jesse Chan}
\affil[1]{Department of Computational and Applied Mathematics, Rice University, Houston, TX}
\affil[2]{Center for Computing Research, Sandia National Laboratories, Albuquerque, NM}
\date{}

\begin{document}
\maketitle

\begin{abstract}
High order schemes are known to be unstable in the presence of shock discontinuities or under-resolved solution features, and have traditionally required additional filtering, limiting, or artificial viscosity to avoid solution blow up. Entropy stable schemes address this instability by ensuring that physically relevant solutions satisfy a semi-discrete entropy inequality independently of discretization parameters. However, additional measures must be taken to ensure that solutions satisfy physical constraints such as positivity. In this work, we present a high order entropy stable discontinuous Galerkin (ESDG) method for the nonlinear shallow water equations (SWE) on two-dimensional (2D) triangular meshes which preserves the positivity of the water heights. The scheme combines a low order positivity preserving method with a high order entropy stable method using convex limiting. This method is entropy stable and well-balanced for fitted meshes with continuous bathymetry profiles. 
\end{abstract}

\section{Introduction}
\label{sec:1}
\begin{subequations}
\numberwithin{equation}{section}

The shallow water equations (SWE) are used to model fluid flows in rivers, open channels, lakes, and coastal regions. One of the most prominent practical applications of the shallow water equations is the numerical prediction of storm surges under extreme weather conditions like hurricanes near coastal regions \cite{dawson2011discontinuous, akbar2013hybrid, wintermeyer2018entropy}. All of these flows share the fact that vertical scales of the motion are much smaller than horizontal scales. Under certain assumptions \cite{dawson2008shallow}, the incompressible Navier–Stokes equations of fluid dynamics can be simplified to the 2D shallow water equations:
\label{eq:SWE}
\begin{align}
h_t+(hu)_x+(hv)_y&=0,         			 \label{eq:SWEa} \\
(hu)_t+(hu^2+gh^2/2)_x+(huv)_y&=-ghb_x,  \label{eq:SWEb} \\ 
(hv)_t+(huv)_x+(hv^2+gh^2/2)_y&=-ghb_x, \label{eq:SWEc}
\end{align}
where $h = h(x,y,t)$ denotes the height of the water measured from the bottom, $u = u(x,y,t)$ and $v = v(x,y,t)$ denote the velocity in the $x$ and $y$ directions respectively, $g$ denotes the gravitational constant. The bathymetry height is denoted by $b=b(x,y)$. The subscript $(.)_t$ denotes the time derivative, while the subscripts $(.)_x$ or $(.)_y$ denote directional derivatives along the $x$ and $y$ directions, respectively.

Researches on numerical methods for the solution of the shallow water equations have achieved considerable attention in the past two decades. Common numerical methods including the finite difference, the finite volume, and the finite element method have been applied to the shallow water system. This paper focuses on the discontinuous Galerkin (DG) method. Discontinuous Galerkin methods combine advantages of both finite element and finite volume methods, and have been successfully applied to a variety of settings \cite{hesthaven2007nodal}. Discontinuous Galerkin methods can achieve high order accuracy, parallel efficiency, and enjoy more straightforward hp-refinement \cite{dawson2002discontinuous, giraldo2002nodal}. However, high order discontinuous Galerkin methods also introduce issues of stability for problems with shocks or under-resolved solutions. Recent works on entropy stable high order discontinuous Galerkin methods \cite{wintermeyer2017entropy, chan2018discretely, wen2020entropy} provide a way to address such instabilities. Entropy stable discontinuous Galerkin methods can also be extended to curved meshes, which can be necessary when dealing with complex geometries.

One of the important numerical properties for the shallow water equations is the preservation of steady state solutions, also known as the ``lake at rest'' condition \cite{wintermeyer2018entropy}:
\begin{align}
H=\text{constant},\qquad u=v=0.
\end{align}
A numerical scheme which preserves this steady state is said to be well-balanced \cite{gassner2016well,noelle2007high}. Schemes that are not well-balanced can generate spurious waves in the presence of varying bottom topographies. Accomplishing this discretely can be challenging because special discretizations of the source terms are required \cite{wintermeyer2017entropy, fjordholm2011well, xing2014exactly}.

Another important issue often experienced in numerical simulations of the shallow water equations is the treatment of dry areas, where the water is absent. Dry areas are common in realistic applications of the shallow water equations, such as dam breaks, flooding, and shoreline simulations. Some numerical methods can fail when encountering zero or negative water heights. Various wetting and drying treatments for the shallow water equations have been developed \cite{nielsen2003parameters, xing2010positivity, audusse2005well, xing2013positivity, jiang2016invariant}. Some use mesh adaption to track dry fronts and modify the meshes to fit the wet area. Others exclude dry elements from computation and include them when they are wet again; however, this approach may introduce oscillations and a loss of mass and momentum conservation. Positivity-preserving limiters that average the water height in partially dry elements have been also explored in \cite{xing2010positivity}. 

In this paper, we introduce a high order entropy stable numerical method that preserves the positivity of water heights for the shallow water equations through convex limiting. The sections in this paper are organized as follows. We first present a high order entropy stable numerical scheme and entropy stable numerical fluxes for the shallow water equation in 1D and on 2D triangular meshes \cite{chan2018discretely, chen2020review, chen2017entropy, wu2021high, wu2021entropy} in Sec. \ref{sec:2}. Then, we present a low order positivity preserving method \cite{guermond2016invariant} which preserves the positivity of the water heights by using a graph viscosity in Sec. \ref{sec:3}. We combine the high order scheme with the low order scheme through a technique called convex limiting \cite{guermond2019invariant} in Sec. \ref{sec:4}. If the high order scheme produces negative water heights, the convex limiting procedure blends them with the low order solution to preserve positivity while maintaining high order accuracy in wet areas. Our approach is also different in the sense that we only limit in dry areas without limiting near shocks. We also prove this scheme is well-balanced for meshes which are fitted to dry-wet interfaces. We demonstrate our new method with five numerical experiments that illustrate different properties of our new method in Sec. \ref{sec:5}. Conclusions are presented in Sec. \ref{sec:6}.

\end{subequations}

\section{High order entropy stable summation-by-parts formulations}
\label{sec:2}
\begin{subequations}
\numberwithin{equation}{section}
\subsection{Mathematical assumption and notations}
We first introduce some underlying mathematical assumptions and notations for our DG method. For consistency, we reuse notations from \cite{chan2018discretely}, with slight modifications to provide a cleaner discrete formulation. We denote the triangular reference element by $\widehat{D}$ with boundary $\partial{\widehat{D}}$. In 2D, the reference element is a triangle with vertices at $(-1,-1)$, $(-1,1)$ and $(1,-1)$.  We use $\widehat{n}_i$ to represent the $i$th component of the outward normal vector scaled by the face Jacobian on the boundary of the reference element. 
The space of polynomials up to degree $N$ on the reference element is defined as
\begin{align}
P^N(\widehat{D}) = \{ \widehat{x}^i\widehat{y}^j,\quad (\widehat{x},\widehat{y}) \in \widehat{D}, \quad 0 \leq i+j \leq N \}.
\end{align}
Finally, we denote the dimension of $P^N$ as $N_p = \rm{dim}(P^N(\widehat{D}))$.

\subsection{SBP quadrature rules}
We require SBP-quadrature rules to satisfy certain properties detailed in \cite{chen2017entropy, wu2021high}. In our numerical experiments, we consider two sets of SBP quadrature points. The first uses 1D Gauss-Legendre quadrature on the edges while the second uses 1D Gauss-Lobatto quadrature on edges. These are shown in the Fig. \ref{fig:SBP_GLE_nodes} and \ref{fig:SBP_GLO_nodes}, respectively. The Gauss-Legendre SBP quadrature rules were introduced in \cite{chen2017entropy}, and the Gauss-Lobatto SBP quadrature rules were introduced in \cite{wu2021high}.
\begin{figure}[H]
\begin{center}
\includegraphics[width=.23\textwidth]{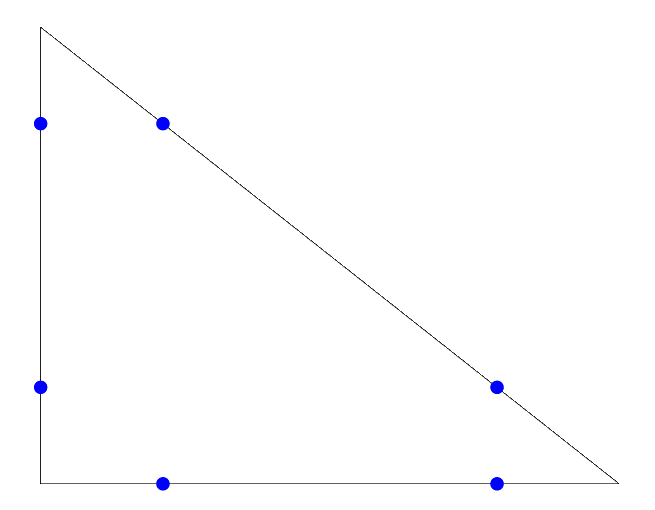}
\includegraphics[width=.23\textwidth]{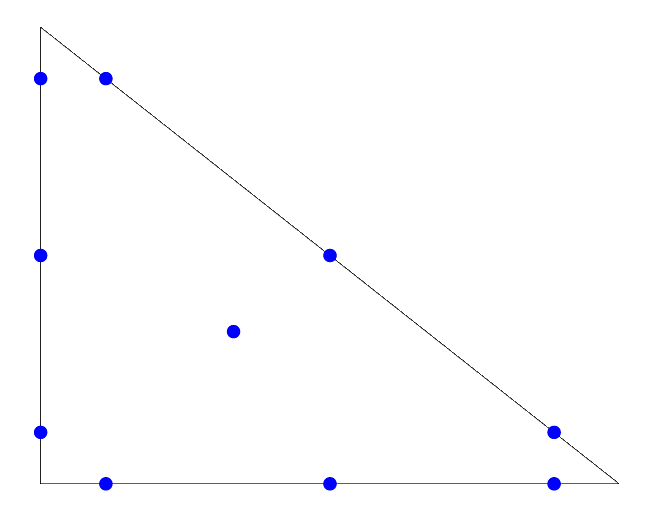}
\includegraphics[width=.23\textwidth]{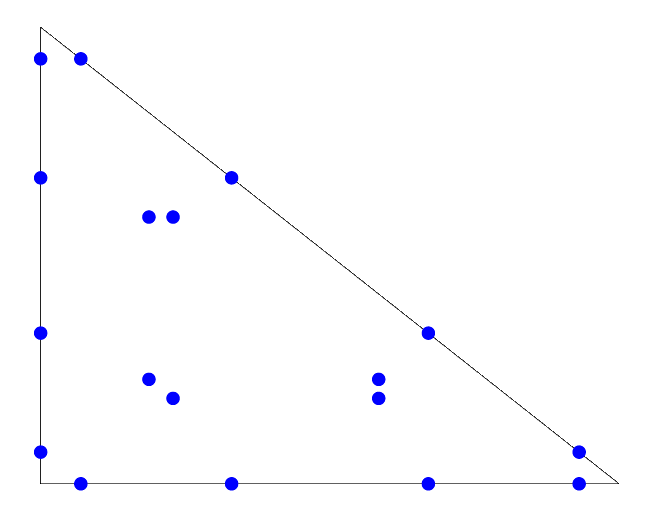}
\includegraphics[width=.23\textwidth]{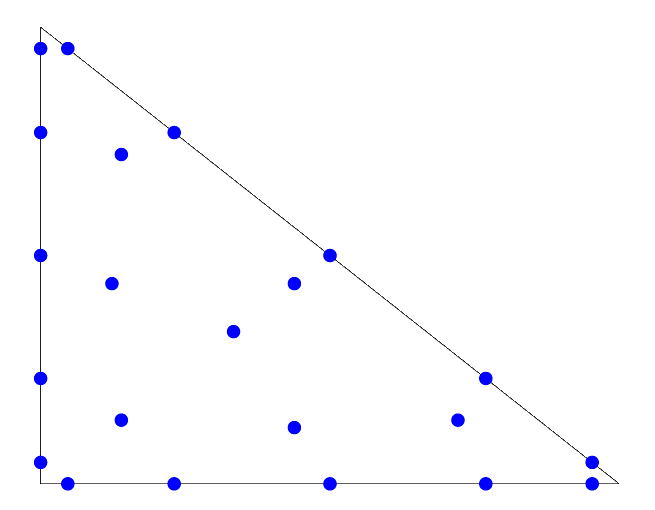}
\caption{SBP nodes based on Gauss-Legendre quadrature for $N = 1,2,3,4$}
\label{fig:SBP_GLE_nodes}
\end{center}
\end{figure}
\begin{figure}[H]
\begin{center}
\includegraphics[width=.23\textwidth]{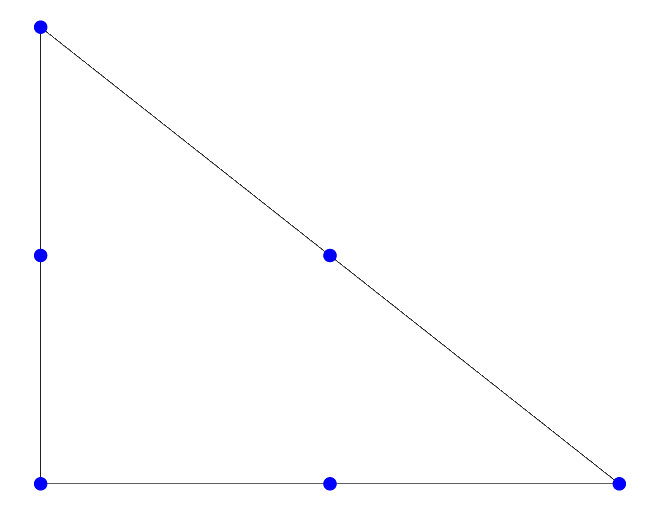}
\includegraphics[width=.23\textwidth]{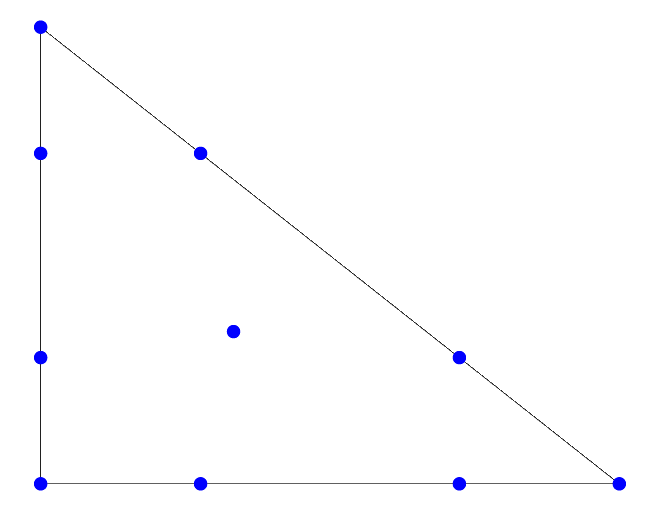}
\includegraphics[width=.23\textwidth]{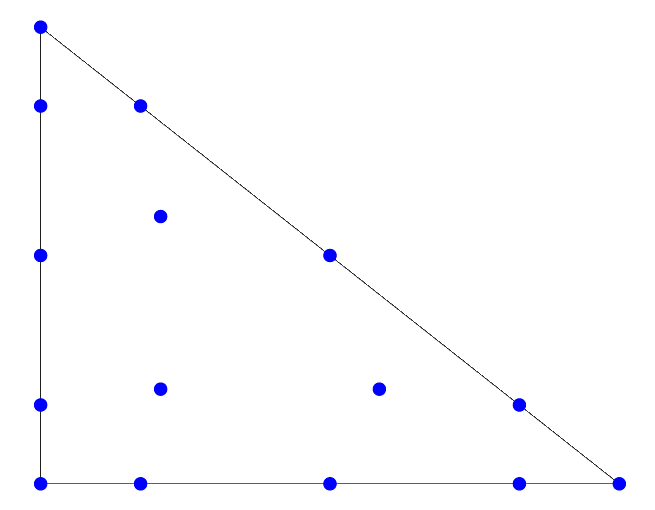}
\includegraphics[width=.23\textwidth]{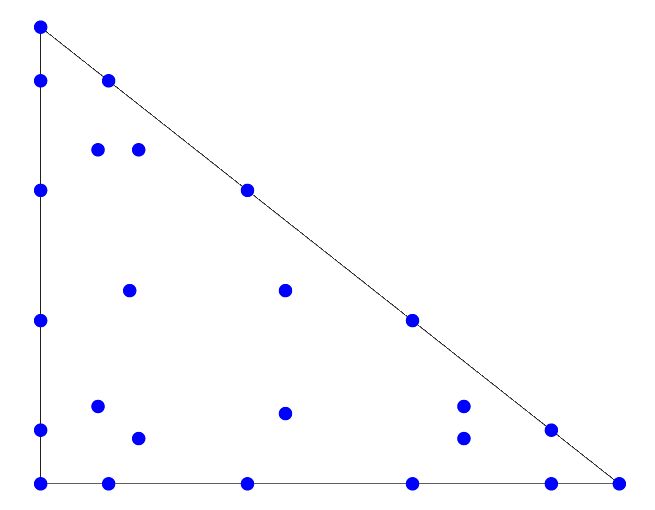}
\caption{SBP nodes based on Gauss-Lobatto quadrature for $N = 1,2,3,4$}
\label{fig:SBP_GLO_nodes}
\end{center}
\end{figure}

\subsection{Multi-dimensional SBP operators}
We now describe properties satisfied by nodal SBP operators. We assume we are given a degree $2N-1$ quadrature rule for the polynomial space $P^N(\widehat{D})$, with $N_q$ nodes, $\{\fnt{x_i}\}_{i=1}^{N_q}$, and positive weights $\{w_i\}_{i=1}^{N_q}$. The face quadrature rule is degree $2N$ and is embedded in the volume nodes. The nodal values of the function $u(x)$ on the quadrature points is denoted by:
\begin{align}
\fnt{u} = [u(\fnt{x}_1),  ... , u(\fnt{x}_{N_q})]^T.
\end{align}
The SBP mass matrix is defined as a diagonal matrix with quadrature weights on the diagonal.
\begin{align}
\fnt{M} = \diag{[w_1, ..., w_{N_q}]}.
\end{align}
We define $\fnt{D}_x$ and $\fnt{D}_y$ to be nodal differentiation matrices associated with the $x$ and $y$ derivatives \cite{chen2020review, wu2021high}. We also define $\fnt{B}_x$ and $\fnt{B}_y$ to be diagonal surface matrices whose entries are surface quadrature weights scaled by the $x$ and $y$ components of the outward normal vector.

We have the following definition \cite{chen2017entropy}.
\begin{Def}{}
Consider the diagonal mass matrix consisting of quadrature weights
\begin{align}
\fnt{M} = \diag{[w_1, ..., w_{N_q}]}.
\end{align}
A 2D operator $\fnt{Q}_i$ is said to be high order accurate and have SBP property if for $i=x,y$, the following properties hold:
\begin{itemize}
\item Let $\fnt{D}_i = \fnt{M}^{-1}\fnt{Q}_i$. Then $
(\fnt{D}_i\fnt{u})_j=\frac{\partial \fnt{u}}{\partial x_i}\bigg\rvert_{x = x_j}
$ for any $\fnt{u}\in P^N(\widehat{D})$.
\item $\fnt{Q}_i+(\fnt{Q}_i)^T = \fnt{B}_i.$
\end{itemize}
\end{Def}

Given this equivalence, we can now construct an entropy conservative DG-SBP form:
\begin{align}
\fnt{M}\td{\fnt{u}}{t} + \sum_{i=x,y} \LRp{2\fnt{Q}_i\circ \fnt{F}_i}\fnt{1} + \fnt{B}_i\LRp{\bm{f}_{S,i}\LRp{\fnt{u}^+, \fnt{u}}-\fnt{f}_i(\bm{u})} = \fnt{S},
\end{align}
where the matrix operations $\circ$ denotes the Hadamard product. $\fnt{F}_i$ denotes the flux matrix and each entry $(\fnt{F}_i)_{jk} = \bm{f}_i(\bm{u}_j, \bm{u}_k)$, where $\bm{f}_i$ is an entropy conservative numerical flux that we will introduce in the next section. $\fnt{S}$ denotes the source term. The construction of the high order operators are presented in \cite{wu2021entropy}. The stability analysis of the SBP formulation follows from the results from \cite{chen2017entropy}, and the extension to curved elements can be found in \cite{crean2018entropy}. 

\subsection{Entropy conservative fluxes for the shallow water equations}
The entropy function for the shallow water equations is the total energy of the system \cite{fjordholm2011well, wintermeyer2017entropy}:
\begin{align}
S(\bm{u}) = \frac{1}{2}h(u^2+v^2) + \frac{1}{2}gh^2+ghb.
\label{eq:entropy_function}
\end{align}
We also define the entropy variable $\bm{v} = S'(\bm{u})$. The convexity of the entropy function guarantees that the mapping between $\bm{u}$ and $\bm{v}$ is invertible.  The entropy variables for the shallow water equations are given explicitly as:
\begin{align}
v_1 &= \frac{\partial S}{\partial h}  = g(h+b) - \frac{1}{2}u^2 - \frac{1}{2}v^2,\\
v_2 &= \frac{\partial S}{\partial (hu)} = u,\\
v_3 &= \frac{\partial S}{\partial (hv)} = v.
\end{align}
It can be shown as in \cite{mock1980systems} that there exists an entropy flux function $F(\bm{u})$ and entropy potential $\psi(\bm{u})$ such that
\begin{align}
\bm{v}(\bm{u})^T\frac{\partial \bm{f}}{\partial \bm{u}} = \frac{\partial F(\bm{u})^T}{\partial \bm{u}}, \hspace{.5cm}
\psi(\bm{u}) = \bm{v}(\bm{u})^T\bm{f}(\bm{u}) - F(\bm{u}), \hspace{0.5cm} \psi'(\bm{u}) = \bm{f}(\bm{u}).
\end{align}

We introduce numerical fluxes for the shallow water equations and describe an entropy conservation discrete formulation \cite{fisher2013high, gassner2018br1, chen2017entropy, gassner2016split}. To construct the entropy stable scheme in higher dimensions, we require entropy conservative fluxes as defined in \cite{tadmor1987numerical}.

\begin{Def}{}
Let $\bm{f}_{S,i}(\bm{u}_L, \bm{u}_R)$ be a bivariate function which is symmetric and consistent with the flux function $\bm{f}_i(\bm{u})$, for $i = 1,...,d$
\begin{align}
\bm{f}_{S,i}(\bm{u},\bm{v}) = \bm{f}_{S,i}(\bm{v},\bm{u}), \hspace{0.5cm}  \bm{f}_{S,i}(\bm{u},\bm{u}) = \bm{f}_i(\bm{u}).
\end{align}
The numerical flux $\bm{f}_{S,i}(\bm{u}_L, \bm{u}_R)$ is entropy conservative if, for entropy variables $\bm{v}_L = \fnt{V}(\bm{u}_L)$ and $\bm{v}_R = \bm{v}(\bm{u}_R)$, $\bm{f}_{S,i}(\bm{u}_L, \bm{u}_R)$ satisfies
\begin{align}
\left( \bm{v}_L - \bm{v}_R\right)^T \bm{f}_{S,i}\left(\bm{u}_L,\bm{u}_R\right) = \psi_{i,L} - \psi_{i,R}, \\
\psi_{i,L} = \psi_i(\bm{v}(\bm{u}_L)), \quad \psi_{i,R} = \psi_i(\bm{v}(\bm{u}_R)).
\end{align}
\label{def:consevative_flux}
\end{Def}
The flux $\bm{f}_{S,i}$ can be used to construct entropy conservative and entropy stable finite volume methods, as well as DG methods. Entropy conservative (EC) fluxes for the 2D shallow water equations \cite{tadmor2003entropy, carpenter2014entropy, gassner2016well, wintermeyer2018entropy} are given by
\begin{align}
\bm{f}_{x,S}\LRp{\bm{u}_L,\bm{u}_R} &=
\begin{bmatrix}
\avg{hu}\\
\avg{hu}\avg{u} + g\avg{h}^2 - \frac{1}{2}g\avg{h^2}\\
\avg{hu}\avg{v}
\end{bmatrix}
\label{eq:SWE_flux1}
, \\
\bm{f}_{y,S}\LRp{\bm{u}_L,\bm{u}_R} &=
\begin{bmatrix}
\avg{hv}\\
\avg{hv}\avg{u}\\
\avg{hv}\avg{v} + g\avg{h}^2 - \frac{1}{2}g\avg{h^2}\\ 
\end{bmatrix}.
\label{eq:SWE_flux2}
\end{align}
\end{subequations}
These fluxes also yield a well-balanced scheme for continuous bathymetry, if the source terms are discretized appropriately \cite{wintermeyer2017entropy}. The source terms $\fnt{S}$ for the shallow water equations with continuous bottom geometry in the discrete form is defined as
\begin{align}
\fnt{S} = -g\begin{bmatrix}
\fnt{0}\\
{\rm diag}\LRp{\fnt{h}}\fnt{Q}_x\fnt{b}\\
{\rm diag}\LRp{\fnt{h}}\fnt{Q}_y\fnt{b}
\end{bmatrix}.
\label{eq:discrete_source_term}
\end{align}
By adding local Lax-Friedrichs interface penalization \cite{wintermeyer2017entropy}, we can make our numerical scheme entropy stable.
\section{Positivity preserving formulation}
\label{sec:3}
\begin{subequations}
\numberwithin{equation}{section}
In this section, we first describe the steps for constructing a low order positivity preserving DG scheme from \cite{guermond2019invariant} for the shallow water equations. Then, we prove that the low order method preserve the positivity of the water heights. Last, we remark on the sparsity of the low order operators.

\subsection{The low order scheme}
We build a low order positivity preserving scheme based on a connectivity graph and the solution of a graph Laplacian problem constructed from the connectivity graph. We present the details of the construction in Appendix A. We assume the same reference element, quadrature points, and discrete operators. We assume the forward Euler time stepping and let $t^n$ denote the time at step $n \in \mathbb{N}$ and $dt$ denote the current time step size. Let $m_i$ denote the $i$th diagonal entry global mass matrix, which is diagonal. Suppose we can construct an inviscid approximation $\fnt{u}^{n+1}$ of the shallow water equations at some time as follow \cite{guermond2019invariant}:
\begin{align}
    \frac{m_i}{dt}(\fnt{u}_i^{n+1} - \fnt{u}_i^{n}) + \sum_{j\in I(i)}\fnt{f}(\fnt{u}_j^n)(\fnt{Q}^{L})_{ij}= m_i\fnt{S}(\fnt{u}_i^{n}),
    \label{eq:low_order_base}
\end{align}
for any $i \in V$, the set of all nodes. For any $i$, the set $I(i)$ is the neighborhood of the node $i$. We assume that the following property holds: $j \in I(i)$ iff $i \in I( j)$. $\fnt{Q}^{L}$ denotes the low order differentiation operator. $\fnt{Q}^{L}$ is described more in details in Sec. \ref{sec:3_sparse}, and its construction is presented in Appendix A. The discrete form of the source term $\fnt{S}$ is defined in Eq. \ref{eq:discrete_source_term}, with replacing the high order differentiation operator with lower order differentiation operator.
We also assume that the entries of the low order differentiation operator, $(\fnt{Q}^{L})_{ij}$, has the following properties:
\begin{align}
    (\fnt{Q}^{L})_{ij} = -(\fnt{Q}^{L})_{ji} \quad {\rm and } \quad \sum_{j\in I(i)}(\fnt{Q}^{L})_{ij} = 0. 
    \label{eq:low_order_req}
\end{align}

To handle shocks and discontinuous in the solutions, we introduce artificial dissipation using the graph Laplacian associated with the connectivity graph $(V, E)$. Let $d^{L,n}_{ij}$ denote the scalar graph viscosity for $(i, j) \in E$ with the following properties:
\begin{align}
    d^{L,n}_{ij} = d^{L,n}_{ji} > 0 \quad {\rm if } \quad i\neq j.
\end{align}
We define $d^{L,n}_{ii} = -\sum_{j \in I(i), j\neq i} d^{L,n}_{ij}$. Now, we can define the low order update $\fnt{u}_i^{L,n+1}$ as follows:
\begin{align}
\frac{m_i}{dt}(\fnt{u}_i^{n+1} - \fnt{u}_i^{n}) + \sum_{j\in I(i)}\bm{f}(\fnt{u}_j^n)(\fnt{Q}^{L})_{ij} -\sum_{j\in I(i) \backslash \{i\}} d^{L,n}_{ij}(\fnt{u}^n_j - \fnt{u}^n_i)= m_i\fnt{S}(\fnt{u}_i^{n}),
\label{eq:low order-base}
\end{align}
for all $1 \leq i \leq N_q$. This scheme is conservative in the following sense \cite{guermond2019invariant}.
\begin{lemma}
Assume that $\bm{S} \equiv \bm{0}$, then the scheme (\ref{eq:low order-base}) is conservative, such that the following identity holds for any $n \in \mathbb{N}$:
\begin{align}
    \sum_{i\in V} m_i\fnt{u}^{L,n+1}_i =  \sum_{i\in V} m_i\fnt{u}^{L,n}_i. 
\end{align}
\end{lemma}
\begin{proof}
See \cite{guermond2019invariant}.
\end{proof}
To enforce the positivity preserving property, we define the graph viscocity coefficients $d^{L,n}_{ij}$ as follows:
\begin{align}
   d^{L,n}_{ij} = {\rm max}(\lambda_{\rm max}(\fnt{u}^n_i, \fnt{u}^n_j) \lVert (\fnt{Q}^{L})_{ij}\rVert_{l^2}, \lambda_{\rm max}(\fnt{u}^n_j, \fnt{u}^n_i)\lVert (\fnt{Q}^{L})_{ji}\rVert_{l^2}),
   \label{eq:low_order_dij}
\end{align}
where 
\begin{align}
    \lambda_{\rm max}(\fnt{u}^n_i, \fnt{u}^n_j)
    \label{eq:shallow_water_wave_speed}
\end{align}
is the upper bound for the maximum wave speed. For the shallow water equations, we use the estimate $\lambda_{\rm max}(\fnt{u}^n_i, \fnt{u}^n_j) = | u_i \bm{n}_{ij} + v_i\bm{n}_{ij}|+\sqrt{gh}$, where $u_i$ and $v_i$ are the $x$ and $y$ velocity respectively at quadrature point $i$ and $\bm{n}_{ij} = (\fnt{Q}^{L})_{ij}/\lVert(\fnt{Q}^{L})_{ij}\rVert$. 

\begin{Thm}
 The numerical scheme (\ref{eq:low order-base}) with $d^{L,n}_{ij}$ defined in (\ref{eq:low_order_dij}) preserves the positivity of the water height for the shallow water equations for appropriate time step size $dt$.
\end{Thm}
\begin{proof}
Since we have $\sum_{j\in I(i)}\bm{f}(\fnt{u}_i^n)(\fnt{Q}^{L})_{ij} = 0$ by property (\ref{eq:low_order_req}), we can rewrite the update on the water height term (\ref{eq:low order-base}) as follows:
\begin{align}
\frac{m_i}{dt}(h_i^{n+1} - h_i^{n}) +\sum_{j\in I(i) \backslash \{i\}} 2d^{L,n}_{ij}h_i^{n} +\LRp{\bm{f}_h(\fnt{u}_j^n) - \bm{f}_h(\fnt{u}_i^n)}(\fnt{Q}^{L})_{ij} -  d^{L,n}_{ij}(h^n_j + h^n_i)= 0,
\label{eq:pos_pf_1}
\end{align}
where $\bm{f}_h$ denotes the flux term that corresponds to the water height. We ignore the source term $\bm{S}(\fnt{u}_i^{n})$ because the topography or friction terms only affect the momentum equations but not water heights. Then we can introduce the intermediate bar states defined by
\begin{align}
\bar{h}_{ij}^{n} = \frac{1}{2}(h_i^n + h_j^n) - \LRp{\bm{f}_h(\fnt{u}_j^n) - \bm{f}_h(\fnt{u}_i^n)}\frac{(\fnt{Q}^{L})_{ij}}{2d^{L,n}_{ij}}.
\label{eq:pos_pf_bar}
\end{align}
Notice that we have $\bar{h}_{ii}^{n} = h_i^n$. We first want to show that if we assume $h_i^n$ is positive for all $i$, then $\bar{h}_{i}^{n}$ is also positive. Notice that, for the shallow water equations, the entropy stable numerical fluxes (\ref{eq:SWE_flux1}) are $\avg{hv}$ for the water height, where $v$ is the velocity of the flow. We can substitute the flux in (\ref{eq:pos_pf_bar}) and rewrite the ``bar state'' as follow:
\begin{align}
\bar{h}_{ij}^{n} = \frac{1}{2}(h_i^n + h_j^n) - \LRp{(\bm{hv})_j^n - (\bm{hv})_i^n}\frac{(\fnt{Q}^{L})_{ij}}{2d^{L,n}_{ij}}.
\label{eq:pos_pf_bar_flux}
\end{align}
By the definition of $d^{L,n}_{ij}$, we know that $\LRb{\bm{v}_i^n\frac{(\fnt{Q}^{L})_{ij}}{2d^{L,n}_{ij}}} = c_i< 1$ and $\LRb{\bm{v}_j^n\frac{(\fnt{Q}^{L})_{ij}}{2d^{L,n}_{ij}}} = c_j< 1$. Since $h_i^n$ and ${h_j^n}$ are positive, we have 
\begin{align}
\bar{h}_{ij}^{n} \geq \frac{1}{2}(h_i^n + h_j^n - c_ih_i^n  - c_jh_j^n )>0
\label{eq:pos_pf_bar_res}
\end{align}
Then we can rewrite the low order update on the water height with the bar state as:
\begin{align}
h_i^{n+1} &= \LRp{1 -  \sum_{j\in I(i) \backslash \{i\}} \frac{2d^{L,n}_{ij}dt}{m_i}} h_i^n + \sum_{j\in I(i) \backslash \{i\}}\frac{2d^{L,n}_{ij}dt}{m_i}\bar{h}_{ij}^{n}\\
&= \LRp{1+2dt\frac{d^{L,n}_{ii}}{m_i}} h_i^n + \sum_{j\in I(i) \backslash \{i\}} \frac{2d^{L,n}_{ij}dt}{m_i}\bar{h}_{ij}^n.
\label{eq:pos_pf_2}
\end{align}
Recall that  $d^{L,n}_{ii} = -\sum_{j \in I(i), j\neq i} d^{L,n}_{ij}$. If time step satisfies that
$dt < \min_i\LRp{\frac{m_i}{2d^{L,n}_{ii}}}$, then, the term $1+2dt\frac{d^{L,n}_{ii}}{m_i}>0$  and the terms in (\ref{eq:pos_pf_2}) form a convex combination of two positive water heights. 
\end{proof}

\subsection{Sparsity of the low order operators}
\label{sec:3_sparse}
One way to build a more accurate lower order numerical scheme is to remove excessive artificial dissipation from the method. In \cite{pazner2021sparse}, the author constructed sparse low order operators for tensor product elements in 2D. In this paper, we construct sparse low order operators on triangular elements in 2D. Sparsification can be accomplished by using a sparse connectivity graph to construct low order operators. The sparsity of the graph adjacent matrix $\fnt{A}$ directly dictates the sparsity of the differentiation operators $\fnt{Q}^{L}$. With sparser low order operators, each quadrature node is connected to a fewer number of neighboring nodes. Sparsity reduces the amount of algebraic dissipation since the dissipation coefficient $d^{L,n}_{ij}$ is only non-zero for connected pairs of nodes in the graph. Details on the construction of sparse connectivity graphs are presented in Appendix \ref{subsec:sparse_con_graph} and numerical results from the low order scheme with different connectivity graphs are presented in Appendix B.
\end{subequations}
\section{Convex limiting}
\label{sec:4}
\begin{subequations}
\numberwithin{equation}{section}
In this section, we describe the application of convex limiting to combine the low order positivity preserving scheme with high order entropy stable discontinuous Galerkin  and apply it to the shallow water equations. This scheme is aligned with the ideas presented in \cite{guermond2019invariant, perthame1994variant, perthame1996positivity}.

In Sec. \ref{sec:3}, we showed that the low order solution satisfies the positivity preserving property. In particular, it preserves the positivity of water heights for the shallow water equations. However, it is low order accurate. To recover a high order accurate numerical scheme that also satisfies positivity requirements, we blend the high order and low order updates using ``convex limiting''. Convex limiting is only applied to elements where the high order scheme produces zero or negative water heights. 

\subsection{An overview of convex limiting}
We use the superscript $^H$ to denote high order terms. By assuming using forward Euler time stepping, we begin by rewriting high order discontinuous Galerkin  methods in the form
\begin{align}
    \frac{m_i}{dt}(\fnt{u}^{H,n+1}_i - \fnt{u}^{n}_i) + \sum_{j\in I(i)}\bm{F}^{H}_{ij}= m_i\fnt{S}(\fnt{u}_i^{n}),
\label{eq:high_order_base}
\end{align}
where $\bm{F}^{H}_{ij} \in \mathbb{R}$ contains the flux contribution between quadrature points $j$ and $i$. We also extend the results to higher order in time with strong-stability preserving Runge-Kutta (SSP-RK). Note that $\bm{F}^{H}_{ij}$ is skew-symmetric, $\bm{F}^{H}_{ij} = -\bm{F}^{H}_{ji}$. Subtracting \ref{eq:low_order_base} from \ref{eq:high_order_base}, we have the following:
\begin{align}
    m_i\fnt{u}^{H,n+1}_i = m_i\fnt{u}^{L,n+1}_i + \sum_{j\in I(i) \backslash \{i\}}dt(\bm{F}^{L}_{ij} - \bm{F}^{H}_{ij}).
\label{eq:high_order_base2}
\end{align}
Let $\bm{A}_{ij} = dt(\bm{F}^{L,n}_{ij} - \bm{F}^{H,n}_{ij})$. Then, we can rewrite (\ref{eq:high_order_base}) as
\begin{align}
    m_i\fnt{u}^{H,n+1}_i = m_i\fnt{u}^{L,n+1}_i + \sum_{j\in I(i) \backslash \{i\}}\bm{A}^n_{ij}.
\label{eq:high_order_update}
\end{align}
Note that we have $\bm{A}^n_{ij} = \bm{A}^n_{ji}$, which implies that $ \sum_{i\in V} m_i\fnt{u}^{L,n+1}_i =  \sum_{i\in V} m_i\fnt{u}^{H,n+1}_i$. Therefore, the high order and low order updates have the same mass. Our goal is to modify the high order update so that the modified high order update satisfies the positive water heights constraint. We introduce symmetric limiting parameters $l_{ij} = l_{ji} \in [0,1]$,  for all $i,j \in V$. We define the limited update $\fnt{u}^{n+1}_i$ as
\begin{align}
    m_i\fnt{u}^{n+1}_i = m_i\fnt{u}^{L,n+1}_i + \sum_{j\in I(i) \backslash \{i\}}l_{ij} \bm{A}^n_{ij}.
\label{eq:convex_limit_general}
\end{align}
The limiting process is conservative for any choice of symmetric coefficients $l_{ij} = l_{ji}$ \cite{guermond2019invariant}. Notice that $\fnt{u}^{n+1}_i  = \fnt{u}^{L,n+1}_i$ if  $l_{ij} = 0$ for all $j \in I(i)\backslash{i}$ and $\fnt{u}^{n+1}_i  = \fnt{u}^{H,n+1}_i$ if  $l_{ij} = 1$ for all $j \in I(i)\backslash{i}$. Therefore, $h>0$ when all $l_{ij} = 0$. To obtain a more accurate approximation, we want the limiting coefficients  $l_{ij}$ to be as close to 1 as possible while enforcing $h>0$. 

We can rewrite equation \ref{eq:convex_limit_general} as
\begin{align}
    m_i\fnt{u}^{n+1}_i = \sum_{j\in I(i) \backslash \{i\}}l_{ij} \lambda_j (\fnt{u}^{L,n+1}_i + l_{i,j}\bm{P}^n_{ij}), \quad \bm{P}^n_{ij} = \frac{1}{m_i\lambda_j}\bm{A}^n_{ij}, 
\label{eq:convex_limit_2}
\end{align}
where $\{\lambda_j\}_{j\in I(i) \backslash \{i\}}$ is any set of strictly positive convex coefficients such that $\sum_{j\in I(i) \backslash \{i\}}\lambda_j = 1$, $\lambda_j>0$ for all $j\in I(i) \backslash \{i\}$. We choose $\lambda_j = \frac{1}{{\rm card}(I(i)) -1}$ as in \cite{guermond2019invariant}. To calculate $l_{ij}$, we first define $l^i_j$ for every $i\in V$ and $j \in I(i)$:
\begin{align}
    l^i_j = \begin{cases}
    1 & {\rm if} \; h^{L,n+1}_i+P^{n,h}_{ij} \geq 0,\\
    {\rm max}\{l \in [0,1] | h^{L,n+1}_i+lP^{n,h}_{ij}\} \geq 0  &{\rm otherwise}, 
    \end{cases}
    \label{eq:lij_algo}
\end{align}
where $P^{n,h}_{ij}$ corresponds to the update on the $h$ component of the solution. We set the node-wise limiting factor $l_{ij} = l_{ji} = \min\{l^i_j, l^j_i\}$. With this definition, we can ensure that the water heights stay positive by Theorem 7.21 in \cite{guermond2019invariant}. While most of time, finding $l^i_j$ require a line search, it remains simple for the shallow water equations since we can express $l^i_j$ explicitly, due to the fact that $h$ is a linear function, as follows:
\begin{align}
    l^i_j = \begin{cases}
    1 & {\rm if} \; h^{L,n+1}_i+P^{n,h}_{ij} \geq 0,\\
    -\frac{h^{L,n+1}}{P^{n,h}_{ij}}  &{\rm otherwise}.
    \end{cases}
    \label{eq:lij_algo_exp}
\end{align}
With convex limiting, we can ensure that our solution remains positive in water heights throughout a simulation.

An alternative limiting scheme, called element-wise limiting, is also considered in this paper. In this scheme, we have a uniform limiting factor across the same element. We define the element-wise limiting factor $l_{e}  = \min_{i,j}\{l_{ij}\}$. This scheme produces solutions that are closer to the low order positivity preserving scheme than the node-wise limiting scheme because it always has smaller local limiting factors. In practice, solutions from the element-wise limiting scheme are usually bounded further away from 0 and therefore, produce smaller maximum wave speed on the shallow water equations. This has added benefit of allowing us to take larger time steps. We can also prove that, in contrast to the node-wise limiting scheme, the element-wise limiting scheme also preserves the entropy stability of the blended scheme \cite{hennemann2021provably}. 

\subsection{Local implementation}
One important and efficient property of our numerical scheme is that it only requires computing limiting factors $l_{ij}$ for the nodes $i$, $j$ within the same element, which was first observed in \cite{pazner2021sparse}. We use the same entropy stable interface numerical fluxes with a Lax-Friedrichs dissipation in both high order and low order schemes. The resulting fluxes on the boundaries of each element are the same. For example, consider the high order and low order flux matrices for a simple periodic domain in 1D with two elements shown in Fig. \ref{fig:cl_demo1}. Suppose there are three nodes in each element. The non-zero pattern of the matrix $\bm{A} = dt(\bm{F}^{L} - \bm{F}^{H})$ from Eq. (\ref{eq:high_order_update}) is shown in Fig. \ref{fig:cl_demo2}. We have non-zero entries from the off diagonal block due to the fluxes calculated on the boundaries of each element, namely the terms $a_{16}$, $a_{34}$, $a_{43}$, and $a_{61}$. However, these boundary fluxes are the same for both high order and low order updates; therefore corresponding entries in the difference matrix $\bm{A}$ vanish.
\begin{figure}
\centering
$\bm{F}^H = \begin{bmatrix}
\begin{array}{ccc|ccc}
f^H_{11} & f^H_{12} &  f^H_{13} & && f^H_{16} \\
f^H_{21} & f^H_{22} &  f^H_{23} & &&  \\
f^H_{31} & f^H_{32} &  f^H_{33} &  f^H_{34} && \\
\hline
 &&f^H_{43} & f^H_{44}&  f^H_{45} &  f^H_{46}\\
 &&& f^H_{54} &  f^H_{55} &  f^H_{56}\\
 f^H_{61} &&& f^H_{64} & f^H_{65} &  f^H_{66}\\
\end{array}
\end{bmatrix}$ \hspace{0.5cm}
$\bm{F}^L = \begin{bmatrix}
\begin{array}{ccc|ccc}
f^L_{11} & f^L_{12} &  & && f^L_{16} \\
f^L_{21} &  &  f^L_{23} & &&  \\
 & f^L_{32} &  f^L_{33} &  f^L_{34} && \\
\hline
 &&f^L_{43} & f^L_{44}&  f^L_{45} &  \\
 &&& f^L_{54} &   &  f^L_{56}\\
 f^L_{61} &&& & f^L_{65} &  f^L_{66}\\
\end{array}
\end{bmatrix}$
\caption{Illustration of high order and low order flux matrices for two one-dimensional elements.}
\label{fig:cl_demo1}
\end{figure}

\begin{figure}
\centering
\[
\bm{A} = \begin{bmatrix}
\begin{array}{ccc|ccc}
a_{11} & a_{12} &  a_{13} & && a_{16} \\
a_{21} & a_{22} &  a_{23} & &&  \\
a_{31} & a_{32} &  a_{33} &  a_{34} && \\
\hline
 &&a_{43} & a_{44}&  a_{45} &  a_{46}\\
 &&& a_{54} &  a_{55} &  a_{56}\\
 a_{61} &&& a_{64} & a_{65} &  a_{66}\\
\end{array}
\end{bmatrix}
= \begin{bmatrix}
\begin{array}{ccc|ccc}
a_{11} & a_{12} & a_{13} &&& \\
a_{21} & a_{22} & a_{23} &&& \\
a_{31} & a_{32} & a_{33} &&& \\
\hline
 &&& a_{44} & a_{45} & a_{46}\\
 &&& a_{54} & a_{55} & a_{56}\\
 &&& a_{64} & a_{65} & a_{66}\\
\end{array}
\end{bmatrix}
\]
\caption{Illustration of the flux correction matrix $\bm{A} =  dt(\bm{F}^L - \bm{F}^H)$ for two one-dimensional elements.}
\label{fig:cl_demo2}
\end{figure}
As a result, we do not need to implement the convex limiting procedure for the pairs of boundary nodes on neighboring elements. Convex limiting is only required for the pairs of nodes within the same element. With this key property, the implementation of convex limiting is purely local, simplifying the resulting numerical scheme.
\end{subequations}
\section{Numerical results}
\label{sec:5}
\begin{subequations}
\numberwithin{equation}{section}
In this section, we present some numerical experiments and results in both 1D and 2D to demonstrate the accuracy and stability of the convex limited entropy stable DG scheme. We also add a local Lax-Friedrichs penalization \cite{wintermeyer2017entropy} to all interface fluxes.

The first experiment is a ``lake-at-rest" condition to test the well-balancedness of our scheme in 1D. The second experiment is the parabolic bowl problem in 1D, designed to investigate the convergence rate when positivity preservation is activated. The third experiment is a dam break simulation in 2D from \cite{wintermeyer2018entropy}. The last experiment is a wave over a bump in 2D \cite{xing2013positivity}. 

All numerical experiments utilize a second order strong-stability preserving Runge-Kutta (SSP-RK) scheme (e.g., Heun's method) \cite{suli2003introduction} for time integration. Following the derivation of stable timestep restrictions in \cite{chan2016gpu}, we define the timestep $dt^H$ for the high-order scheme to be
\begin{align}
dt^H = CFL \times \frac{h}{C_N},\hspace{1cm} C_N = \frac{(N+1)(N+2)}{2},
\end{align}
where $C_N$ is the degree dependent constant in the trace inequality space \cite{warburton2003constants}, and CFL is a user-defined constant. We use CFL = 0.125 for all experiments. For the low order scheme, the timestep $dt^L$ is calculated via \begin{align}
dt^L = \min_{i\in V} \frac{w_i}{2\lambda_{\max}},
\end{align}
where $\lambda_{\max}$ is the maximum wave speed calculated with Eq. (\ref{eq:shallow_water_wave_speed}) over the entire domain. We define the timestep $dt = {\min}\{ dt^H, dt^L \}$. To avoid division by zero, we enforce the water height to be above a small positive tolerance $tol = 10^{-14}$ in dry areas. 

\subsection{Lake at rest}
The test problem is chosen to verify that the convex limited DG scheme indeed preserves the still water steady state with a non-flat bottom containing a wet/dry interface. The bottom topography is given by the depth function \cite{xing2010positivity, liang2009numerical}:
\begin{align}
    b(x) = {\max} (0, -20(x-1/8)(x+1/8)+2) \quad {\rm on} \; -1\leq x \leq 1.
\end{align}
The initial data is the stationary solution:
\begin{align}
    h(x)+b(x)  = {\max}(2,b),\quad hu = 0, \quad {\rm on} \; -1\leq x \leq 1.
\end{align}
A periodic boundary condition is used in this experiment. We compute the solution until T = 1 using 128 uniform cells. The computed surface level $h + b$ and the bottom $b$ are plotted in Fig. \ref{fig:CL_WB}. In order to demonstrate that the lake at rest solution is preserved up to round-off error, we compute the $L^2$ error of the solution. In double precision, the $L^2$ error is 3.56E-15. 

We note that well-balancedness is not preserved if the mesh element boundary is not aligned with the wet/dry interface, and that designing a sub-cell well-balanced scheme remains challenging. For example, \cite{azerad2017well} presents two numerical schemes, and the scheme which preserves sub-cell well-balanceness is restricted to second order accuracy. 
\begin{figure}
\begin{center}
\includegraphics[width=.65\textwidth]{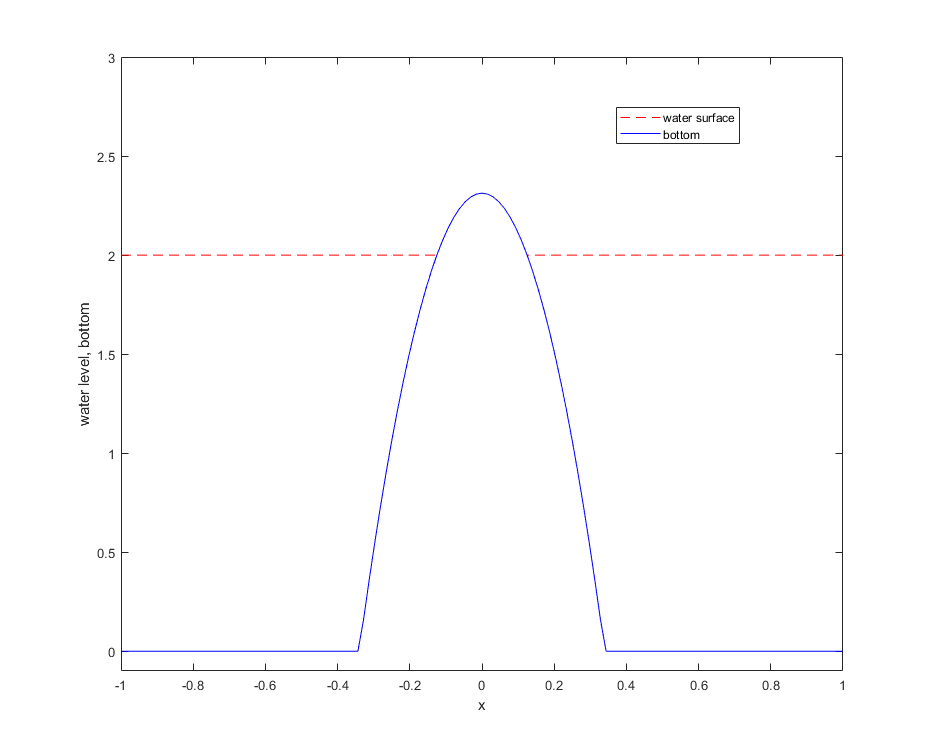}
\caption{The surface level h + b and the bottom b for the lake at rest}
\label{fig:CL_WB}
\end{center}
\end{figure}

\subsection{Translating vortex}
We now consider the translating vortex test. We set the domain to be $[-10,10]\times[-5,5]$. The exact solution for the vortex at any time $t$ is given by \cite{kang2020imex, ricchiuto2009stabilized}:
\begin{align}
&h = h_{\infty}-\frac{\beta^2}{32\pi^2}e^{-2(r^2-1)},\\
&u = u_{\infty}-\frac{\beta}{2\pi}e^{-2(r^2-1)}y_t,\\
&v = v_{\infty}+\frac{\beta}{2\pi}e^{-2(r^2-1)}x_t,\\
&b = 0,
\end{align}
where
\begin{align}
x_t &= x - x_c - u_{\infty}t,\\
y_t &= y - y_c - v_{\infty}t,\\
r^2 &= x^2 + y^2.
\end{align}
In this example, we set 
\begin{align}
h_{\infty} = 1,\quad \beta = 5, \quad  g = 2 \quad {\rm{and}}\quad (u_{\infty}, v_{\infty}) = (1,0).
\end{align}
Initially the vortex is located at $(x_c, y_c) = (0, 0)$. In this setup, the vortex propagates to the right along the x-axis.  The domain size is chosen to be large enough such that periodic boundary conditions can be used without affecting accuracy.

We use affine triangular meshes for this experiment. We compute the $L^2$ error for the SBP formulation using operators based on Gauss-Legendre quadrature nodes. Recall that the SBP-DG discretization does not correspond to a polynomial approximation space.  Thus, to calculate  $L^2$ errors for the SBP-DG discretization, we first project the final numerical solution to polynomials of degree $N$. We then use this projection to evaluate the $L^2$ error using a quadrature rule which is exact for at least degree $2N+2$ polynomials. The convergence results are presented in Fig. \ref{fig:vortex_error}. We observe an $O(h^{N+1})$ rate of convergence for the high order method on affine triangular meshes. The method can also be extended to curved triangular meshes \cite{wu2021high}.
\begin{figure}
\begin{center}
\includegraphics[width=.6\textwidth]{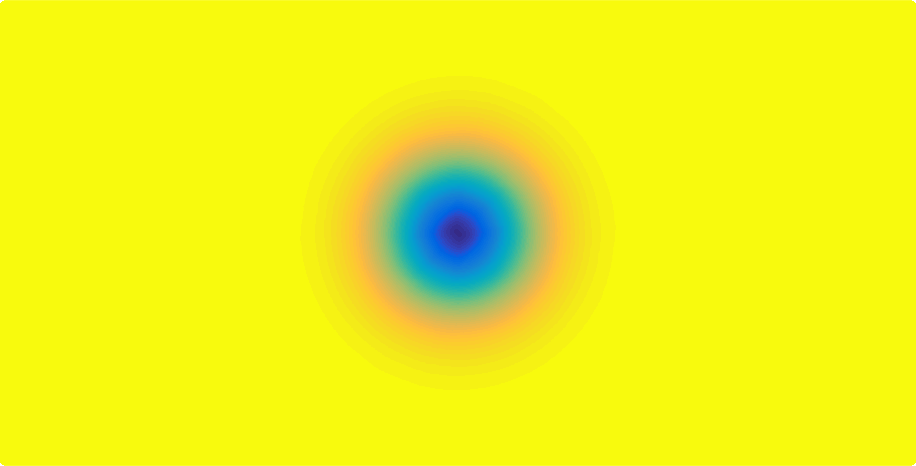}
\caption{A translating vortex in 2D}
\label{fig:vortex}
\end{center}
\end{figure}
\begin{figure}
\begin{center}
\includegraphics[width=.6\textwidth]{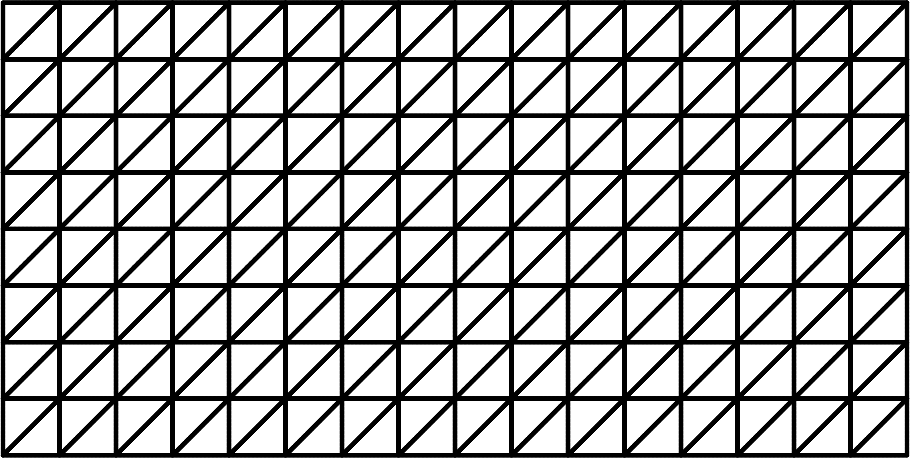}
\caption{Meshes in 2D for translating vortex problem}
\label{fig:mesh vortex}
\end{center}
\end{figure}

\begin{figure}
\begin{center}
\includegraphics[width=0.8\textwidth]{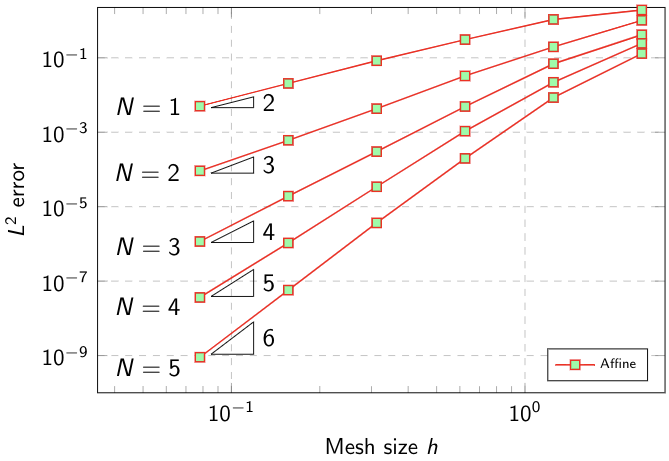}
\caption{$L^2$ errors for the translating vortex after 0.5 seconds using nodal DG and Gauss-Legendre SBP operator schemes on various mesh sizes and polynomial degrees.}
\label{fig:vortex_error}
\end{center}
\end{figure}

\subsection{Parabolic bowl}
The second numerical experiment is used to verify the accuracy and convergence rate of the convex limited DG scheme. For the one-dimensional shallow water equations with a parabolic bottom topography, analytic solutions have been derived by Sampson
et al. \cite{sampson2005moving}. For this example, the parabolic bottom is given by
\begin{align}
    b(x) = h_0(x/a)^2,
\end{align}
where we set $h_0 = 8$ and $a = 2$. The computational domain is $[-5,5]$. The analytical water height for the shallow water equations is given by
\begin{align}
    h(x,t) + b(x) = h_0 - \frac{B^2}{4g}\cos(2\omega t) - \frac{B^2}{4g} - \frac{Bx}{2a} \sqrt{\frac{8h_0}{g}}\cos(\omega,t),
    \label{eq:p_bowl}
\end{align}
where $\omega = \sqrt{2gh_0}/a$ and $B=2$. The exact location of
the wet/dry fronts take the form
\begin{align}
   x_0 = -\frac{B\omega a^2}{2gh_0}\cos(\omega t) \pm a
\end{align}
The initial water height is then defined by (\ref{eq:p_bowl}). The velocity is initialized to zero. We can use either periodic or wall boundary conditions for this experiment because the water will never reach the boundary. We set the polynomial degree to $N = 3$ and run the simulation until $T = 1$ with various numbers of uniform elements. We plot the water surface at different times using the node-wise limiting and the element-wise limiting in Fig. \ref{fig:CL_PB_node} and Fig. \ref{fig:CL_PB_elem}, respectively. We also include the analytical solution to provide a comparison. We observe that the numerical solution matches well with the analytic solution. The element-wise limiting scheme produces smoother solutions than the node-wise limiting scheme. The node-wise limiting scheme produces some spurious oscillations near the dry front, though these oscillations are localized and vanish as we refine the mesh.

We plot the $L^2$ norm of the errors for different polynomial degrees and mesh sizes for the node-wise and element-wise limiting schemes in Fig. \ref{fig:CL_Berror_node} and Fig. \ref{fig:CL_Berror_elem}. We observe both schemes have a convergence rate of between $O(h)$ and $O(h^2)$. The node-wise limiting produces larger errors near the dry-wet interface, where the convex limiting scheme pushes towards the $O(h)$ low order scheme. The element-wise limiting, on the other hand, produces smaller errors throughout the entire wet domain for lower polynomial degree $N$. Element-wise limiting errors also converge slightly faster than node-wise limiting errors. 

\begin{figure}
\begin{center}
\includegraphics[width=0.48\textwidth]{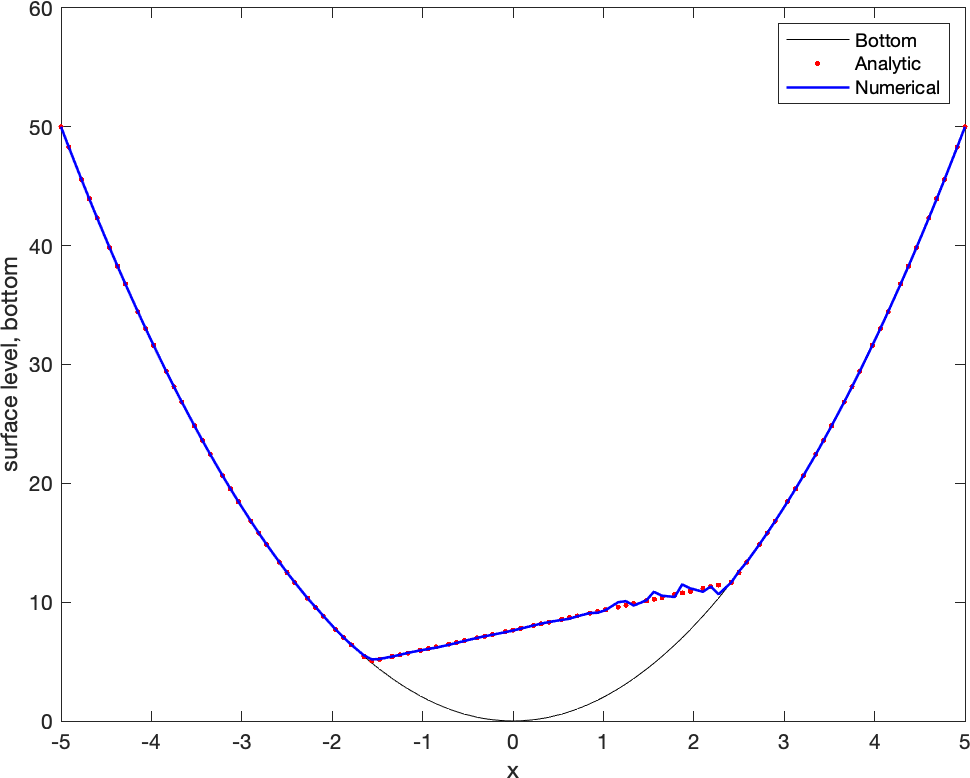}
\includegraphics[width=0.48\textwidth]{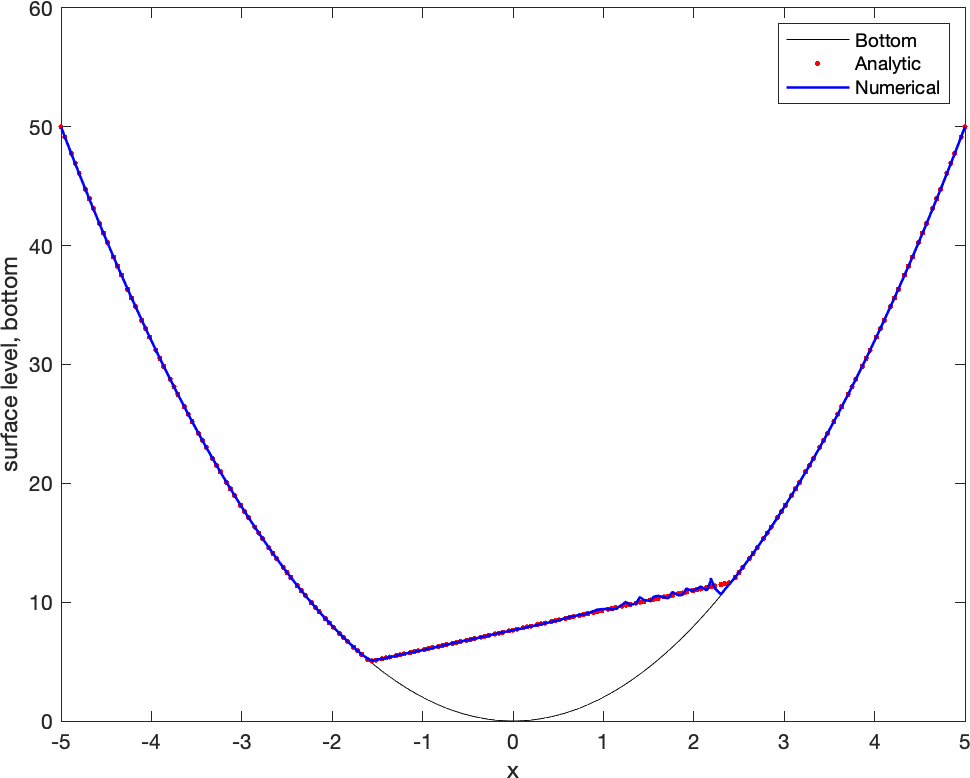}\\

\includegraphics[width=0.48\textwidth]{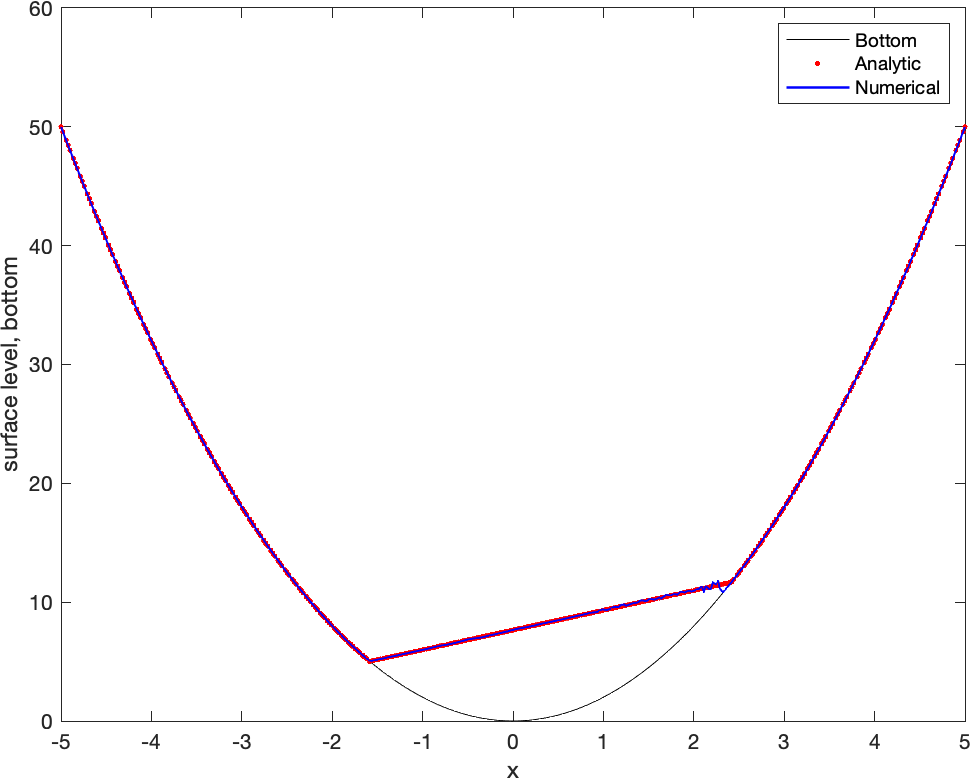}
\includegraphics[width=0.48\textwidth]{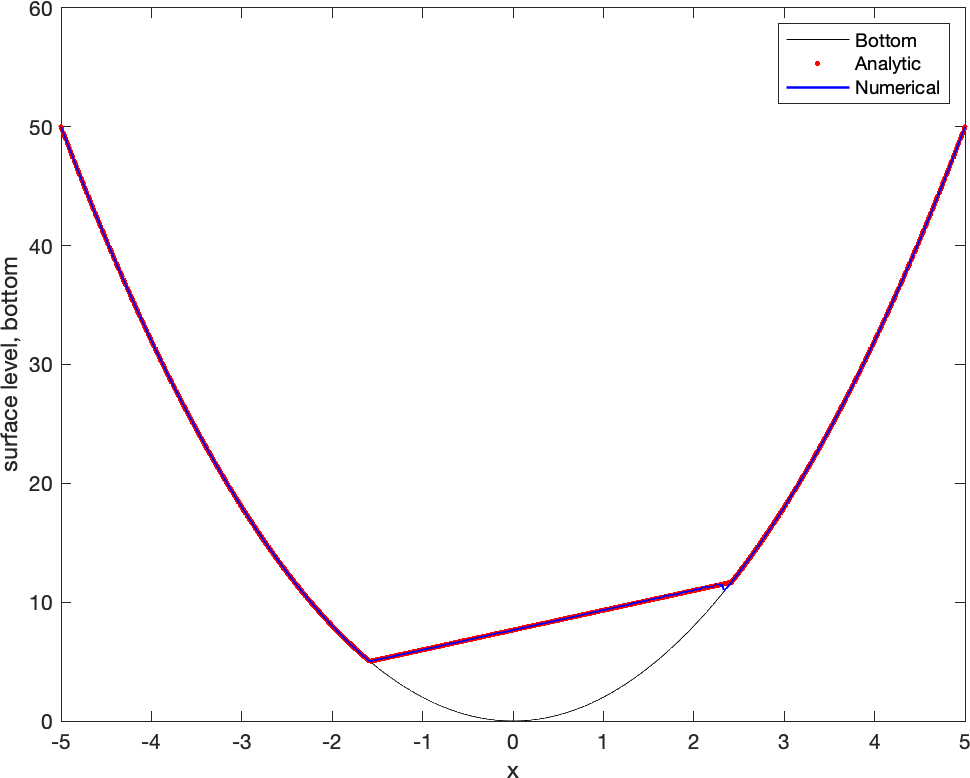}

\caption{The surface level h + b and the bottom b for the parabolic bowl at $T=1$ with node-wise convex limiting: $N=3$; top left: $K = 32$; top right: $K = 64$; bottom left: $K = 128$; bottom right: $K = 256$}
\label{fig:CL_PB_node}
\end{center}
\end{figure}

\begin{figure}
\begin{center}
\includegraphics[width=0.48\textwidth]{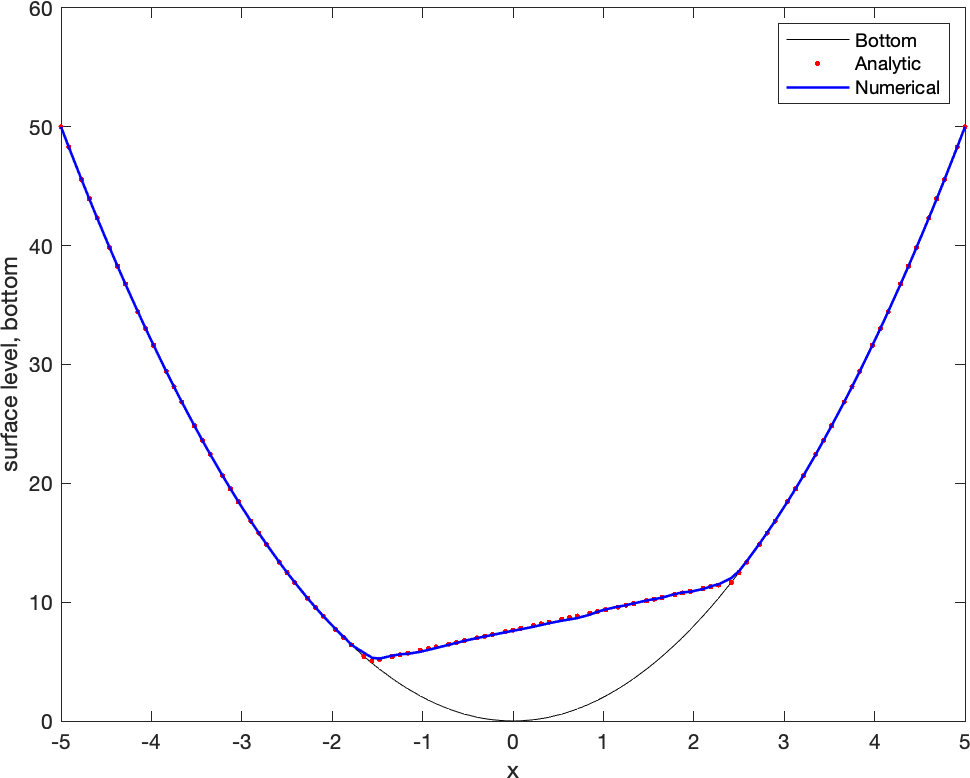}
\includegraphics[width=0.48\textwidth]{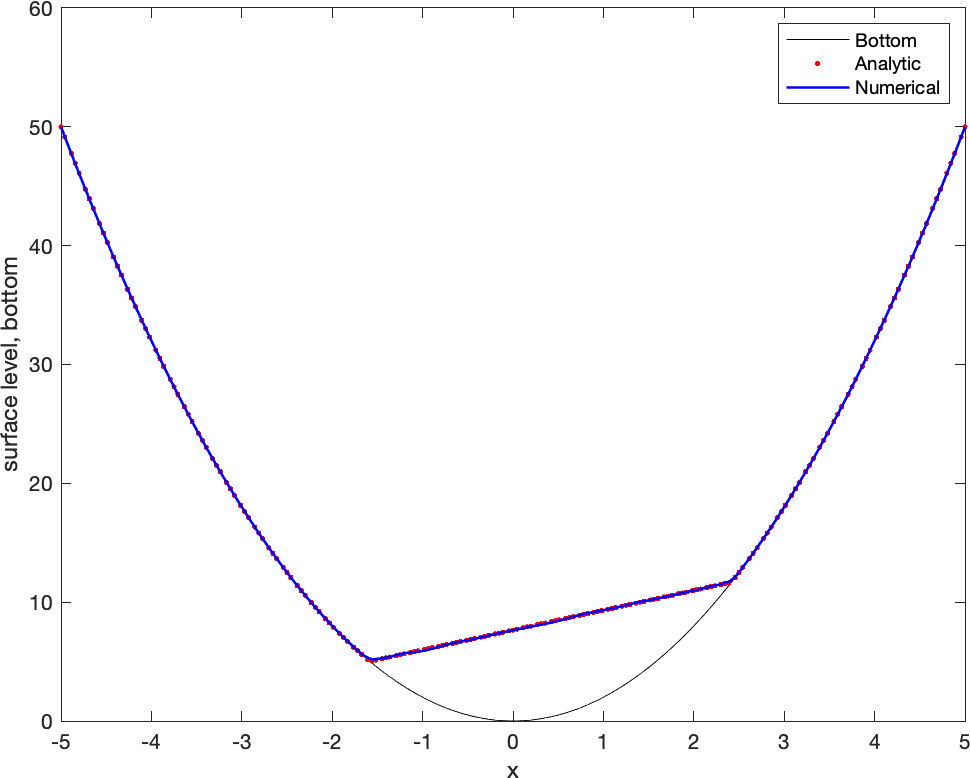}\\

\includegraphics[width=0.48\textwidth]{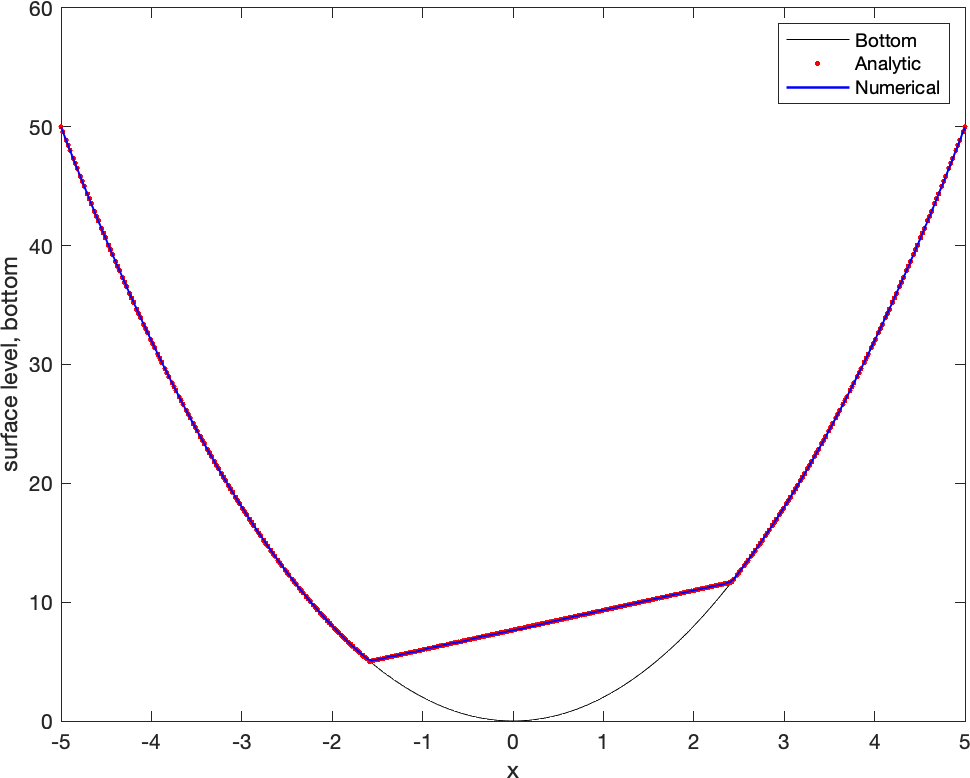}
\includegraphics[width=0.48\textwidth]{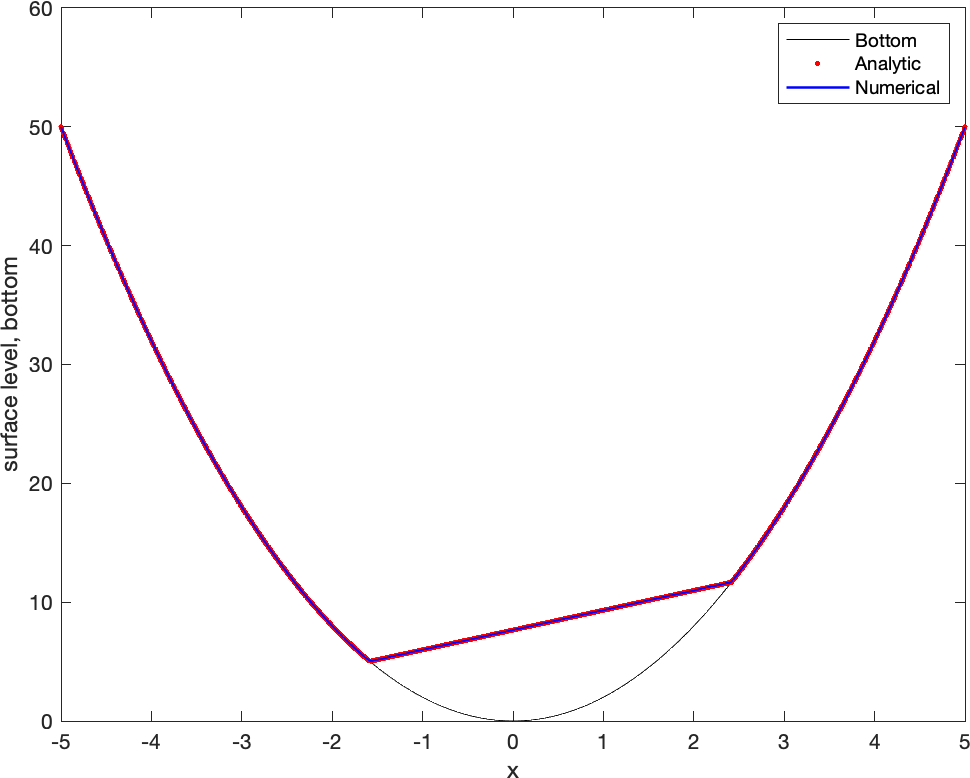}

\caption{The surface level h + b and the bottom b for the parabolic bowl at $T=1$: $N=3$ with element-wise convex limiting; top left: $K = 32$; top right: $K = 64$; bottom left: $K = 128$; bottom right: $K = 256$}
\label{fig:CL_PB_elem}
\end{center}
\end{figure}

\begin{figure}
\begin{center}
\includegraphics[width=0.8\textwidth]{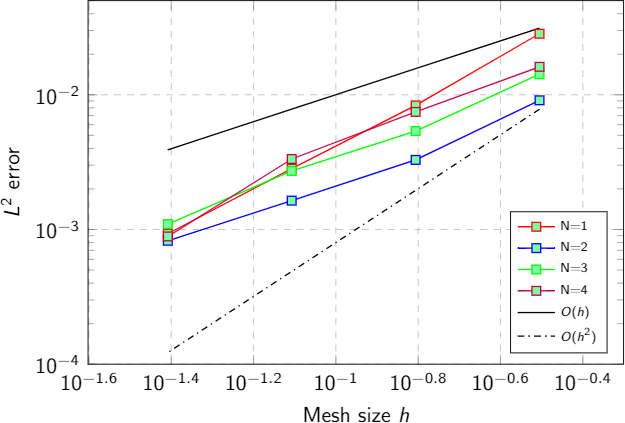}

\caption{$L^2$ errors for the 1D oscillating bowl at T=0.5 using node-wise convex limiting with various mesh sizes and polynomial degrees}
\label{fig:CL_Berror_node}
\end{center}
\end{figure}

\begin{figure}
\begin{center}
\includegraphics[width=0.8\textwidth]{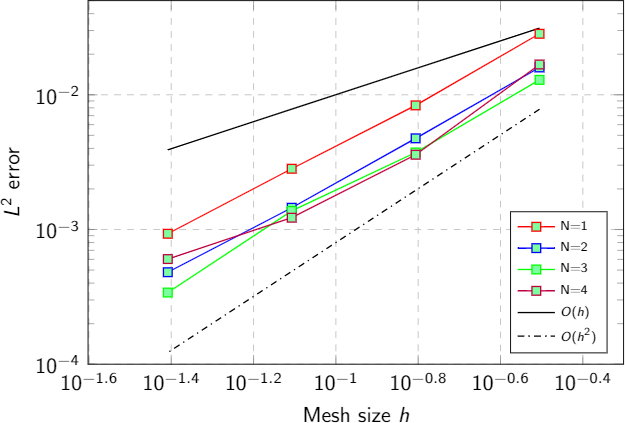}

\caption{$L^2$ errors for the 1D oscillating bowl at T=0.5 using element-wise convex limiting with various mesh sizes and polynomial degrees}
\label{fig:CL_Berror_elem}
\end{center}
\end{figure}

\subsection{Dam break}
\label{sec:5_Dam_break}
This experiment is taken from \cite{wintermeyer2018entropy, lukacova2009entropy}. Similar dam break experiments and results can be found in \cite{wintermeyer2017entropy, fjordholm2011well}. We utilize the same physical setting but use triangular meshes instead of quadrilateral meshes. The domain $[-1,2]\times[-1,1]$ is discretized using a $48\times32$ grid of quadrilaterals, and each of them is split into two triangular elements. We set $N=3$ for the polynomial degree. The dam is modeled by imposing reflective boundary conditions along the y-axis with a break between $y=-0.1$ and $y=0.1$ to allow water to flow through, as shown in Fig. \ref{fig:Dam0}.

We start with a zero velocity and constant water heights on both sides of the dam. The bathymetry is set to $b = 0$ on both sides. We set the initial water height to be $h = 5$ on the left side of the dam and $h=tol$ on the right side as shown in Fig. \ref{fig:Dam0}.
\begin{figure}
\begin{center}
\includegraphics[width=0.48\textwidth]{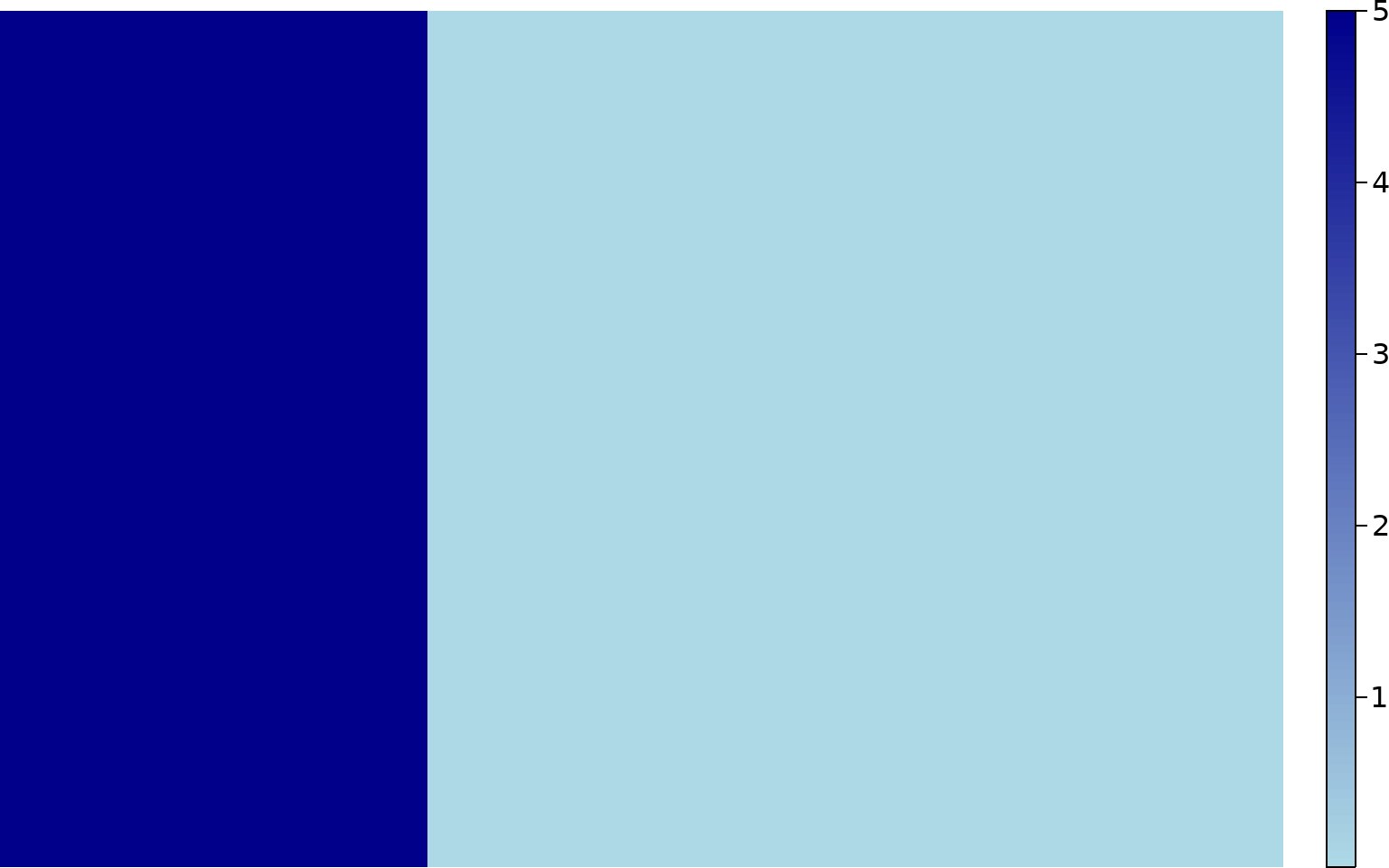}
\includegraphics[width=0.48\textwidth, height = 0.25\textheight]{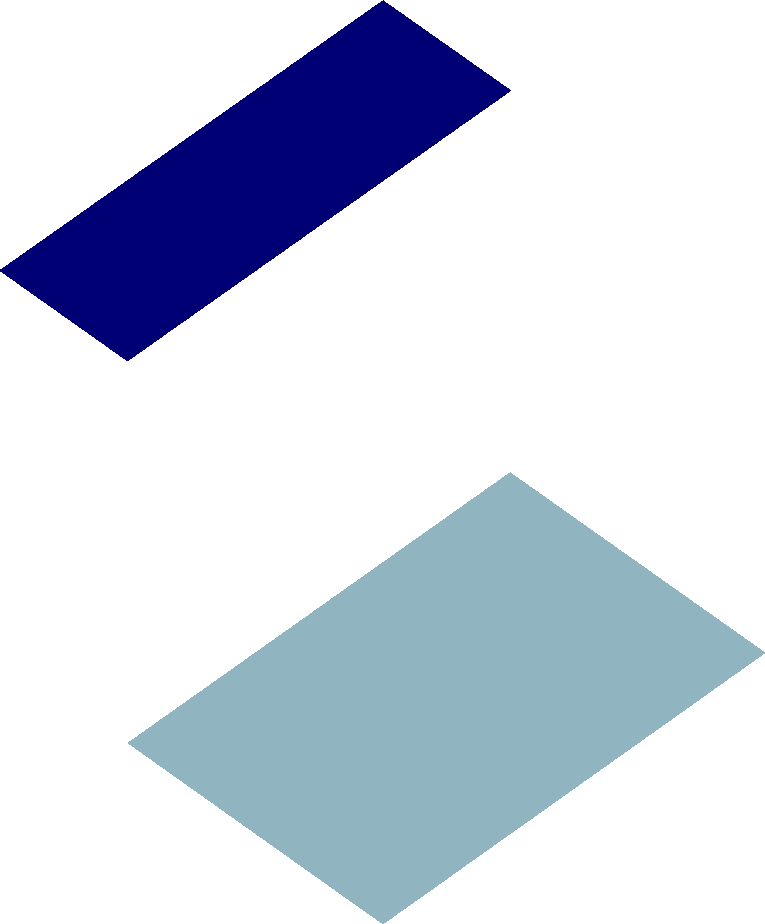}\\
\caption{2D dam break problem with initial water height}
\label{fig:Dam0}
\end{center}
\end{figure}
We plot the solutions from node-wise limiting, element-wise limiting, and the low order positivity preserving method at several different times in Fig. \ref{fig:Dam_n1} through \ref{fig:Dam_n4}. We observe that the water falling from the left side of the dam to the right produces a wavefront in the lower half of the dam. This wave is discontinuous, so we observe some mesh dependent solution oscillations, but they are localized and do not cause the solution to blow up. The water heights remain positive during the entire numerical experiment. We have also tried simulations with polynomial degrees $N= 1$ and $N=2$. Both of these simulations also remain stable throughout the experiment. The numerical solution of the parabolic dam break problem demonstrates that the entropy stable numerical method remains robust in the presence of shock discontinuities, while the use of convex limiting preserves positive water heights.

We compare the solutions from all three methods. We observe that the results from the node-wise limiting and the element-wise limiting schemes are very similar. They both become wavy near the right boundary. The results from the low order positivity preserving method are much more dissipative than the results from the convex limited schemes. The node-wise limiting scheme produces slightly larger oscillations; however, these oscillations may not be physical as they are not observed in other reference solutions in the literature \cite{tadmor2008energy, wintermeyer2018entropy}. The numerical results from the element-wise limited scheme are more consistent with reference solutions from the literature. This may be due to the fact that an element-wise constant limiting factor preserves a semi-discrete entropy inequality \cite{hennemann2021provably}.


\begin{figure}
\begin{center}
\includegraphics[width=0.32\textwidth, height = 0.2\textheight,
trim={0 0 5cm 0},clip]{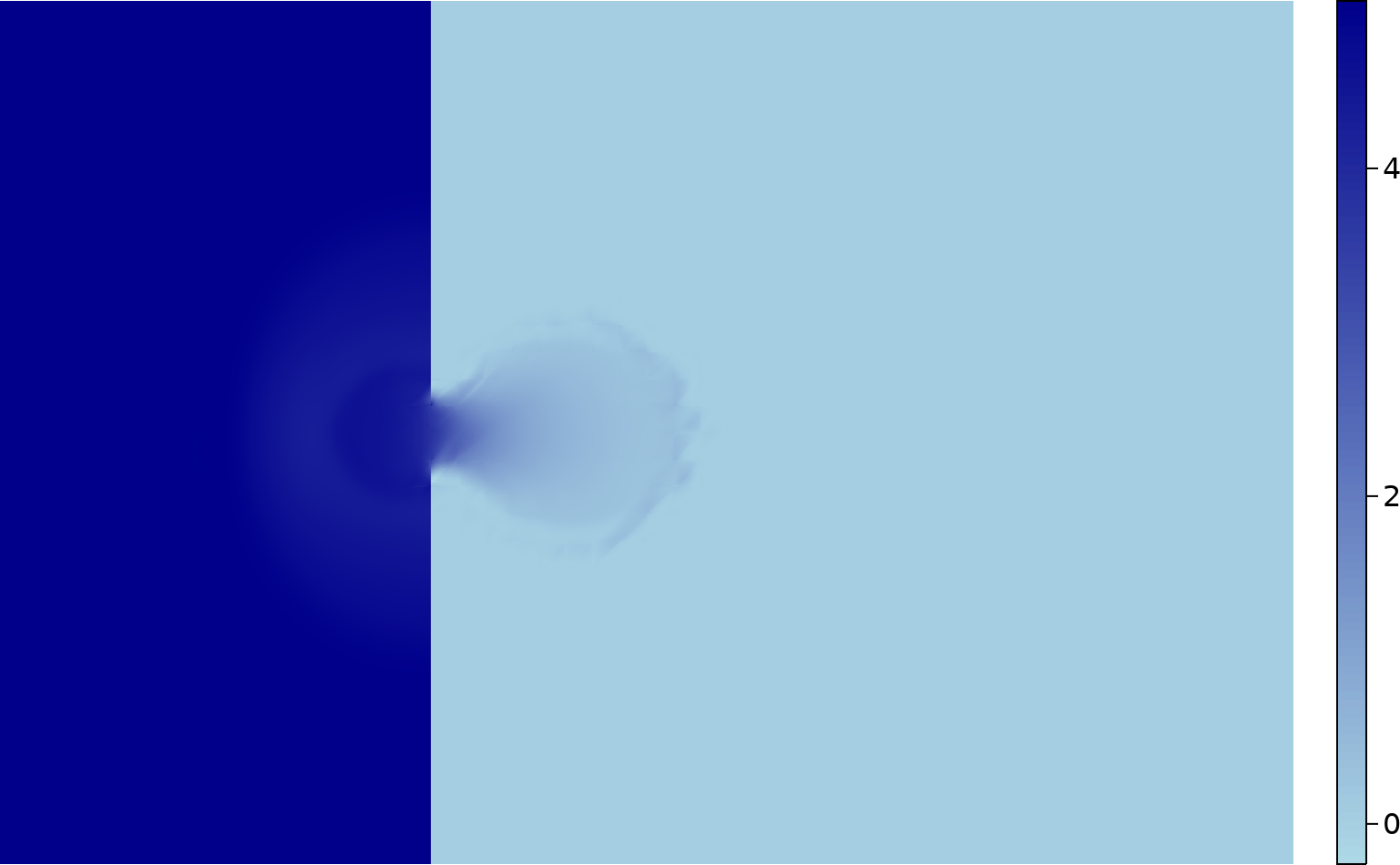}
\includegraphics[width=0.32\textwidth, height = 0.2\textheight, 
trim={0 0 5cm 0},clip]{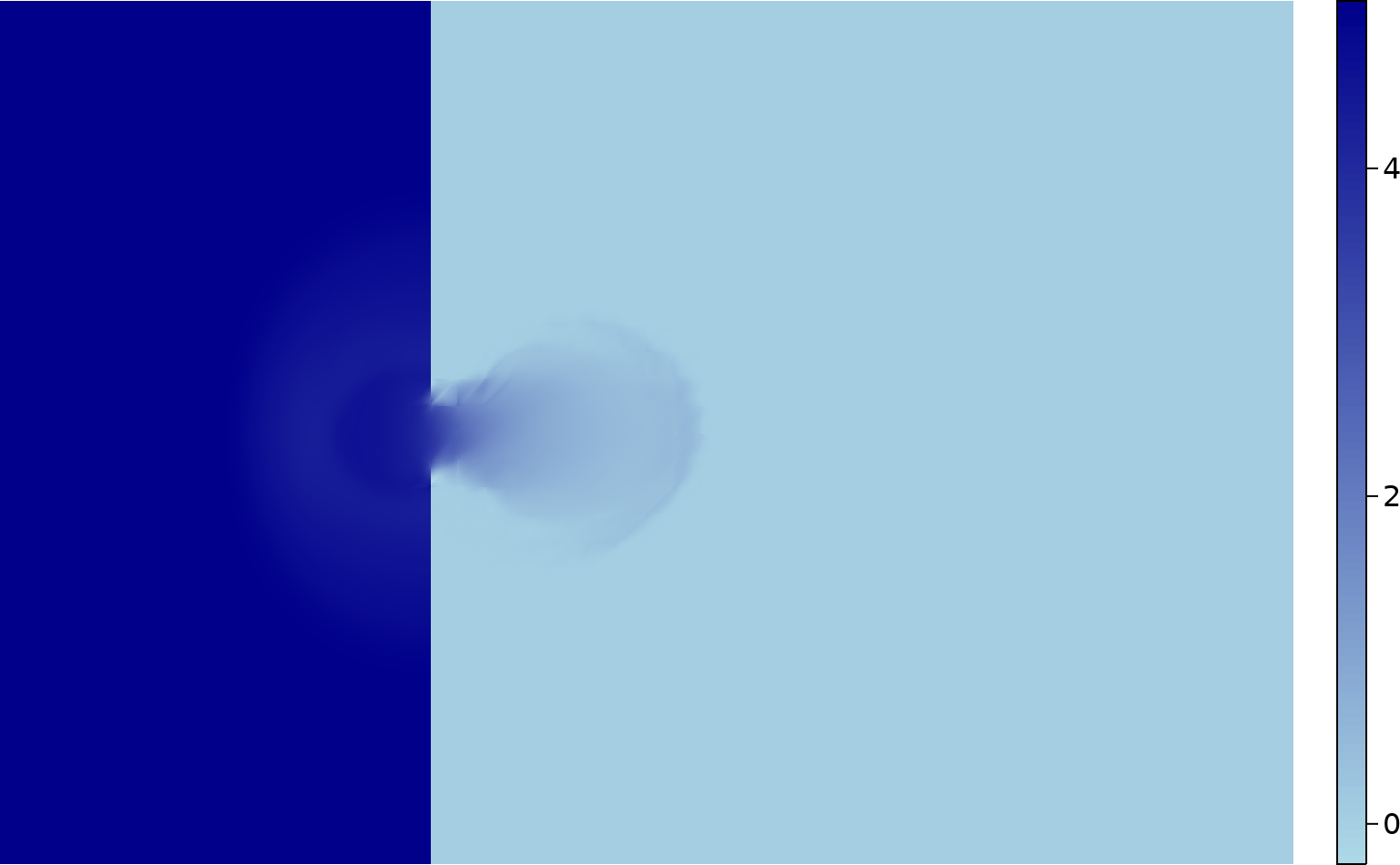}
\includegraphics[width=0.32\textwidth, height = 0.2\textheight, 
trim={0 0 5cm 0.3cm},clip]{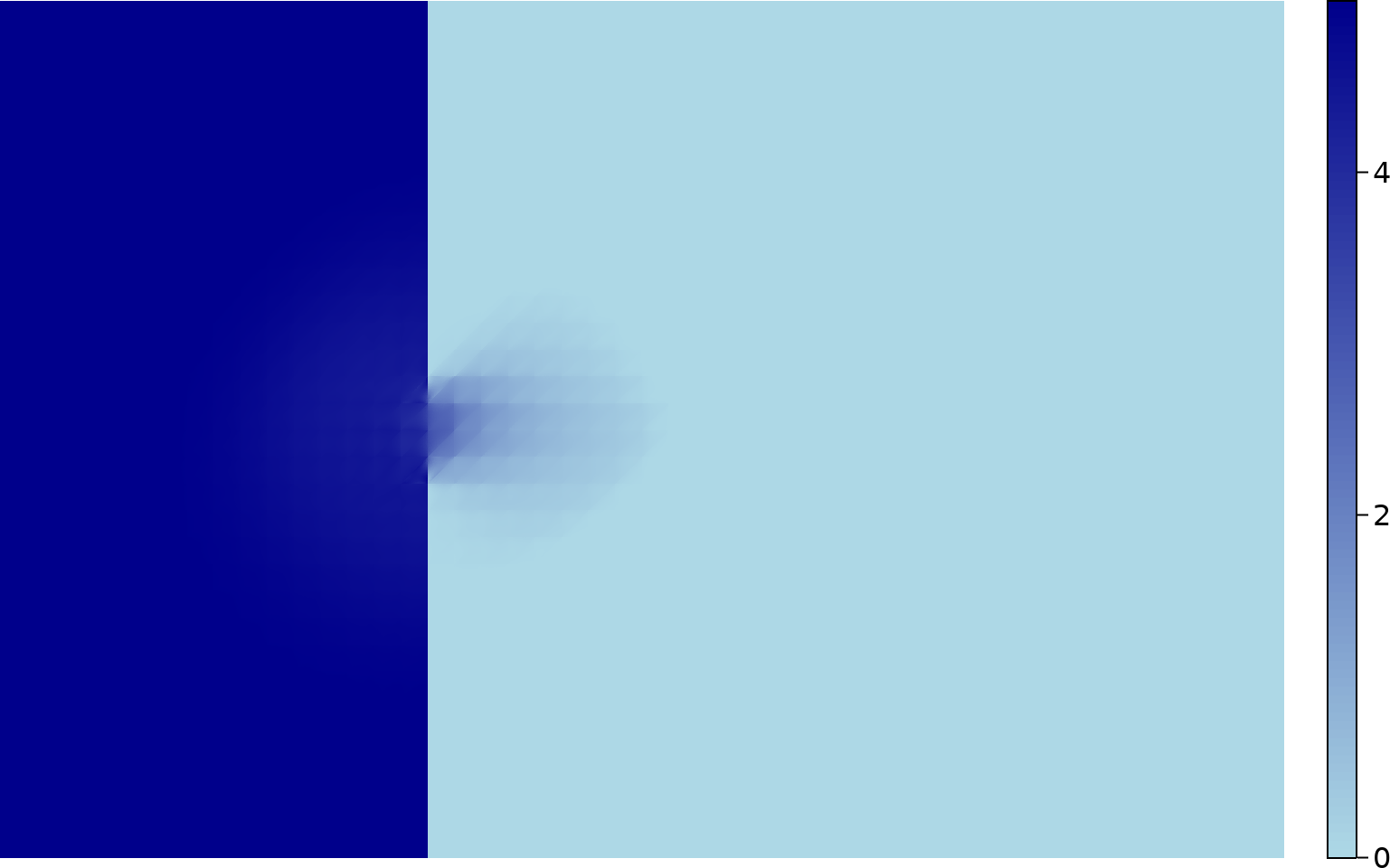}\\

\includegraphics[width=0.32\textwidth, height = 0.21\textheight]{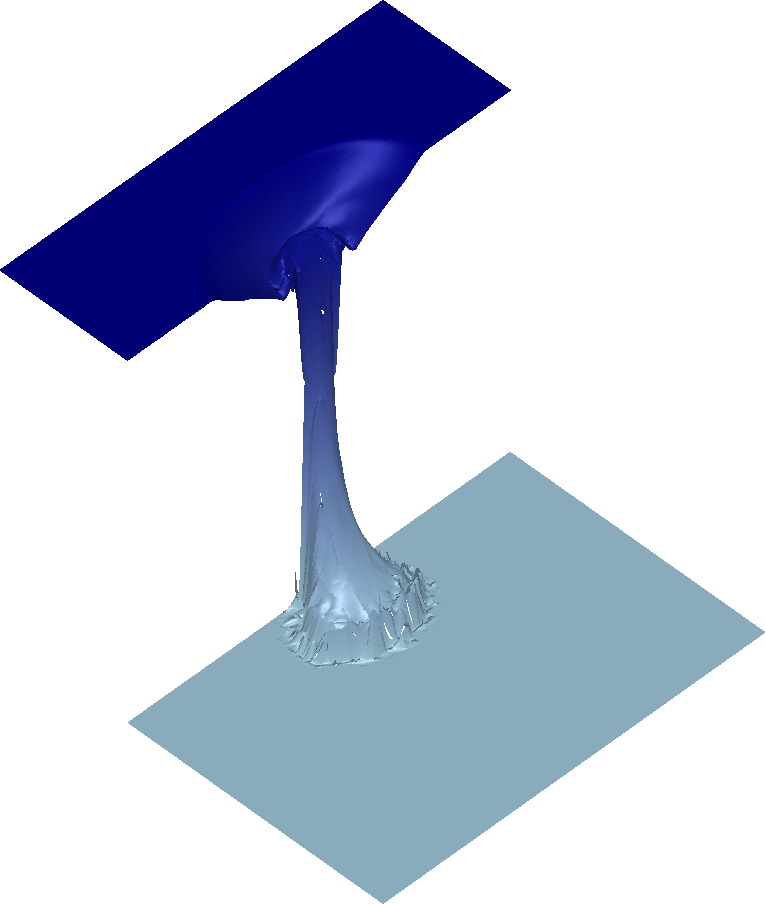}
\includegraphics[width=0.32\textwidth, height = 0.21\textheight]{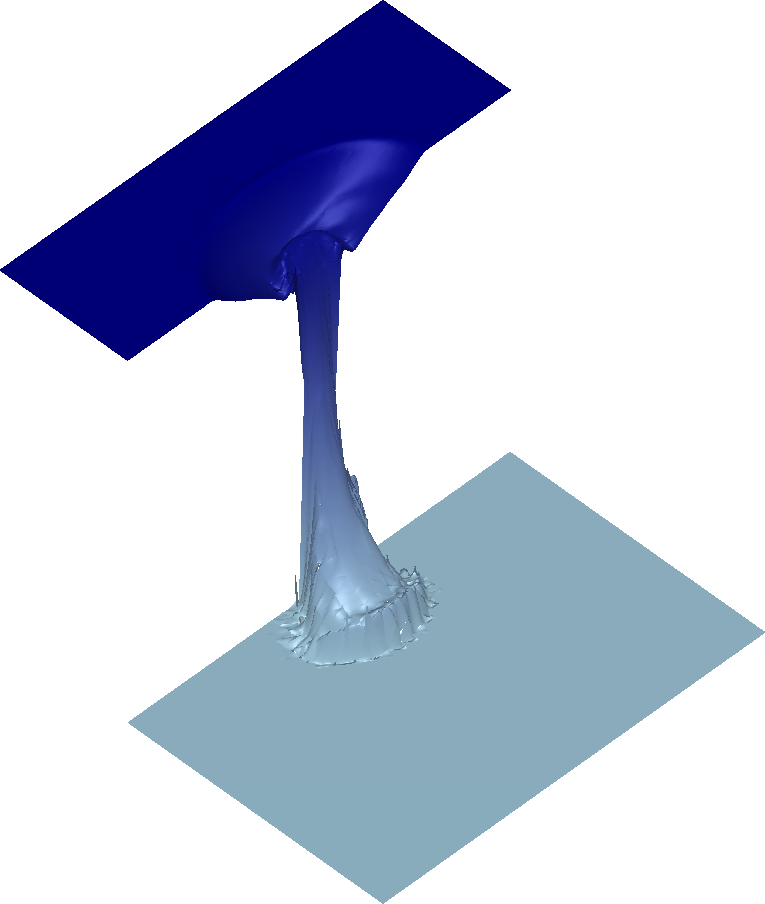}
\includegraphics[width=0.32\textwidth, height = 0.21\textheight]{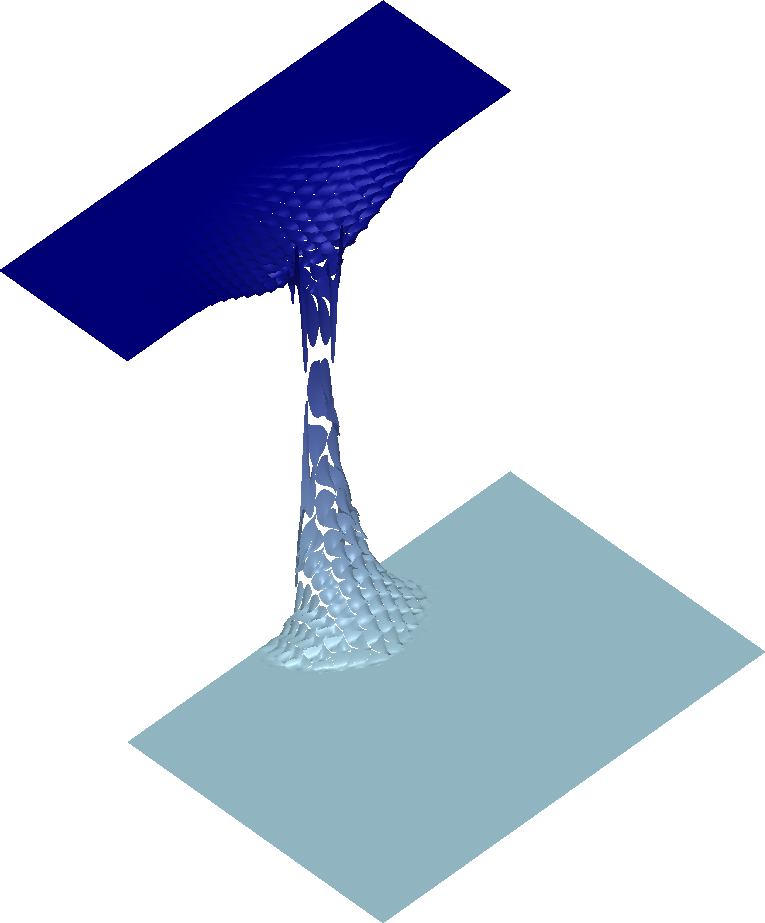}
\caption{Top view (top row) and side view (bottom row) of the dam break problem at T = 0.25. Left: node-wise limiting. Middle: element-wise limiting. Right: low order scheme}
\label{fig:Dam_n1}
\end{center}
\end{figure}





\begin{figure}
\begin{center}
\includegraphics[width=0.32\textwidth, height = 0.2\textheight,
trim={0 0 5cm 0},clip]{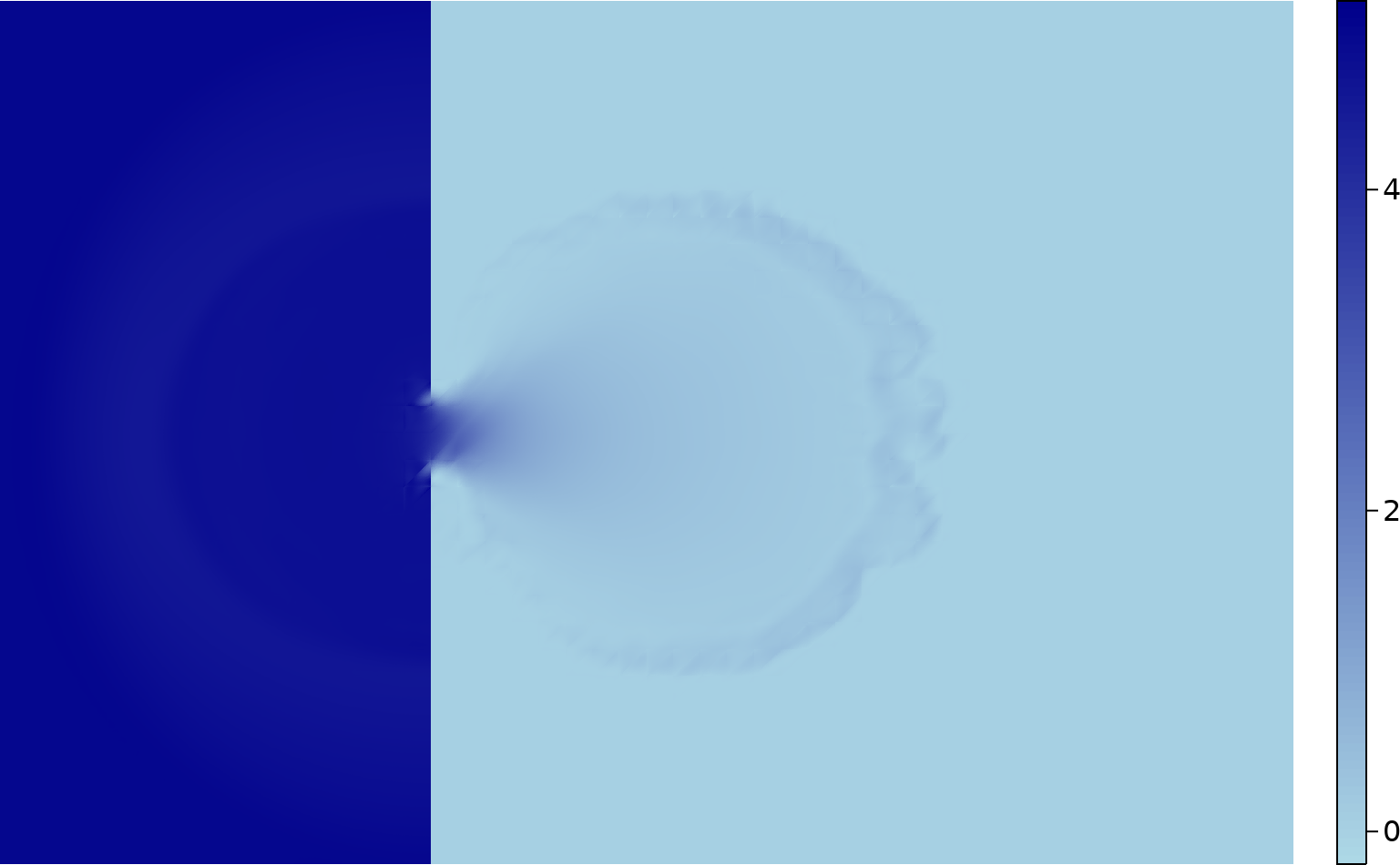}
\includegraphics[width=0.32\textwidth, height = 0.2\textheight, 
trim={0 0 5cm 0},clip]{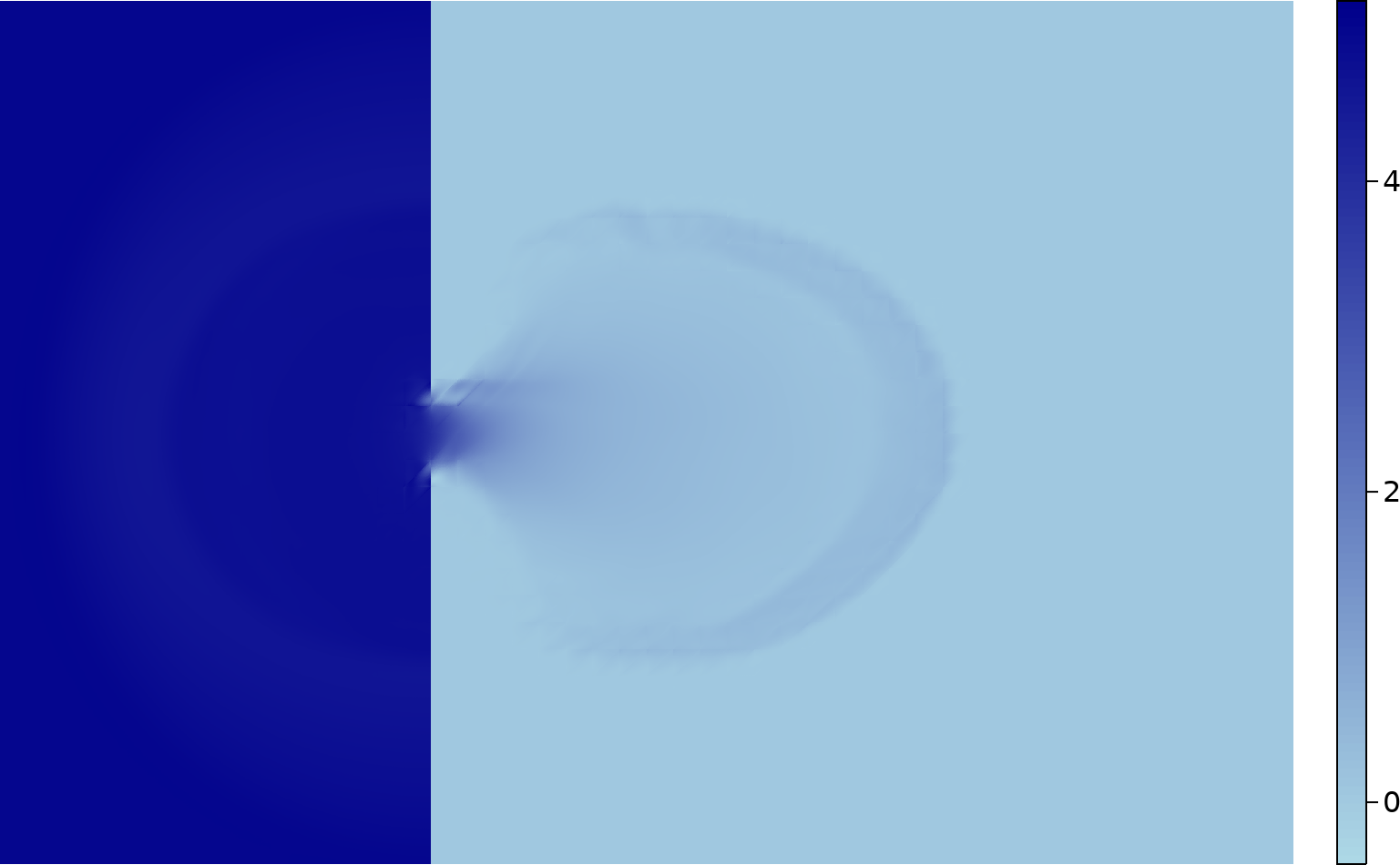}
\includegraphics[width=0.32\textwidth, height = 0.2\textheight, 
trim={0 0 5cm 0.3cm},clip]{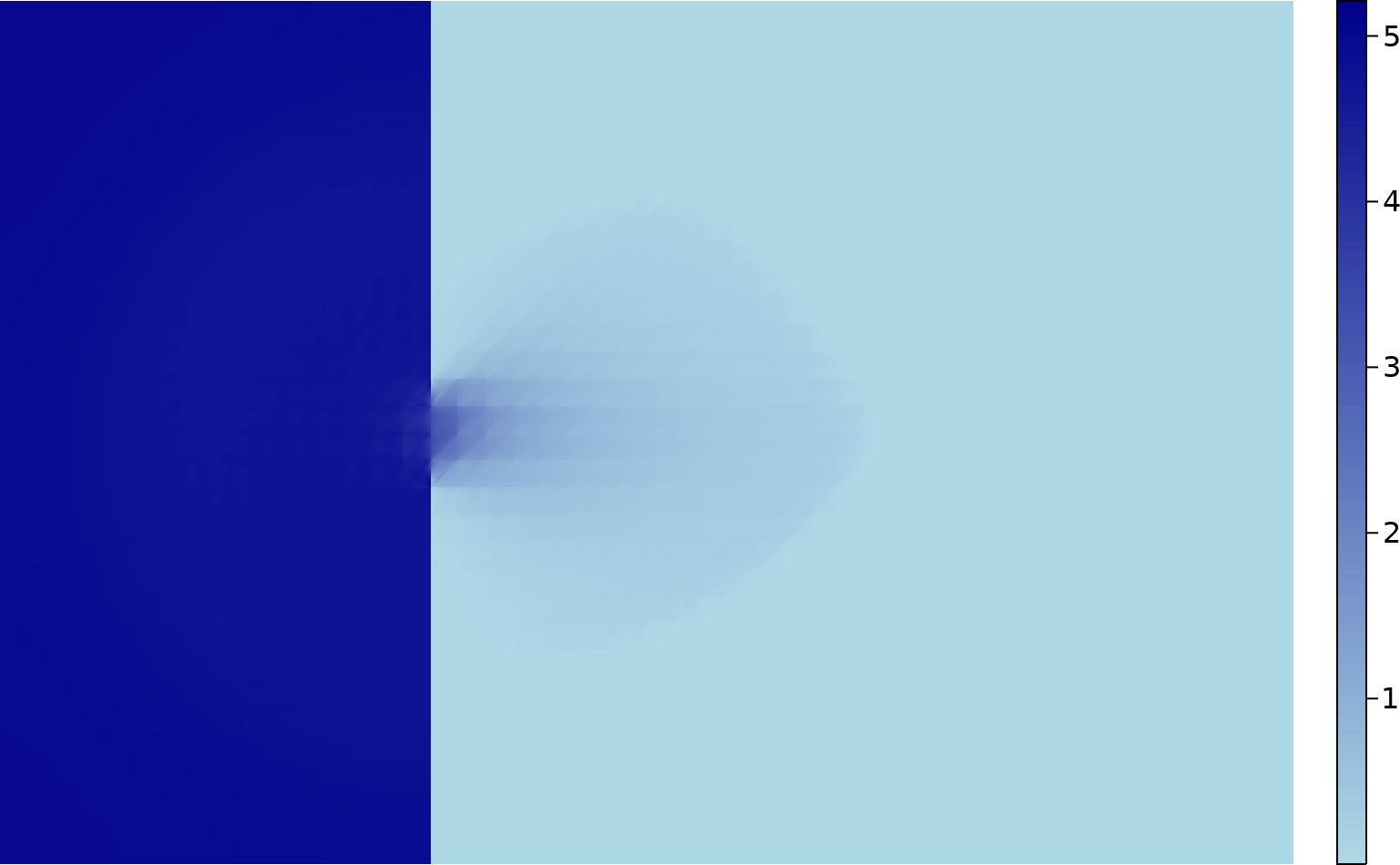}\\

\includegraphics[width=0.32\textwidth, height = 0.21\textheight]{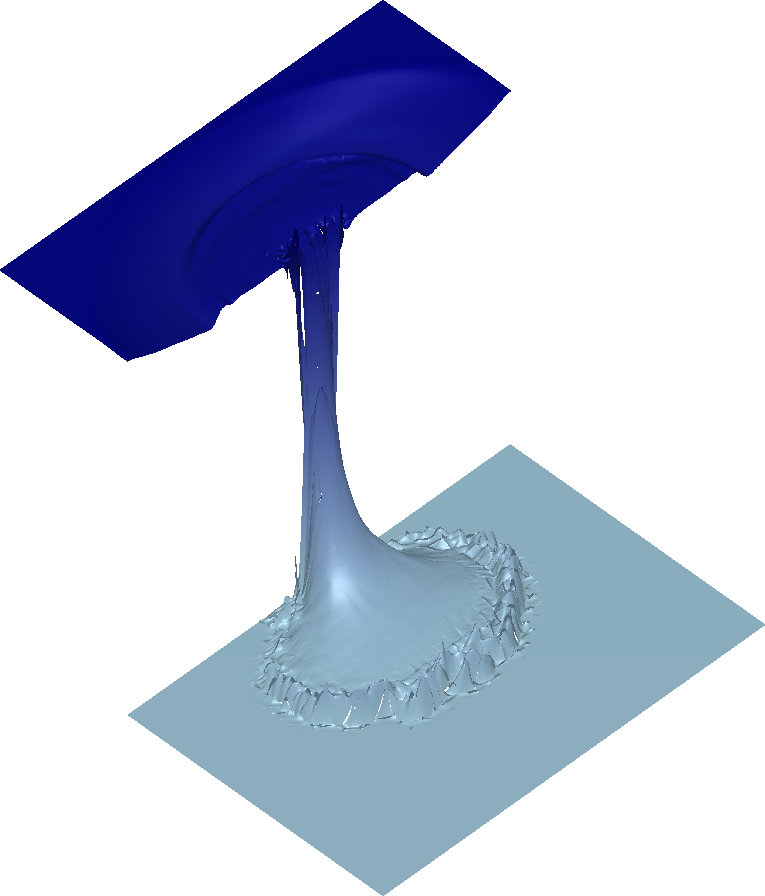}
\includegraphics[width=0.32\textwidth, height = 0.21\textheight]{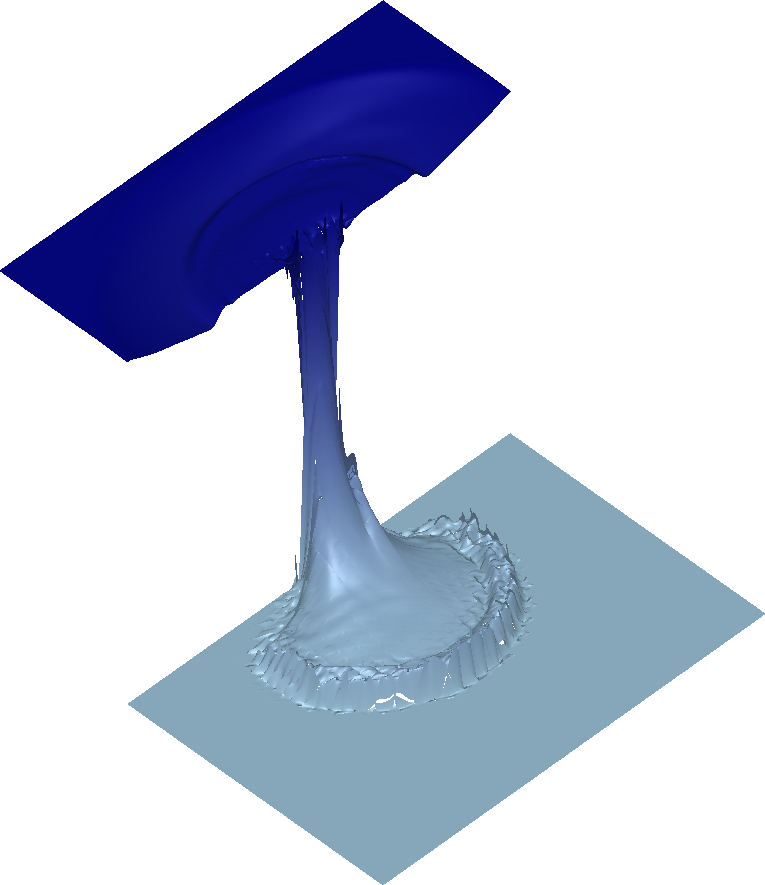}
\includegraphics[width=0.32\textwidth, height = 0.21\textheight]{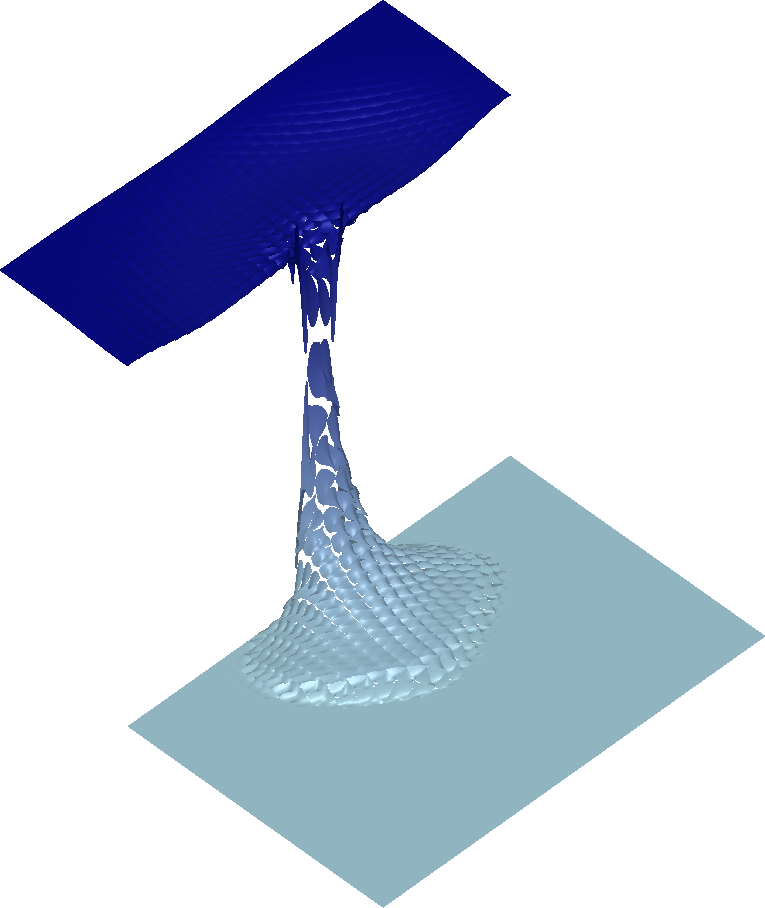}
\caption{Top view (top row) and side view (bottom row) of the dam break problem at T = 0.5. Left: node-wise limiting. Middle: element-wise limiting. Right: low order scheme}
\label{fig:Dam_n2}
\end{center}
\end{figure}





\begin{figure}
\begin{center}
\includegraphics[width=0.32\textwidth, height = 0.2\textheight,
trim={0 0 5cm 0},clip]{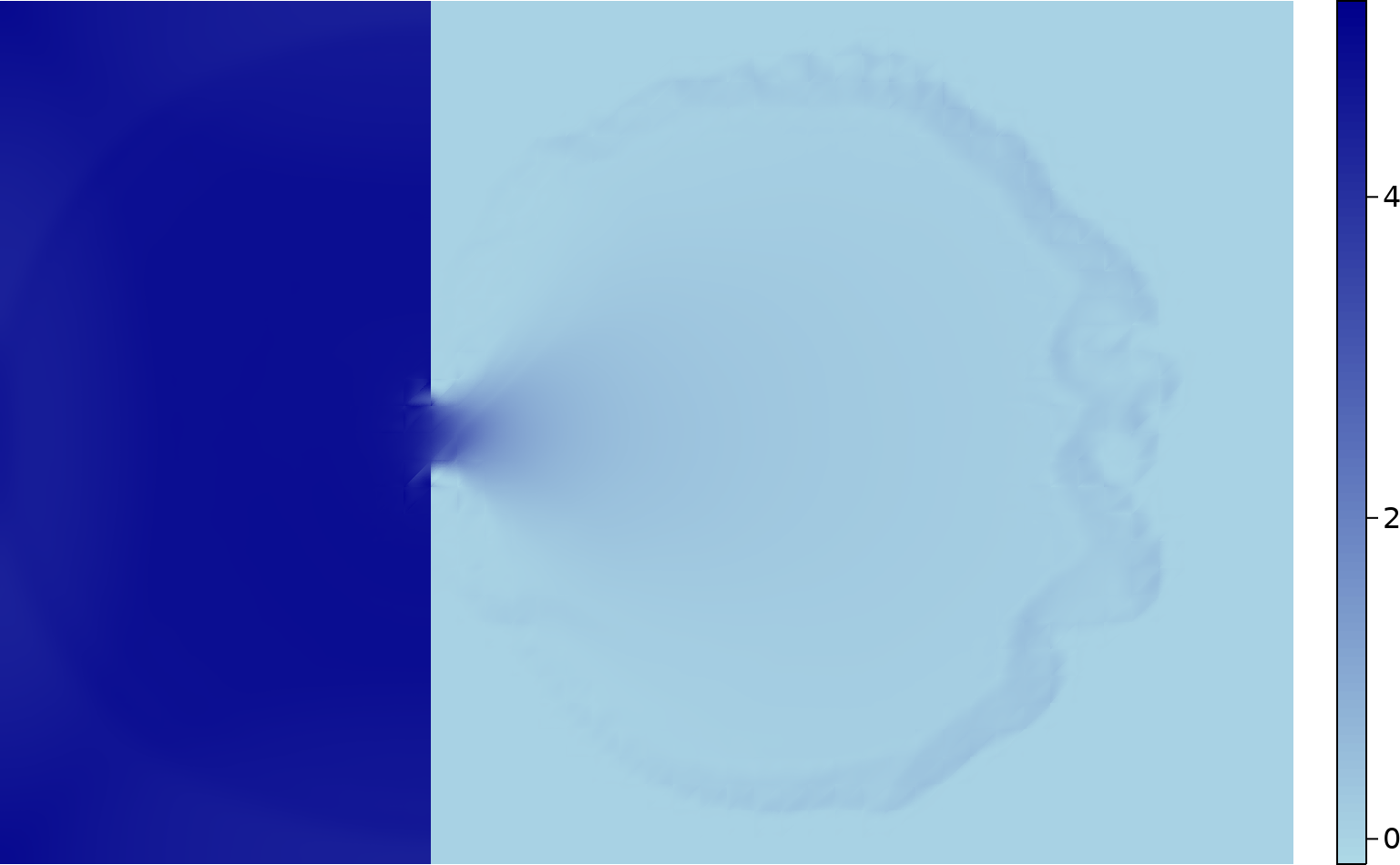}
\includegraphics[width=0.32\textwidth, height = 0.2\textheight, 
trim={0 0 5cm 0},clip]{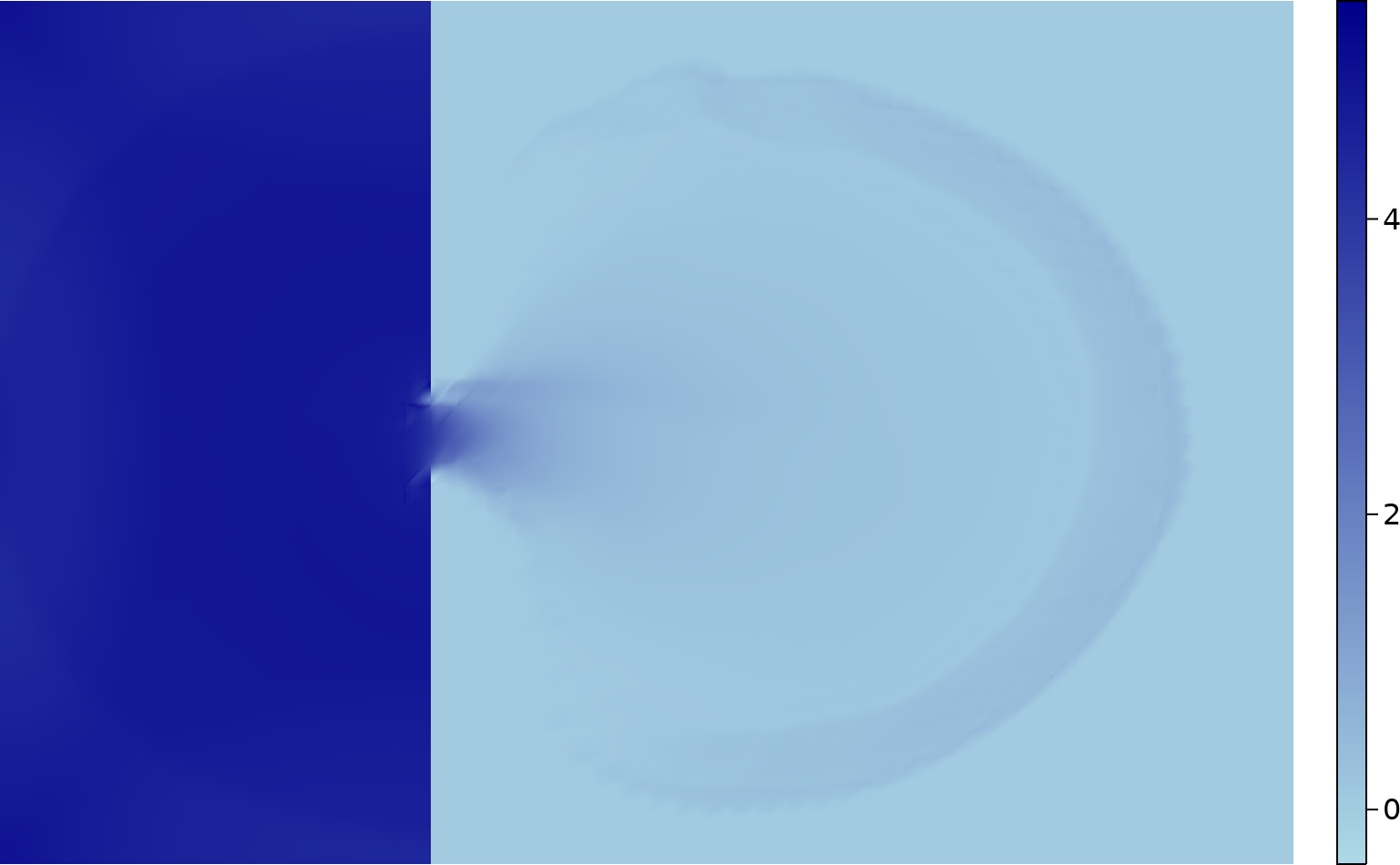}
\includegraphics[width=0.32\textwidth, height = 0.2\textheight, 
trim={0 0 5cm 0.3cm},clip]{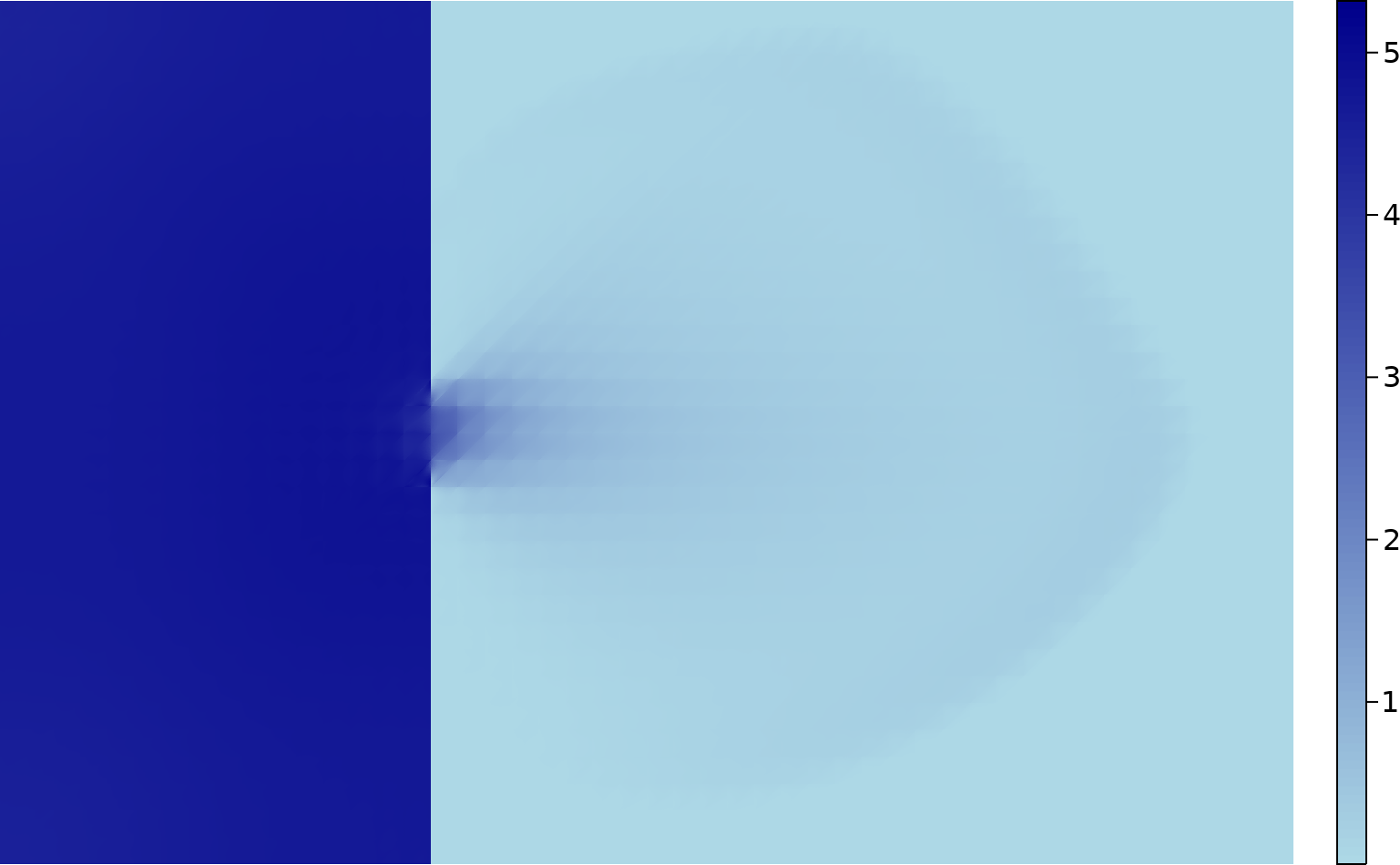}\\

\includegraphics[width=0.32\textwidth, height = 0.21\textheight]{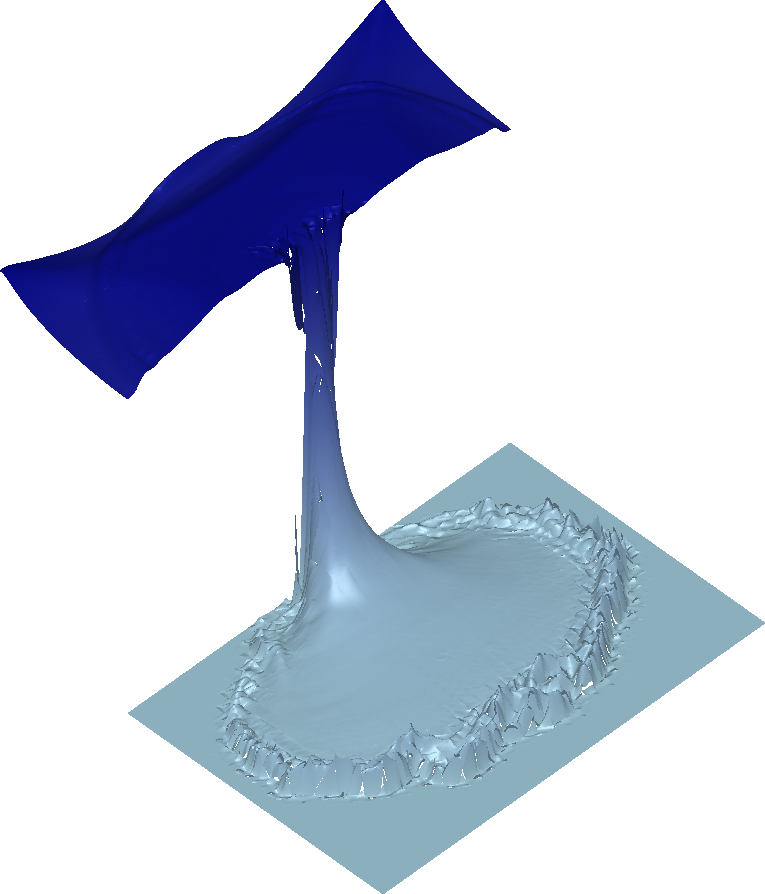}
\includegraphics[width=0.32\textwidth, height = 0.21\textheight]{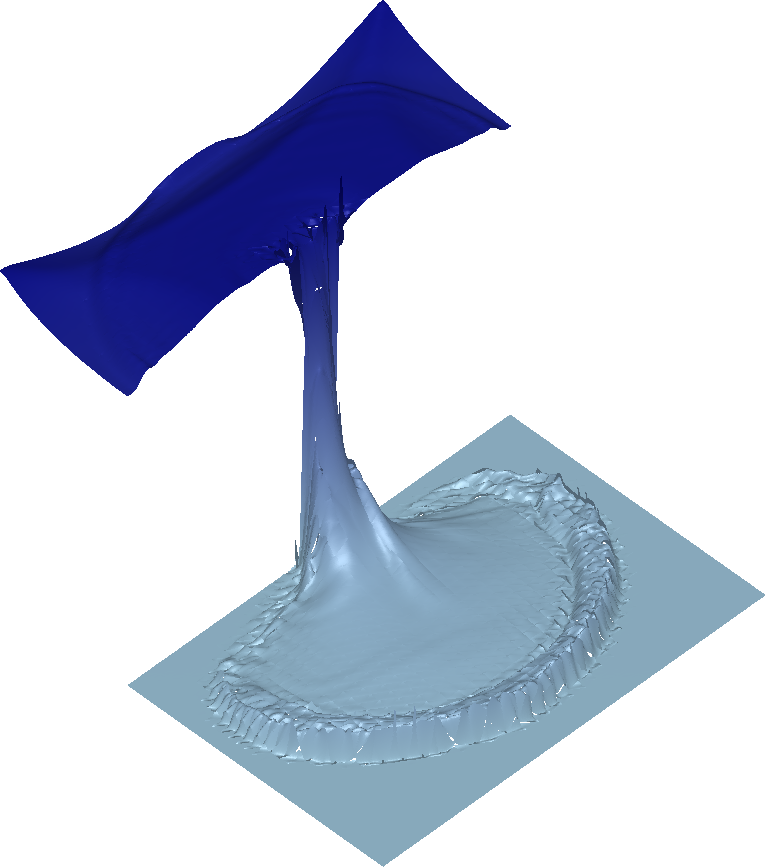}
\includegraphics[width=0.32\textwidth, height = 0.21\textheight]{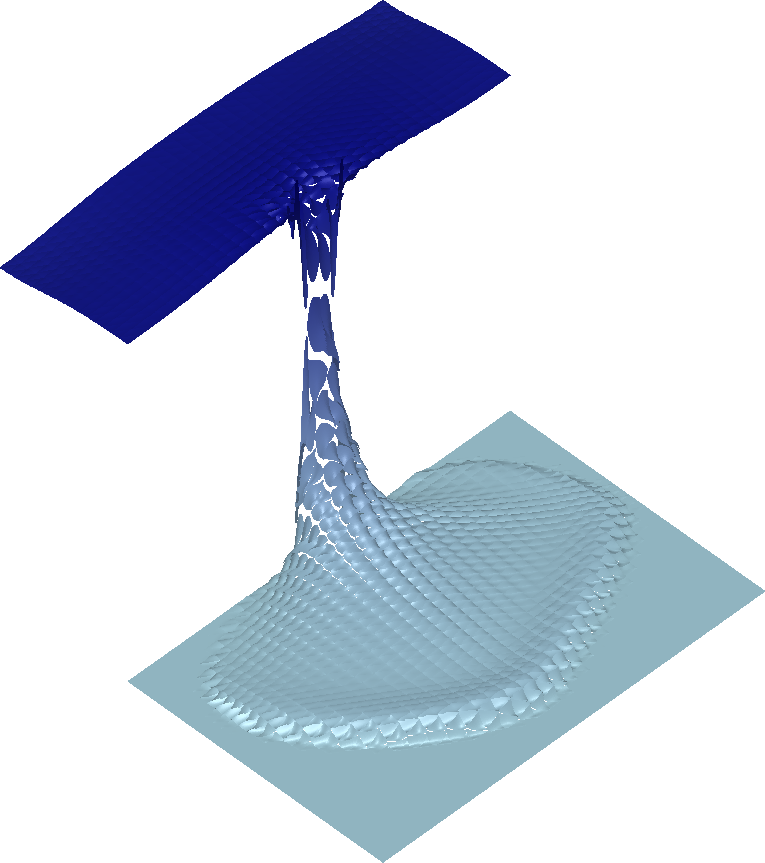}
\caption{Top view (top row) and side view (bottom row) of the dam break problem at T = 0.75. Left: node-wise limiting. Middle: element-wise limiting. Right: low order scheme}
\label{fig:Dam_n3}
\end{center}
\end{figure}





\begin{figure}
\begin{center}
\includegraphics[width=0.32\textwidth, height = 0.2\textheight,
trim={0 0 5cm 0},clip]{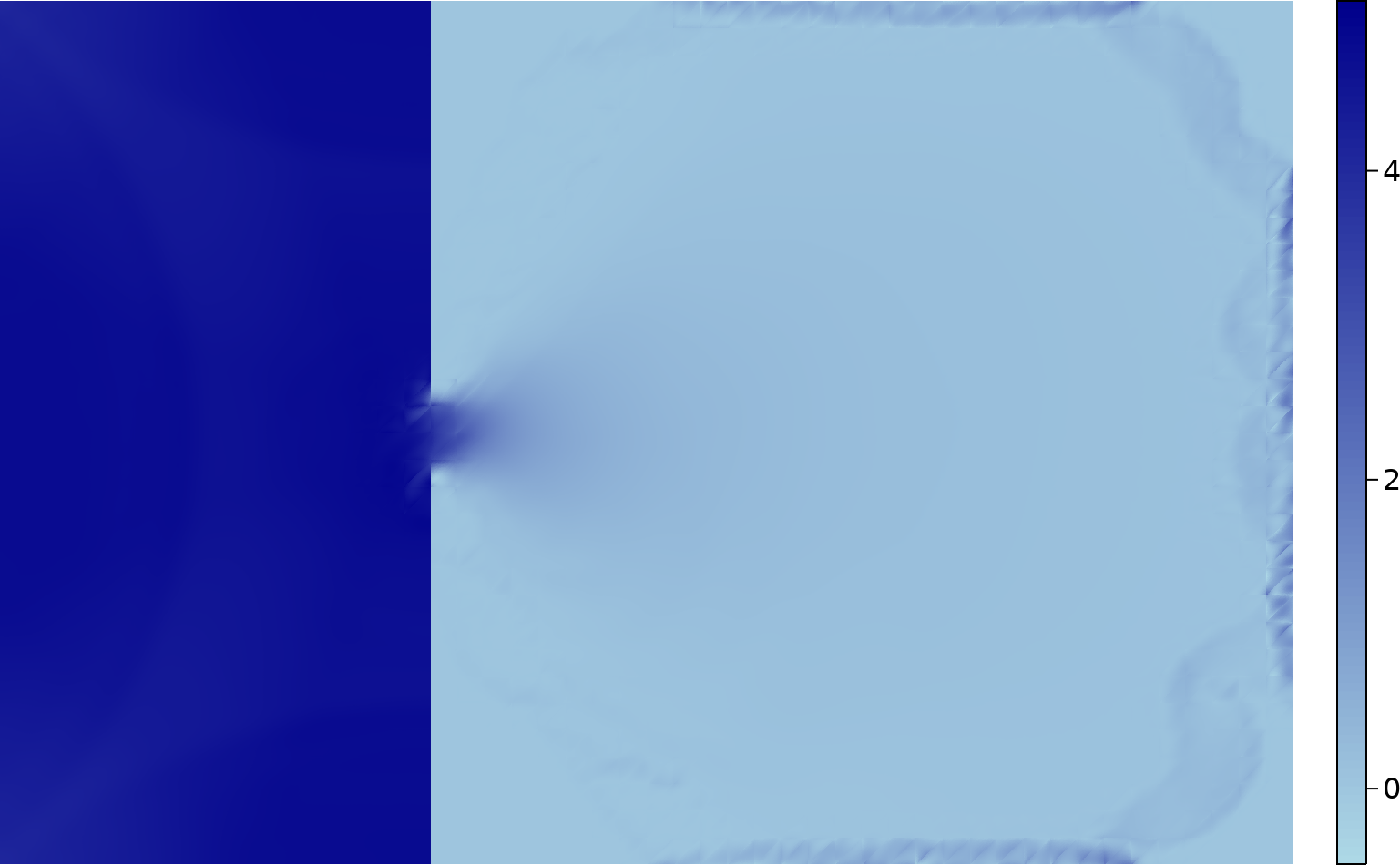}
\includegraphics[width=0.32\textwidth, height = 0.2\textheight, 
trim={0 0 5cm 0},clip]{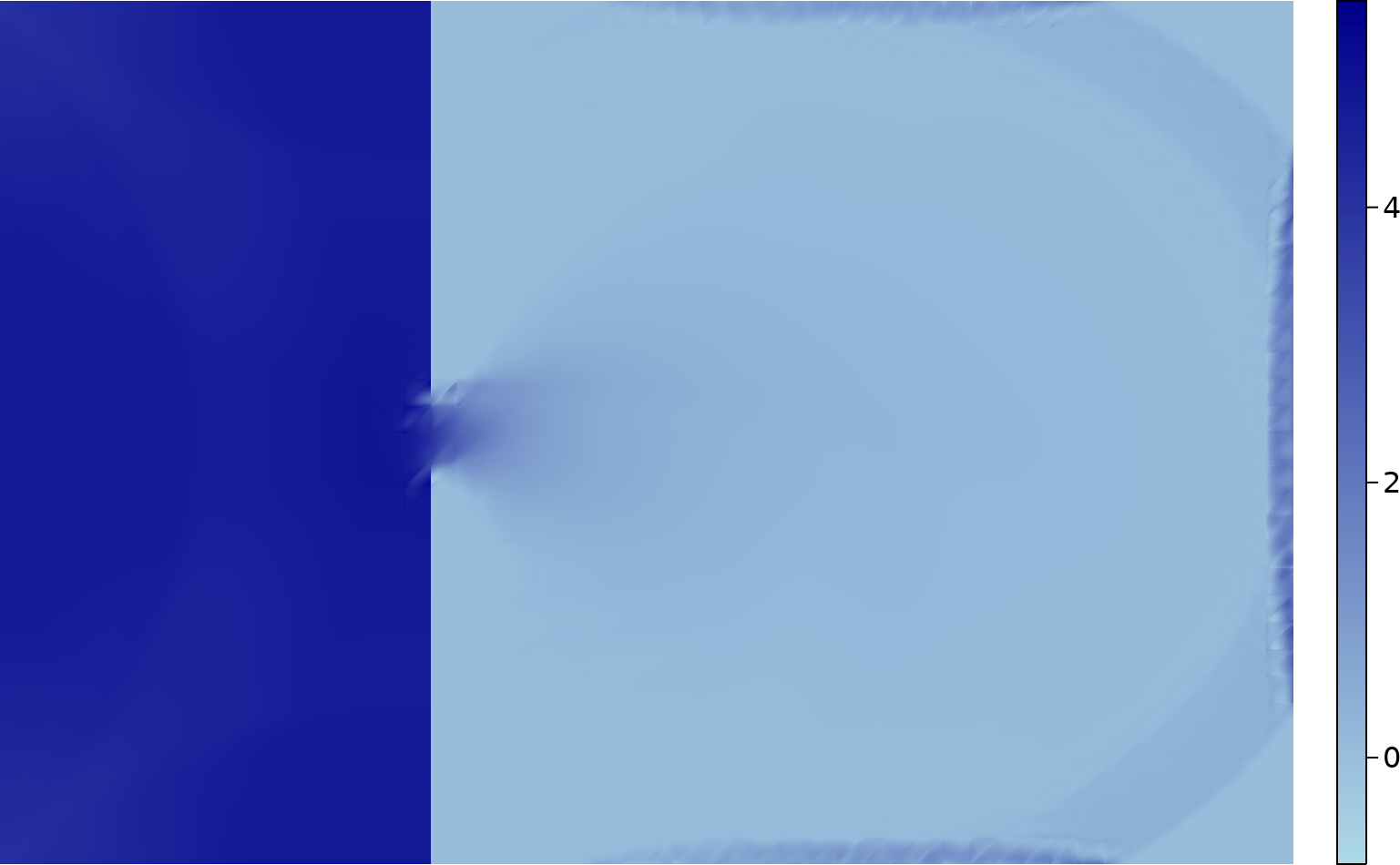}
\includegraphics[width=0.32\textwidth, height = 0.2\textheight, 
trim={0 0 5cm 0.3cm},clip]{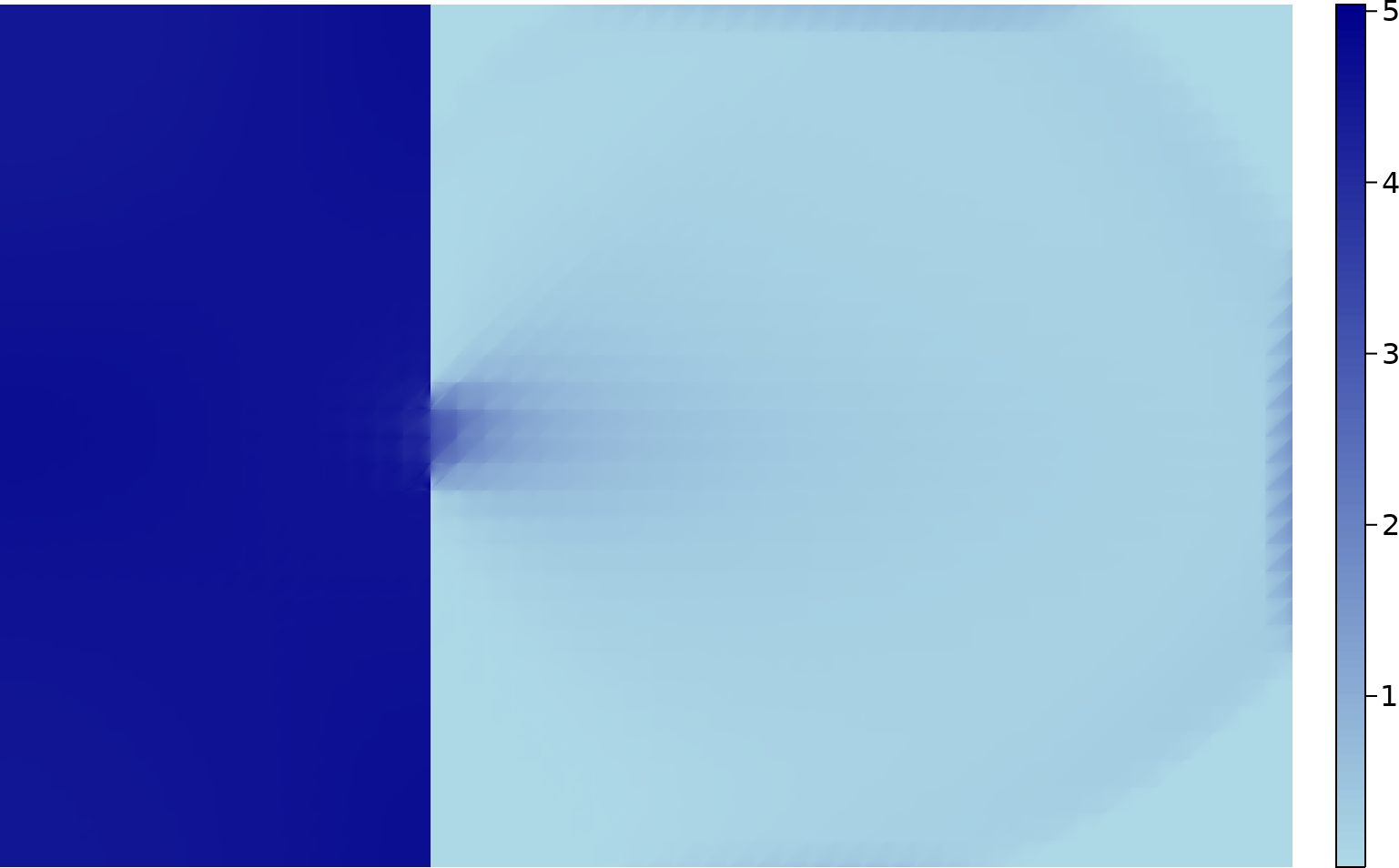}\\

\includegraphics[width=0.32\textwidth, height = 0.21\textheight]{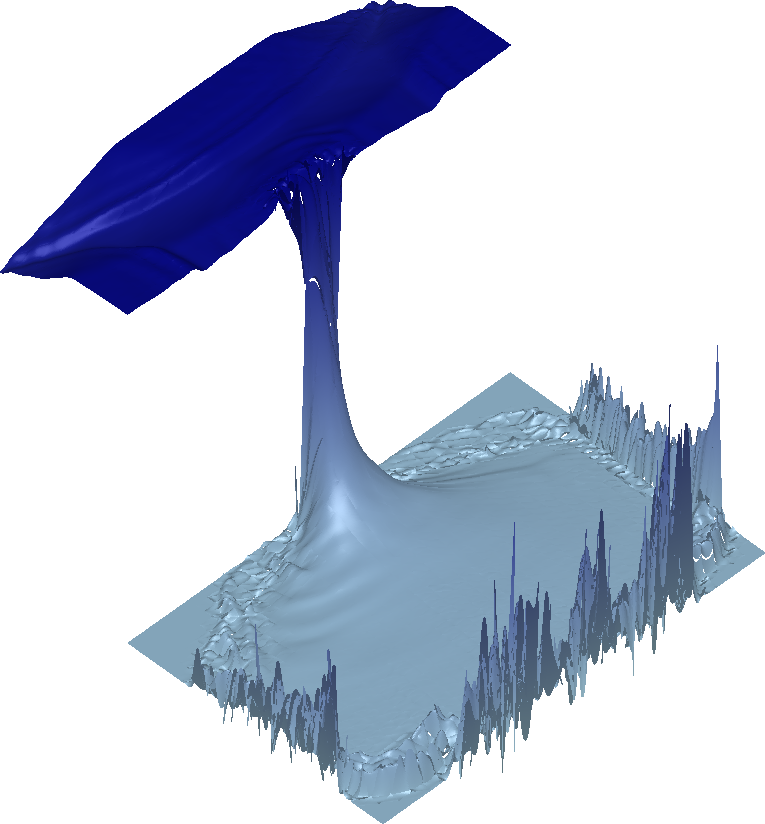}
\includegraphics[width=0.32\textwidth, height = 0.21\textheight]{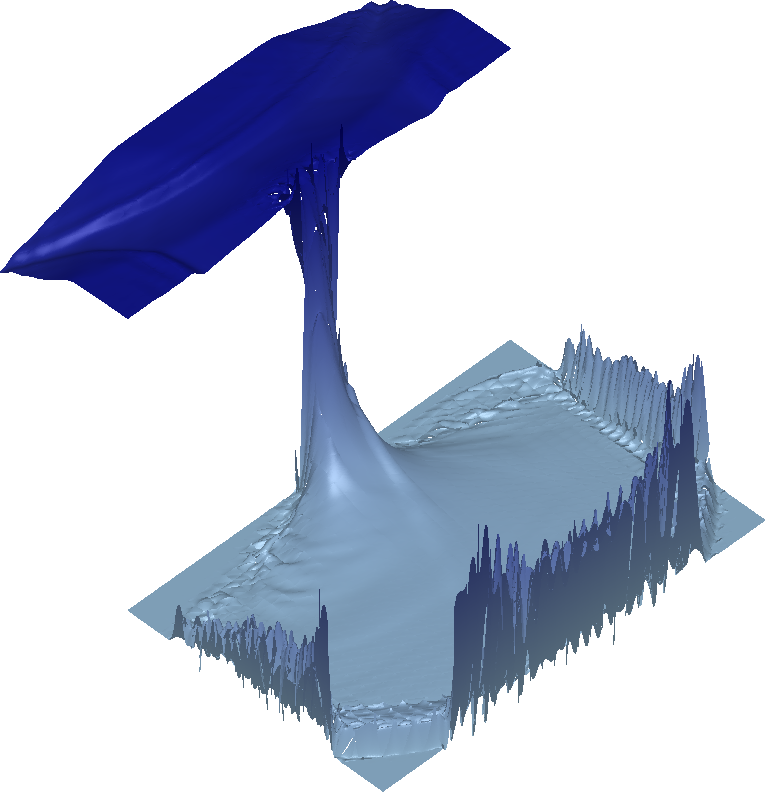}
\includegraphics[width=0.32\textwidth, height = 0.21\textheight]{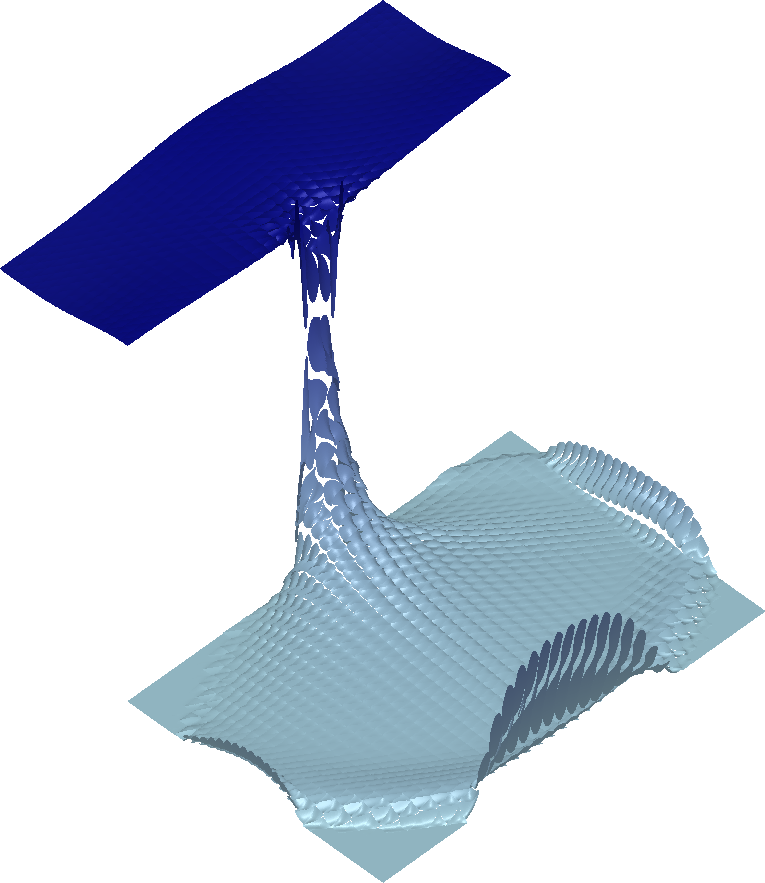}
\caption{Top view (top row) and side view (bottom row) of the dam break problem at T = 1. Left: node-wise limiting. Middle: element-wise limiting. Right: low order scheme}
\label{fig:Dam_n4}
\end{center}
\end{figure}




\subsection{Wave over a bump}
The last experiment we perform is a 2D simulation of a wave passing over a bump. The bottom geometry is modeled by a Gaussian bump with a height such that the bump is never entirely submerged. This experiment provides a test of how our numerical solution handles a persistent dry area. We use polynomial degree $N=3$ with a total of 8192 ($64\times 64 \times 2$) elements. This experiment runs on the domain $[-1,1]^2$ with periodic boundary conditions. The bottom topography is given by the following expression:
\begin{align}
    b(x,y) = 5\exp(-25(x^2+y^2)),
\end{align}
and the initial condition for this experiment is
\begin{align}
    &h(x,y) = \max(0,\exp(-25((x.+0.5).^2))+2 - b),\\
    &u(x,y) = \begin{cases}
    1 & {\rm if } \, h > 2,\\
    0 &{\rm otherwise}, 
    \end{cases} \qquad
    v = 0.
\end{align}

\begin{figure}
\begin{center}
\includegraphics[width=0.48\textwidth]{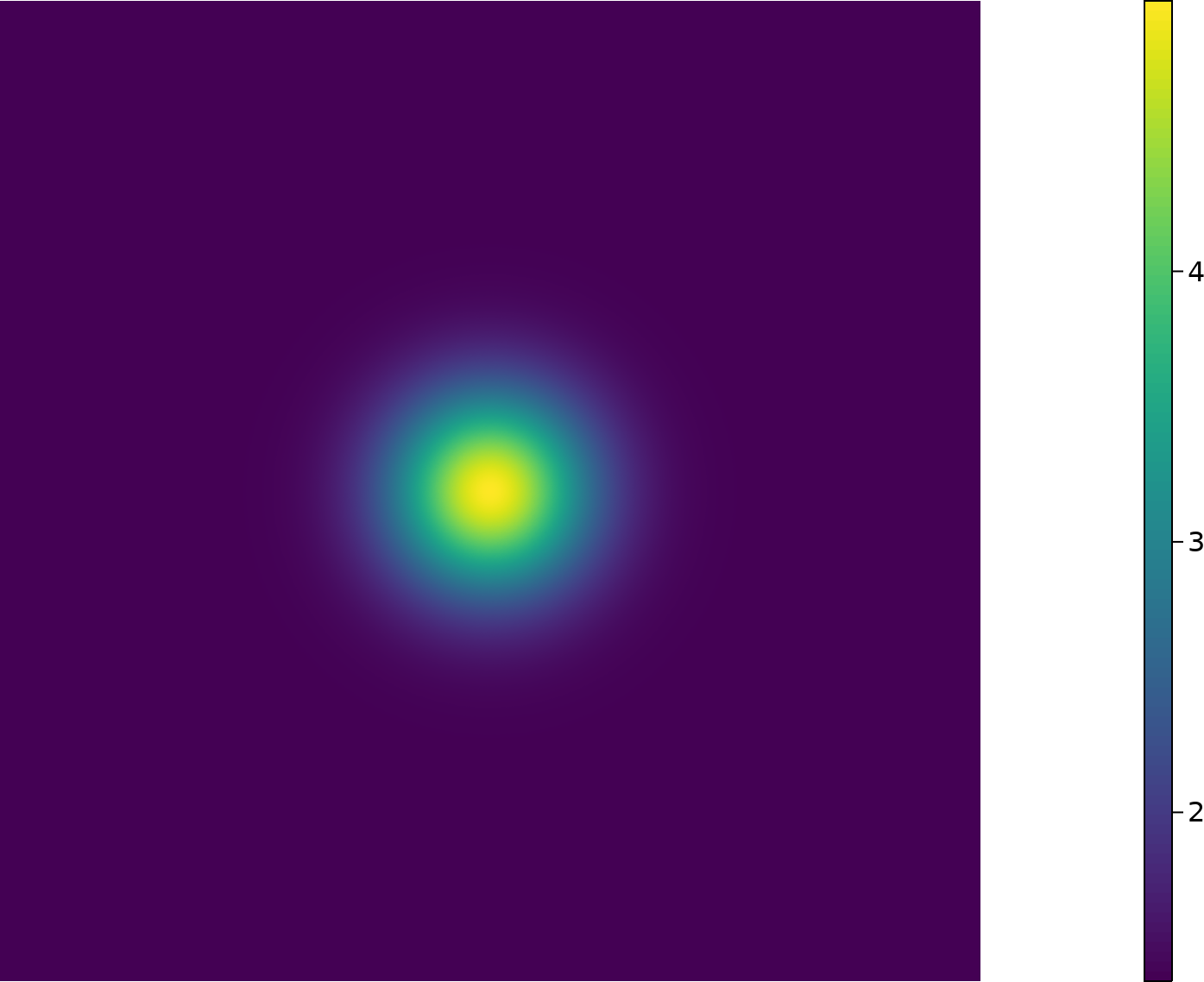}
\includegraphics[width=0.48\textwidth, height = 0.29\textheight] {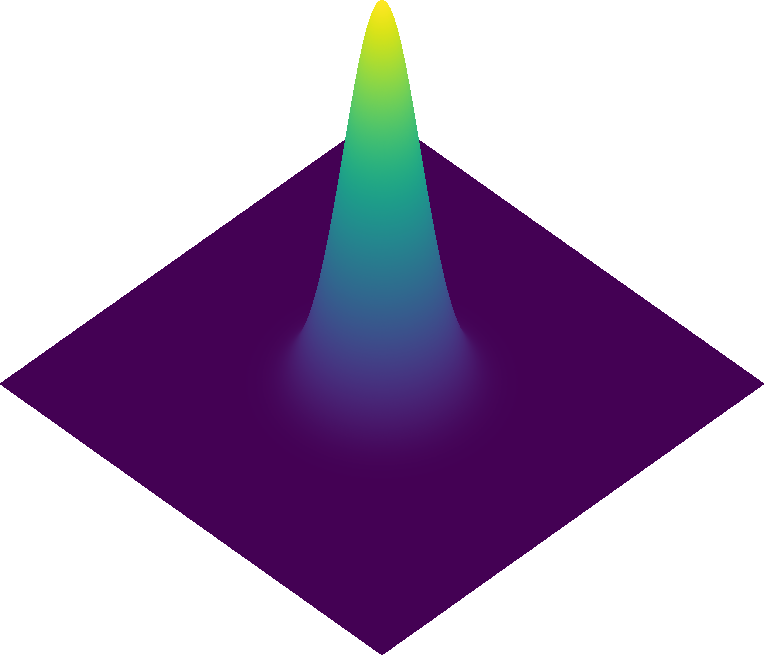}\hspace{0.06\textwidth}\\
\vspace{0.6cm}
\includegraphics[width=0.48\textwidth]{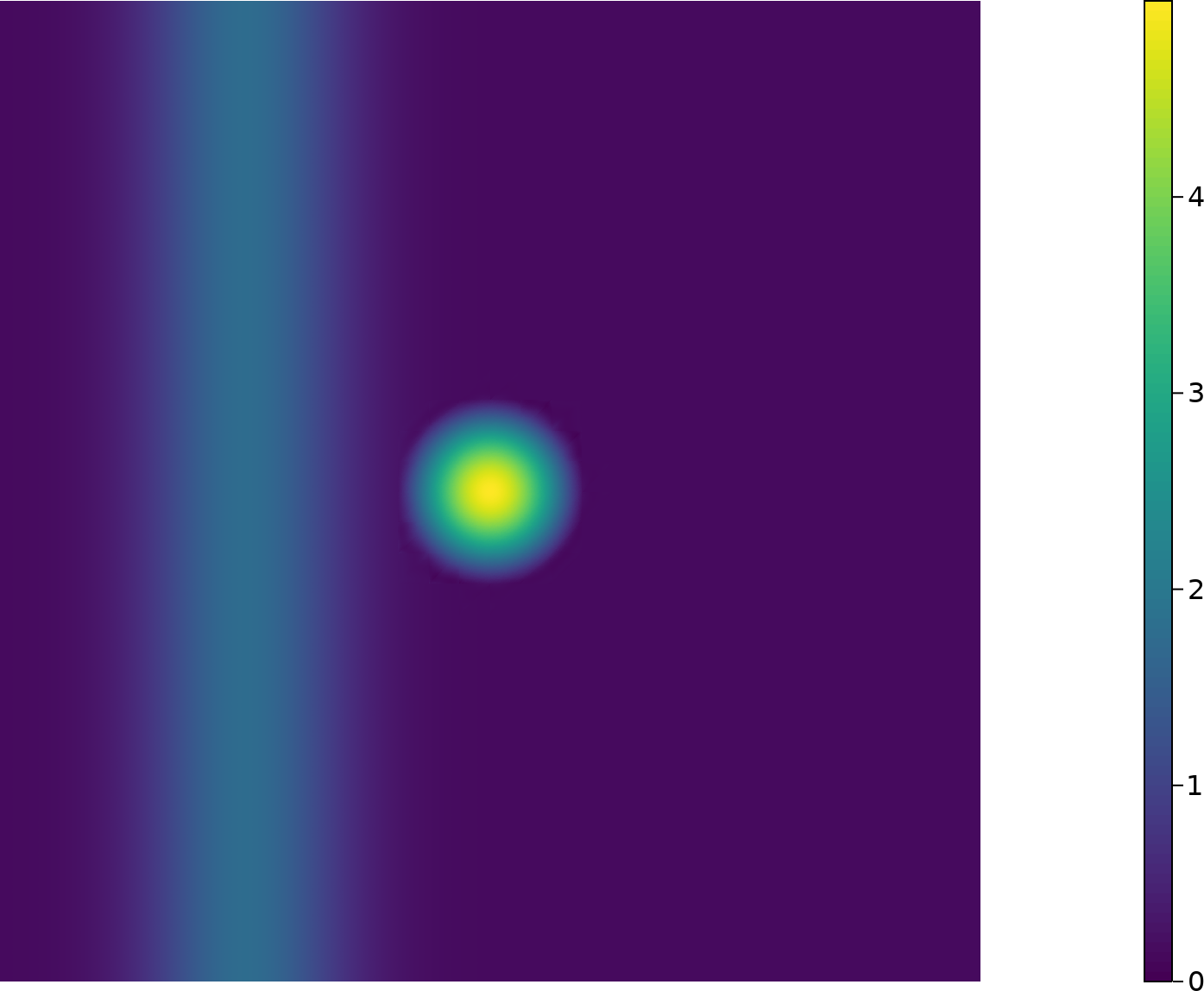}
\includegraphics[width=0.48\textwidth, height = 0.29\textheight]
{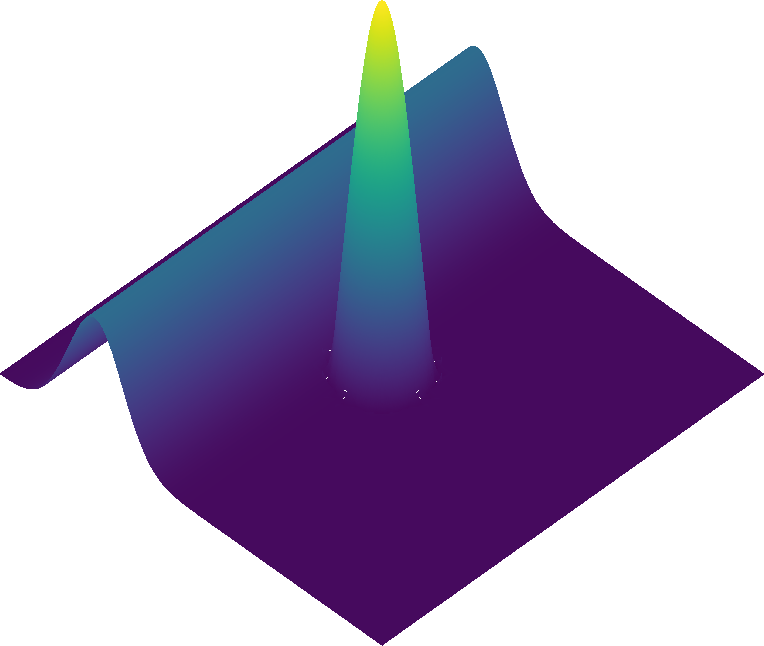}\\
\caption{Bottom geometry and initial condition for wave over a bump: Top: top view and side view of bottom geometry. Bottom; top view and side view of initial water surface and dry area.}
\label{fig:CL_wav0}
\end{center}
\end{figure}
We run this experiment up to time T = 1 and plot numerical solutions from three different schemes. In Fig. \ref{fig:CL_wave_1} through Fig. \ref{fig:CL_wave_4}, we present the solutions from node-wise limiting scheme, element-wise limiting scheme, and the lower order positivity preserving scheme. For all schemes, we observe that the wave propagates along the $x$-axis, hits the bump and splits, and converges after passing it. We notice that our solution contains small oscillations at the wavefront but never blows up. The water height remains positive even in the presence of a dry area throughout the entire experiment. We observe that the results from the node-wise limiting scheme are very similar to the results from the element-wise limiting scheme for this problem.
\begin{figure}
\begin{center}
\includegraphics[width=0.32\textwidth, height = 0.2\textheight,
trim={0 0 6cm 0},clip]{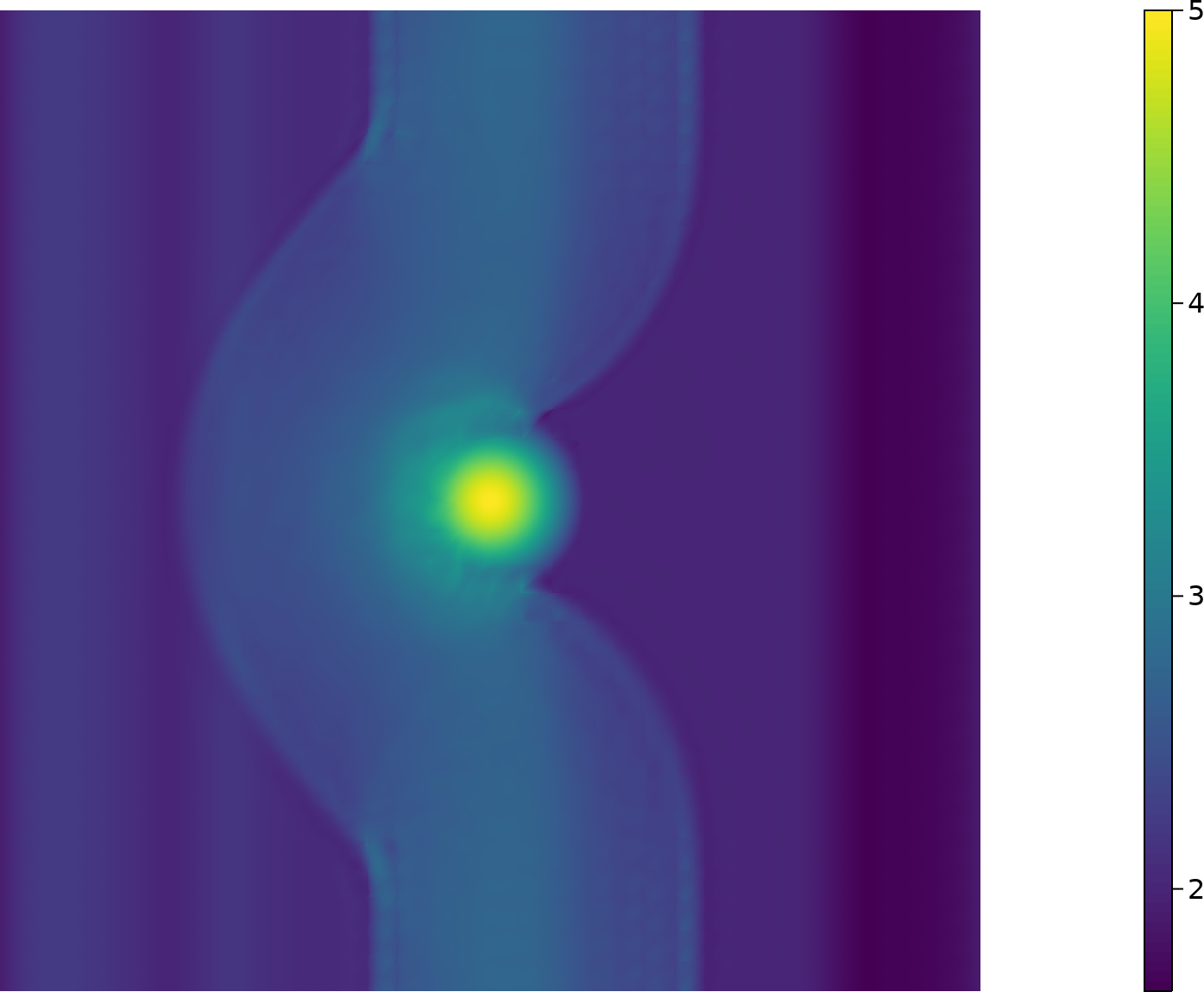}
\includegraphics[width=0.32\textwidth, height = 0.2\textheight, 
trim={0 0 6cm 0},clip]{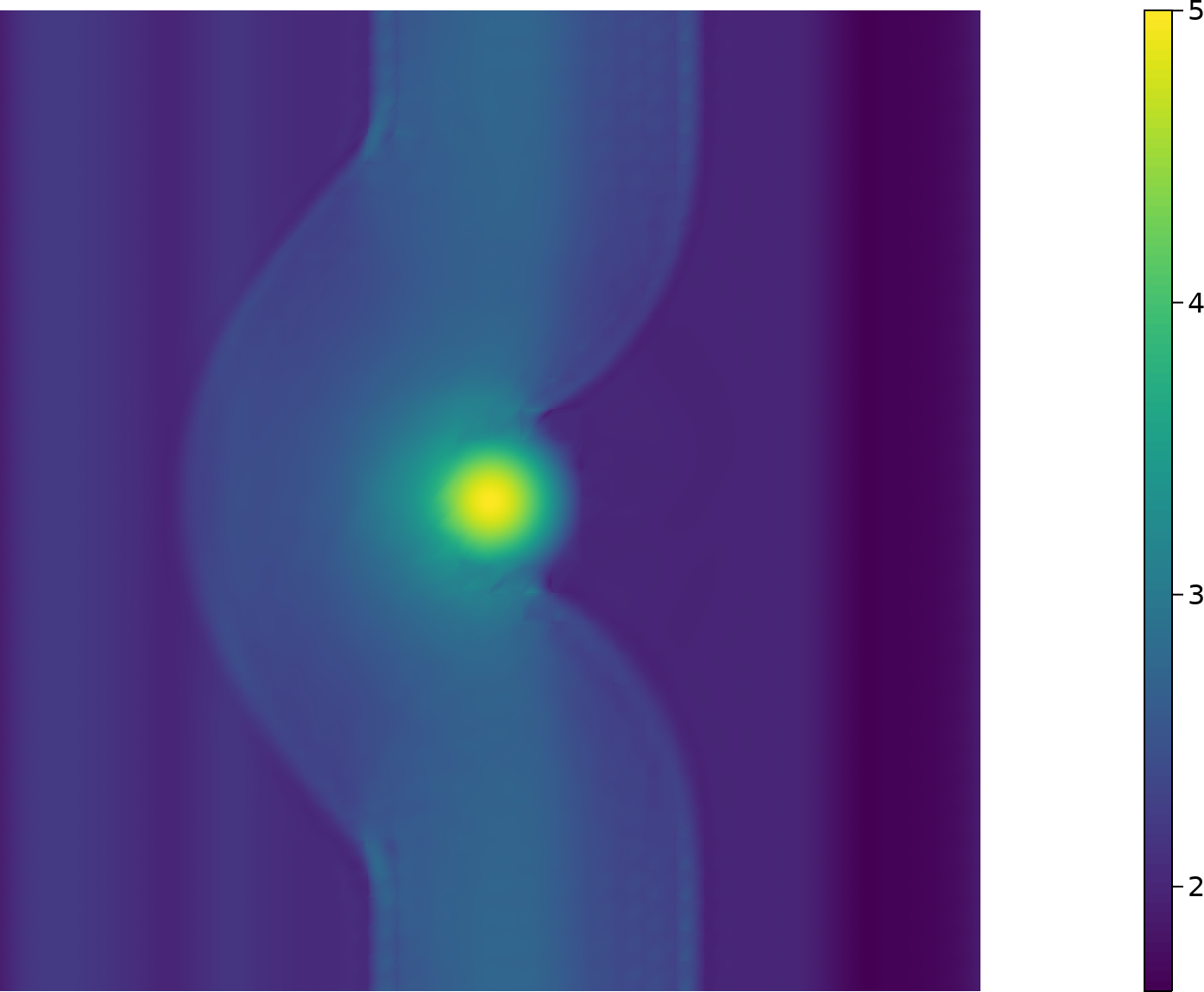}
\includegraphics[width=0.32\textwidth, height = 0.2\textheight, 
trim={0 0 6cm 0 cm},clip]{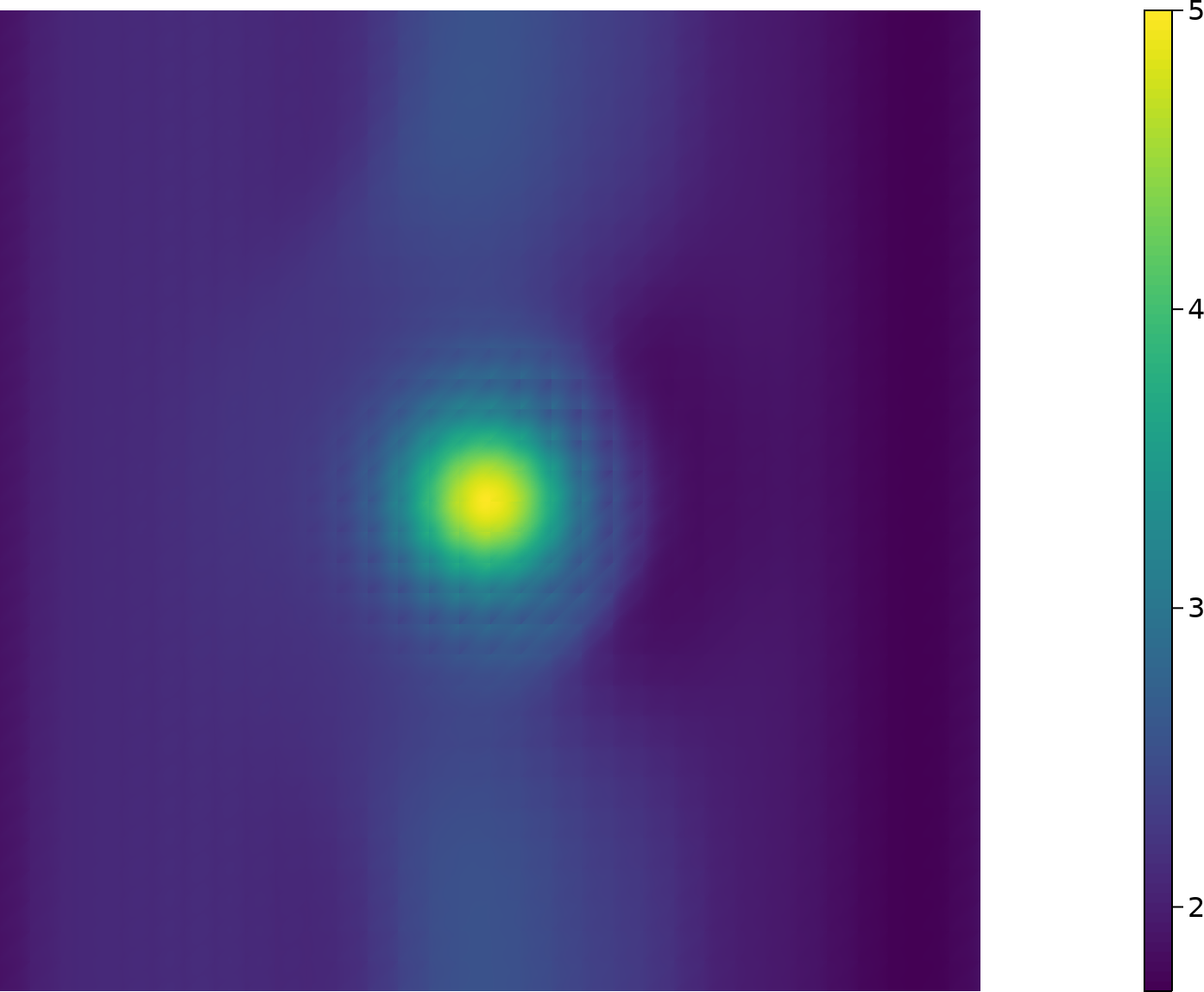}\\

\includegraphics[width=0.32\textwidth, height = 0.19\textheight]{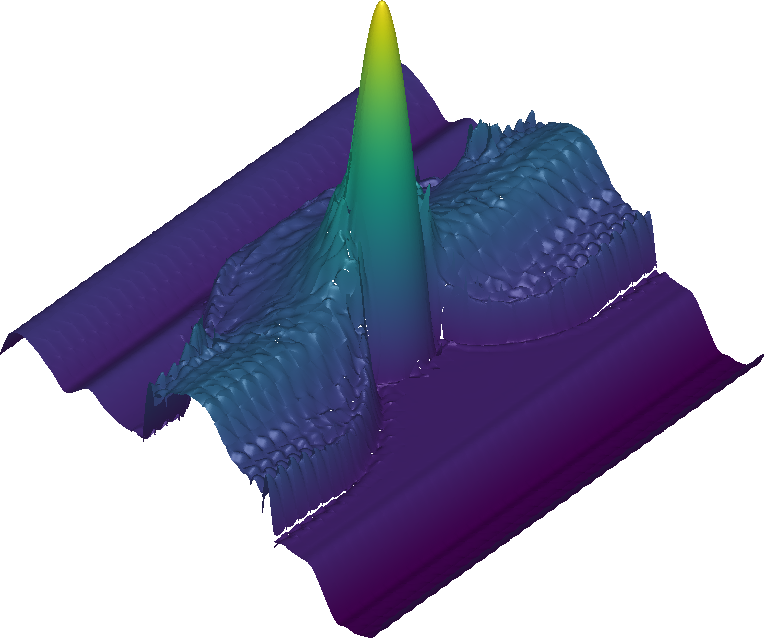}
\includegraphics[width=0.32\textwidth, height = 0.19\textheight]{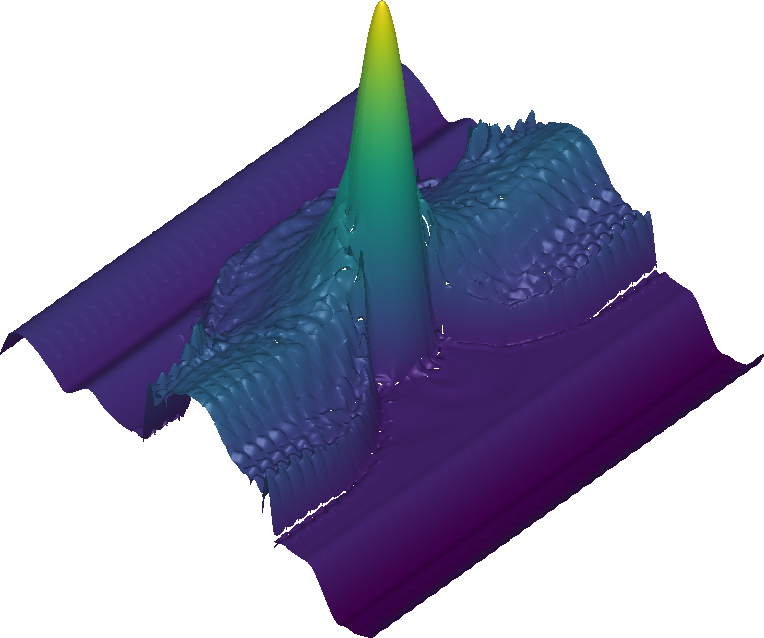}
\includegraphics[width=0.32\textwidth, height = 0.19\textheight]{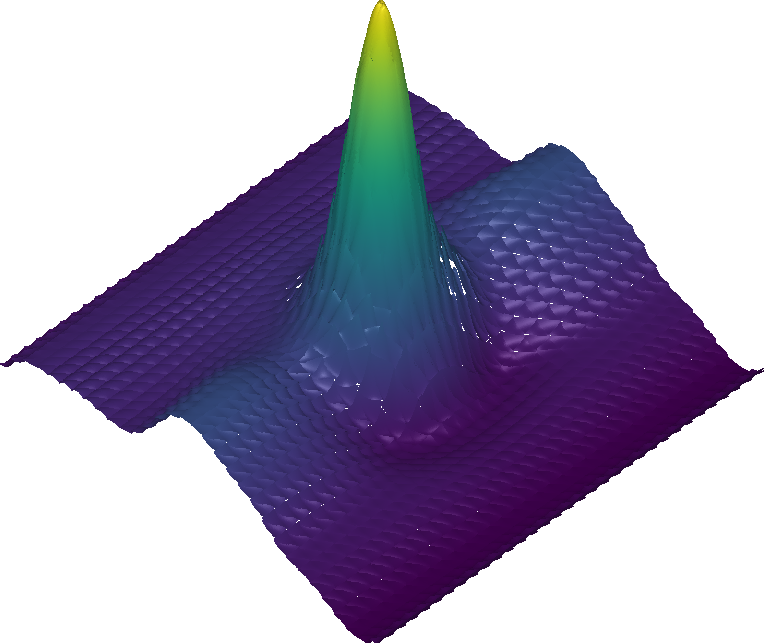}
\caption{Top view (top row) and side view (bottom row) of water surface and dry area for wave over a bump at T = 0.25. Left: node-wise limiting. Middle: element-wise limiting. Right: low order scheme}
\label{fig:CL_wave_1}
\end{center}
\end{figure}




\begin{figure}
\begin{center}
\includegraphics[width=0.32\textwidth, height = 0.2\textheight,
trim={0 0 6cm 0},clip]{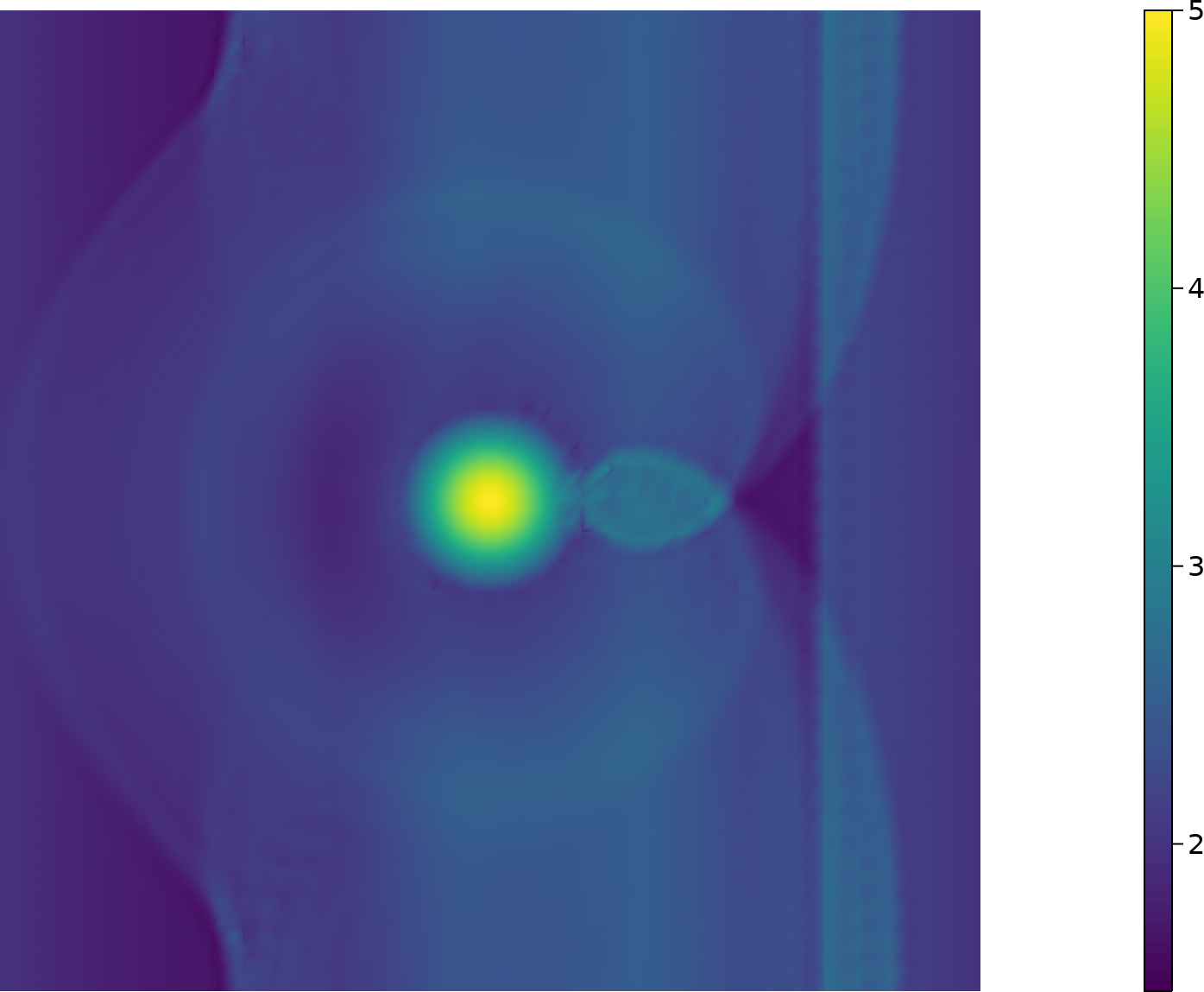}
\includegraphics[width=0.32\textwidth, height = 0.2\textheight, 
trim={0 0 6cm 0},clip]{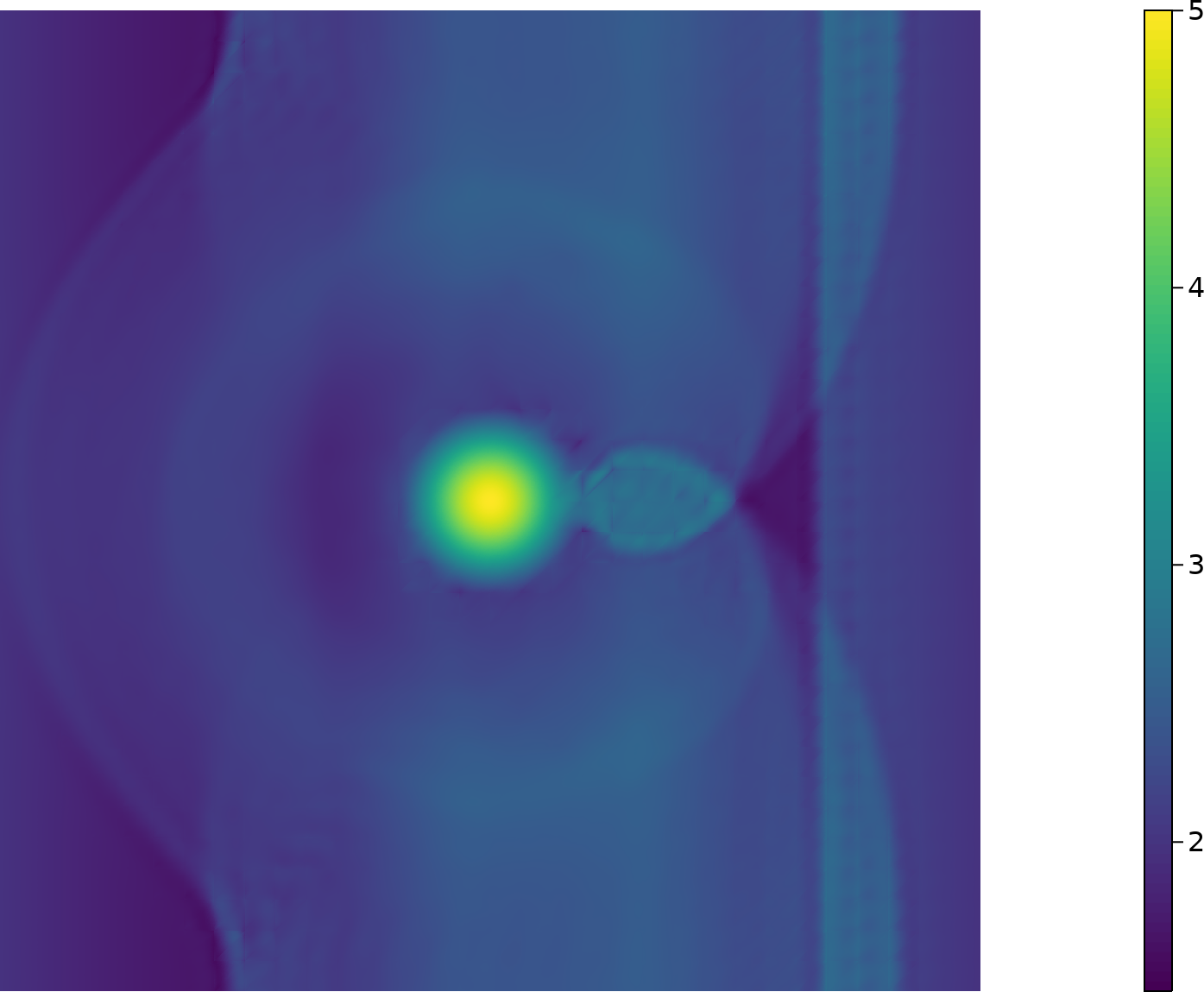}
\includegraphics[width=0.32\textwidth, height = 0.2\textheight, 
trim={0 0 6cm 0 cm},clip]{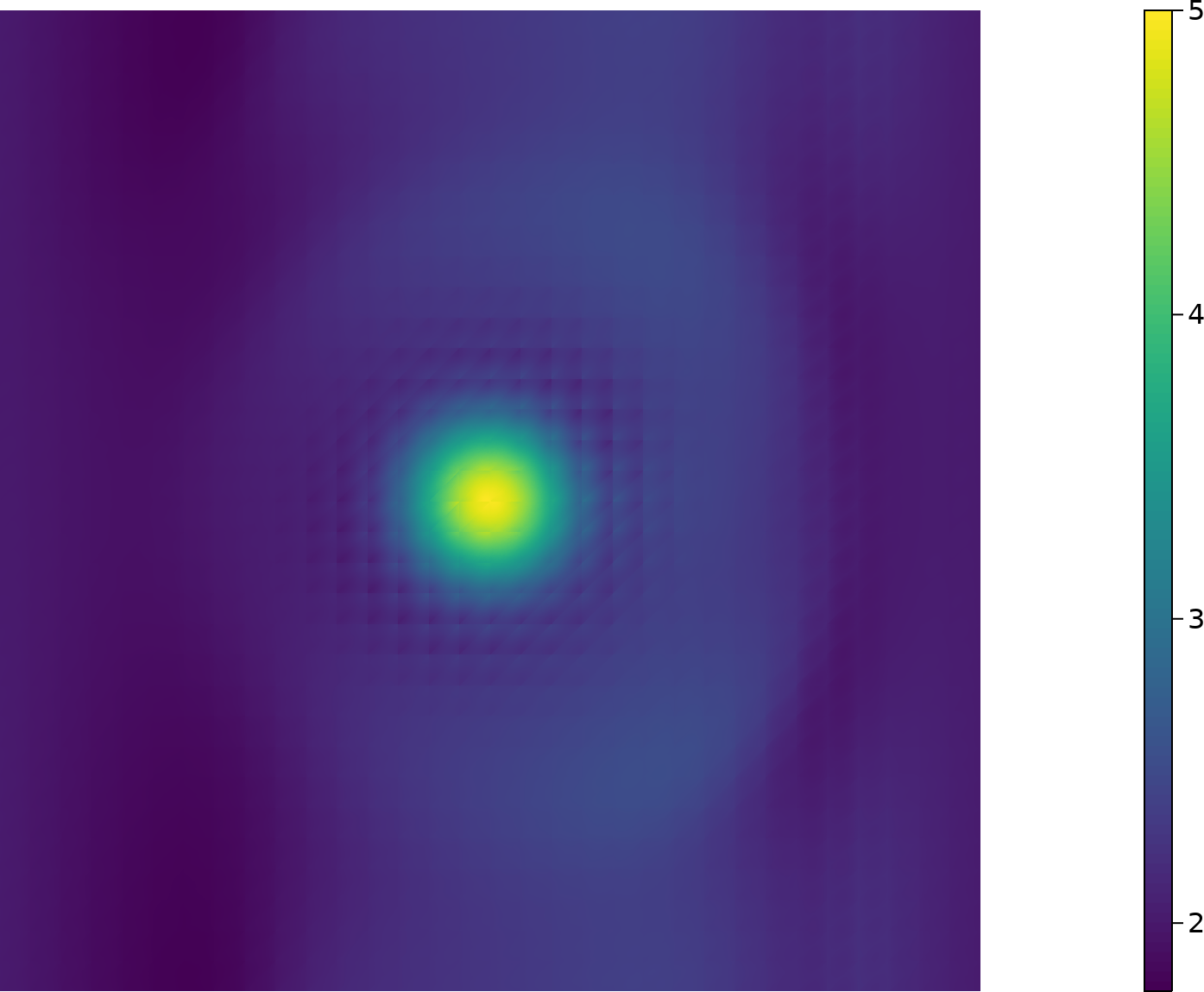}\\

\includegraphics[width=0.32\textwidth, height = 0.19\textheight]{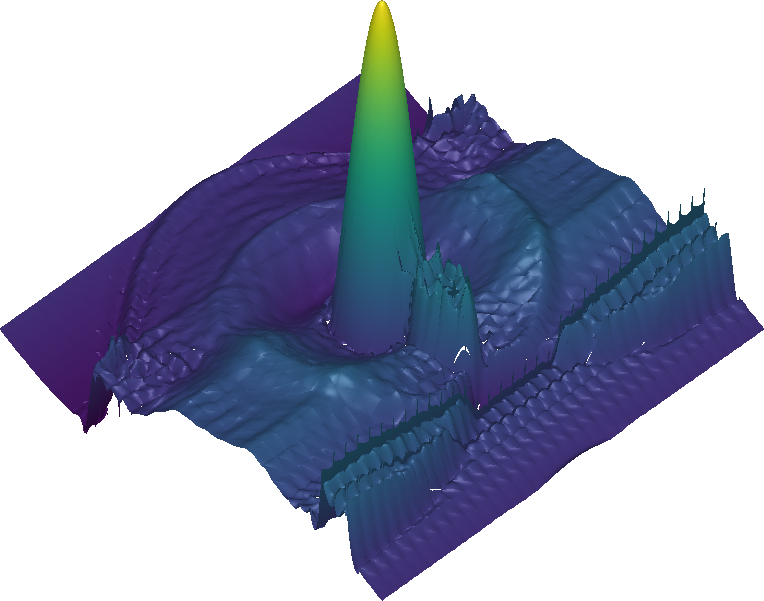}
\includegraphics[width=0.32\textwidth, height = 0.19\textheight]{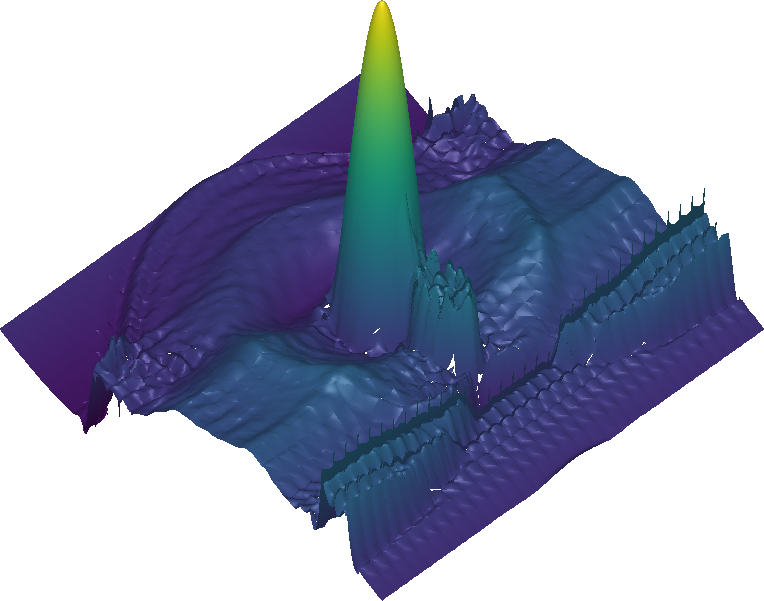}
\includegraphics[width=0.32\textwidth, height = 0.19\textheight]{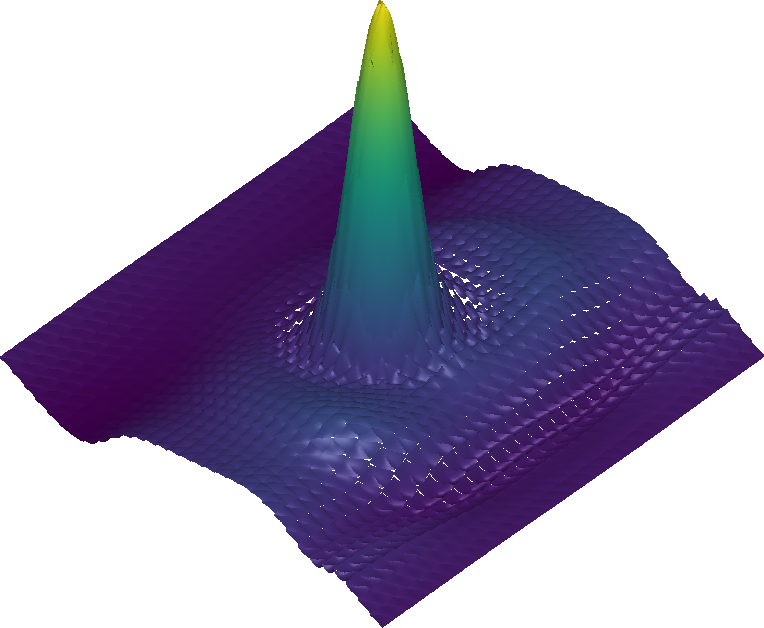}
\caption{Top view (top row) and side view (bottom row) of water surface and dry area for wave over a bump at T = 0.5. Left: node-wise limiting. Middle: element-wise limiting. Right: low order scheme}
\label{fig:CL_wave_2}
\end{center}
\end{figure}




\begin{figure}
\begin{center}
\includegraphics[width=0.32\textwidth, height = 0.2\textheight,
trim={0 0 6cm 0},clip]{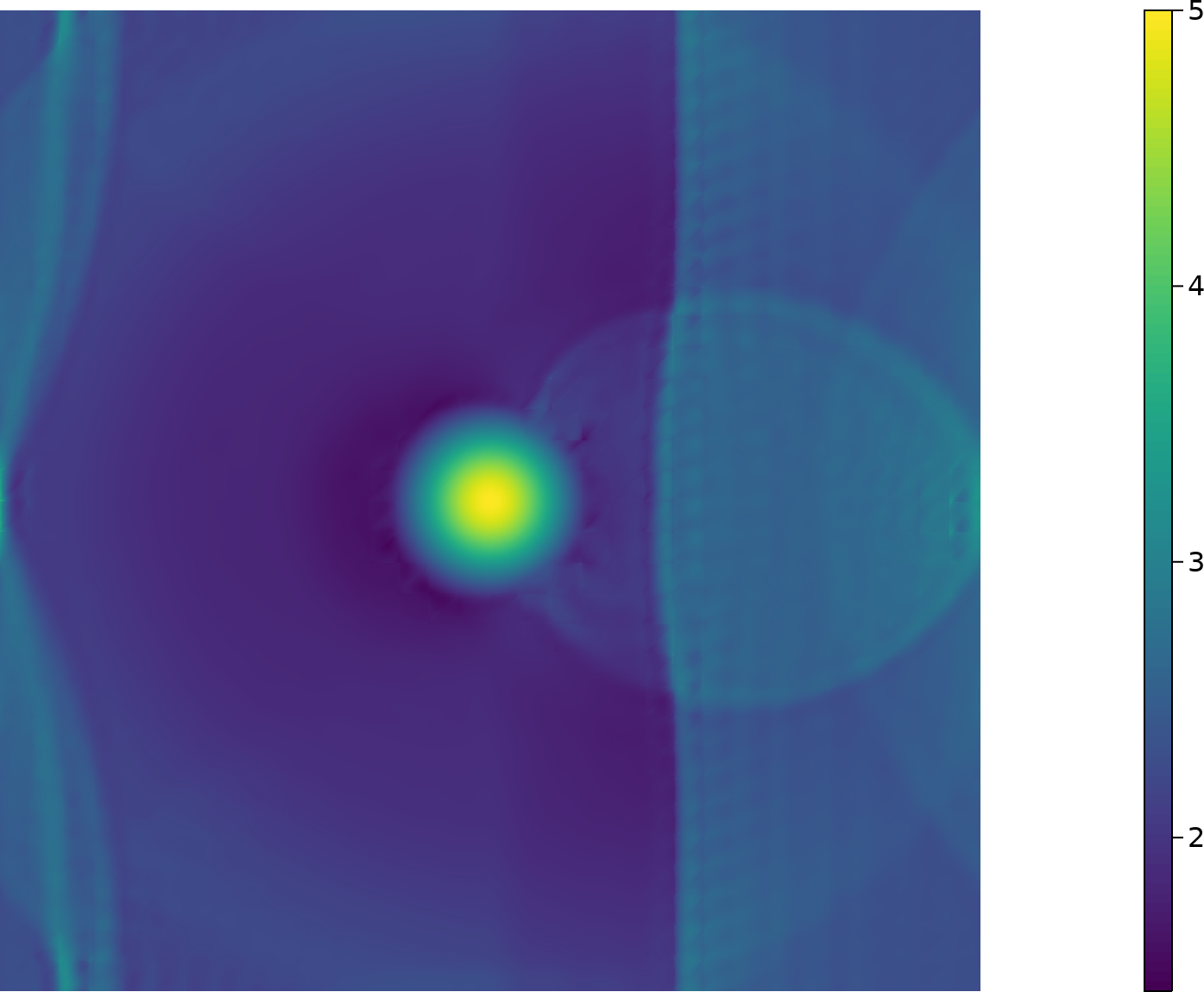}
\includegraphics[width=0.32\textwidth, height = 0.2\textheight, 
trim={0 0 6cm 0},clip]{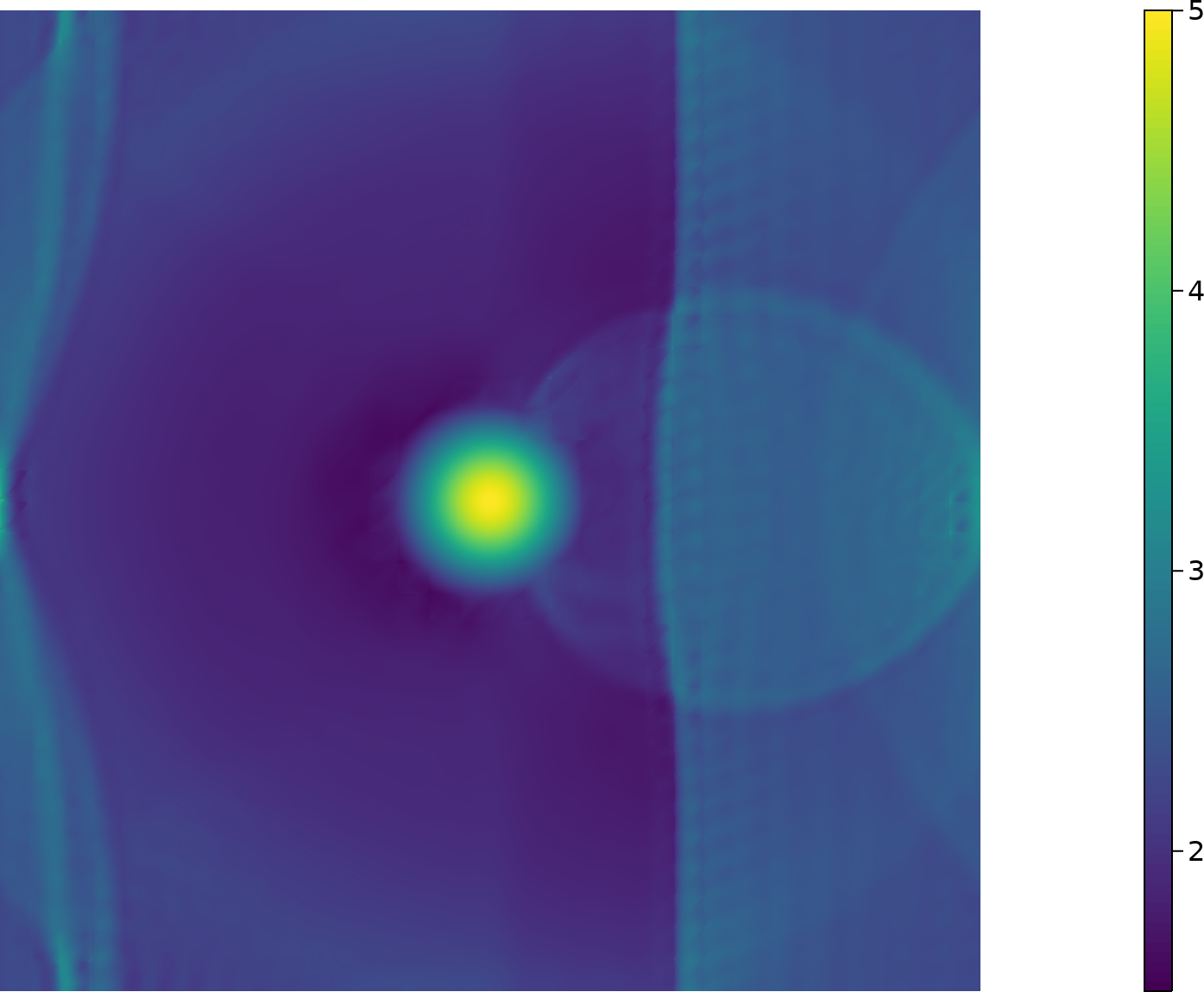}
\includegraphics[width=0.32\textwidth, height = 0.2\textheight, 
trim={0 0 6cm 0 cm},clip]{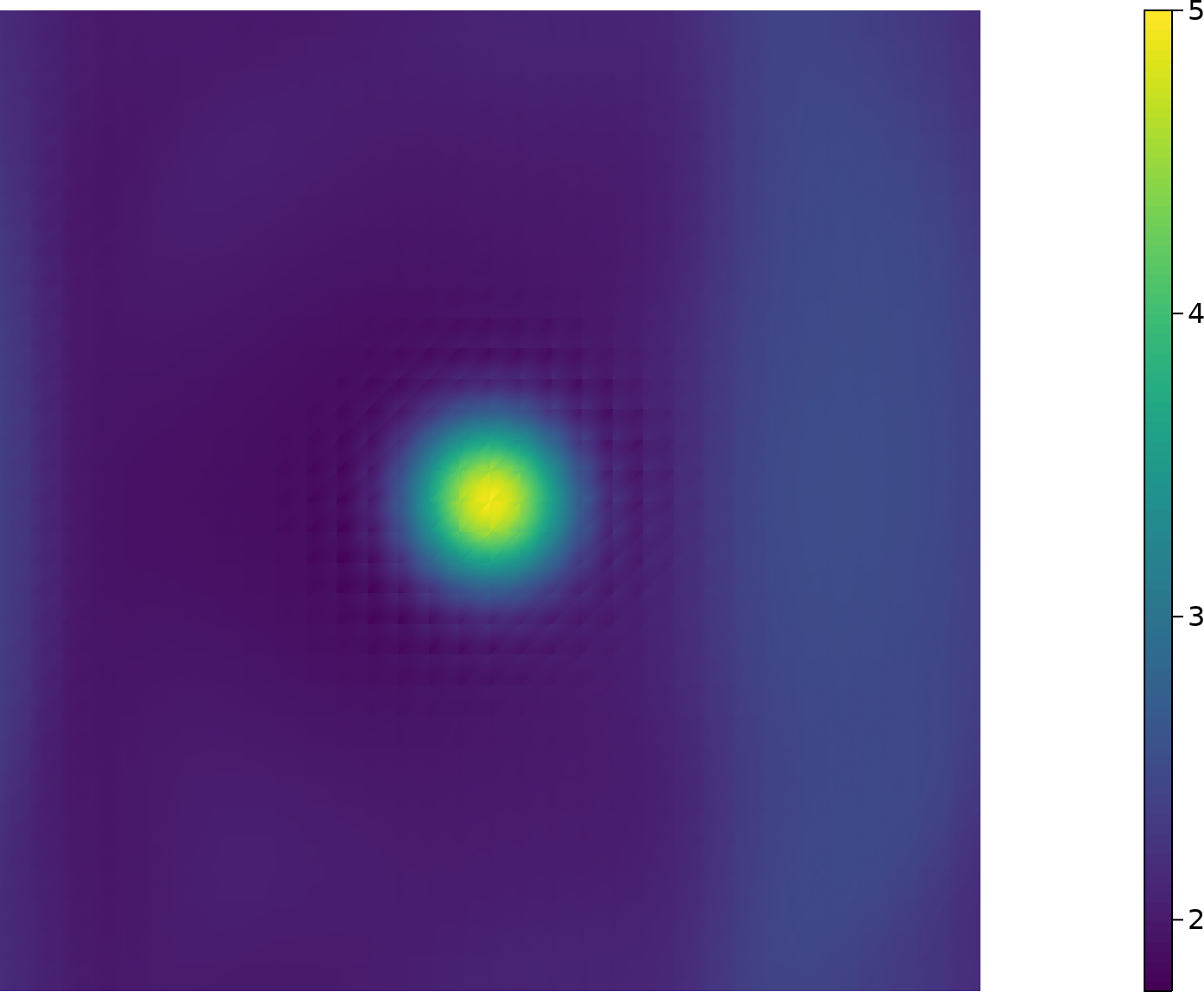}\\

\includegraphics[width=0.32\textwidth, height = 0.19\textheight]{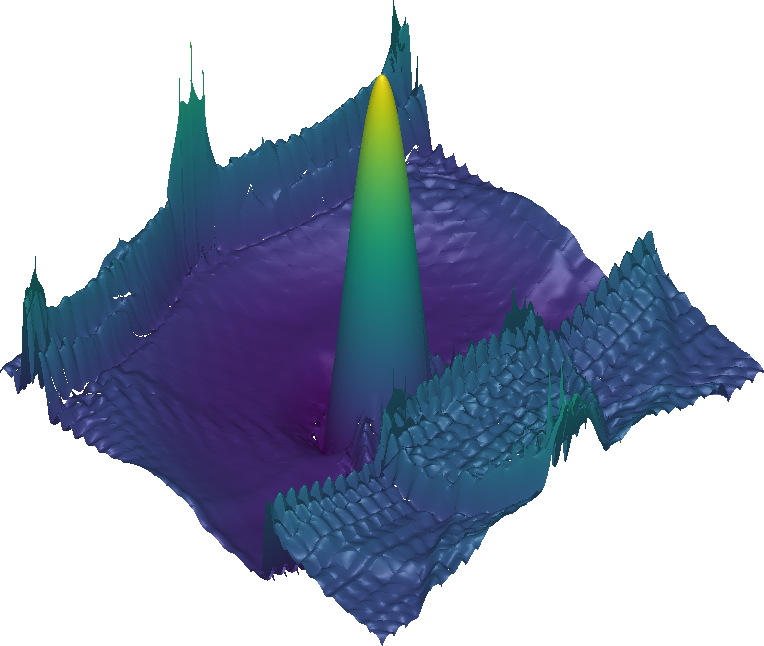}
\includegraphics[width=0.32\textwidth, height = 0.19\textheight]{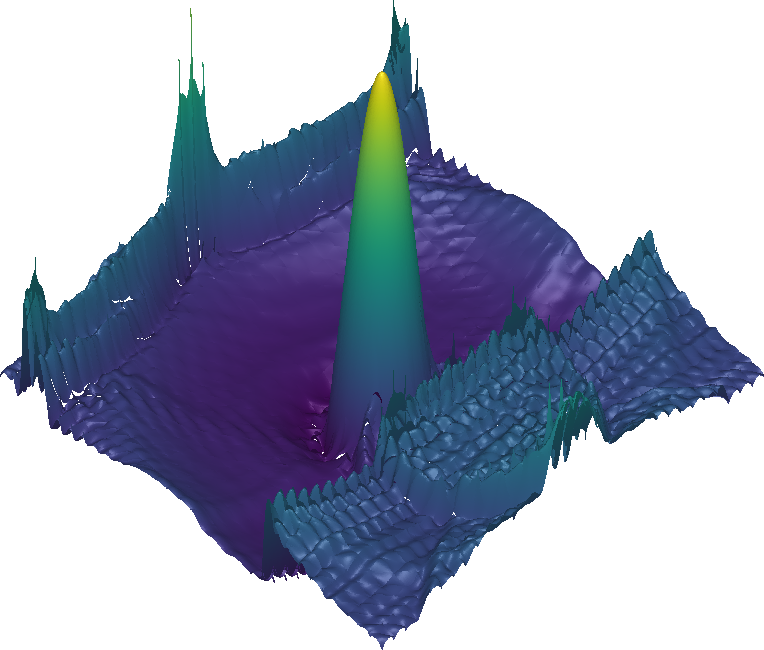}
\includegraphics[width=0.32\textwidth, height = 0.19\textheight]{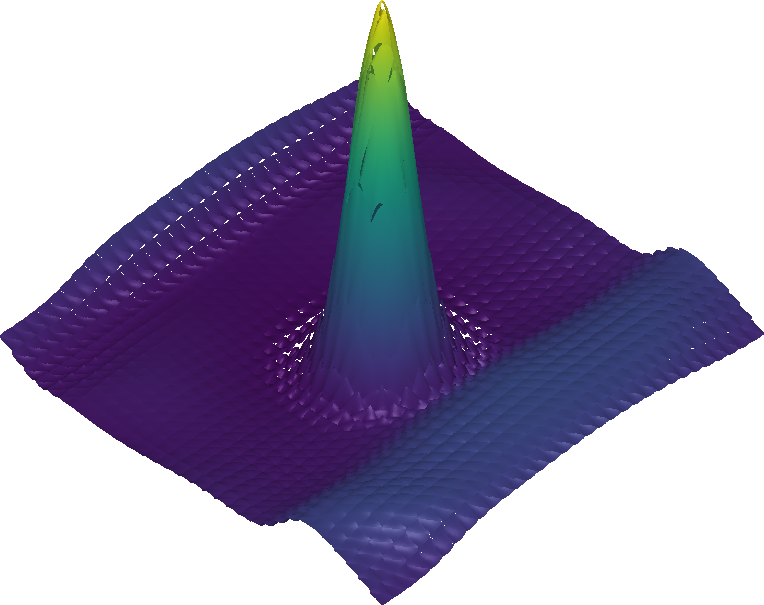}
\caption{Top view (top row) and side view (bottom row) of water surface and dry area for wave over a bump at T = 0.75. Left: node-wise limiting. Middle: element-wise limiting. Right: low order scheme}
\label{fig:CL_wave_3}
\end{center}
\end{figure}




\begin{figure}
\begin{center}
\includegraphics[width=0.32\textwidth, height = 0.2\textheight,
trim={0 0 6cm 0},clip]{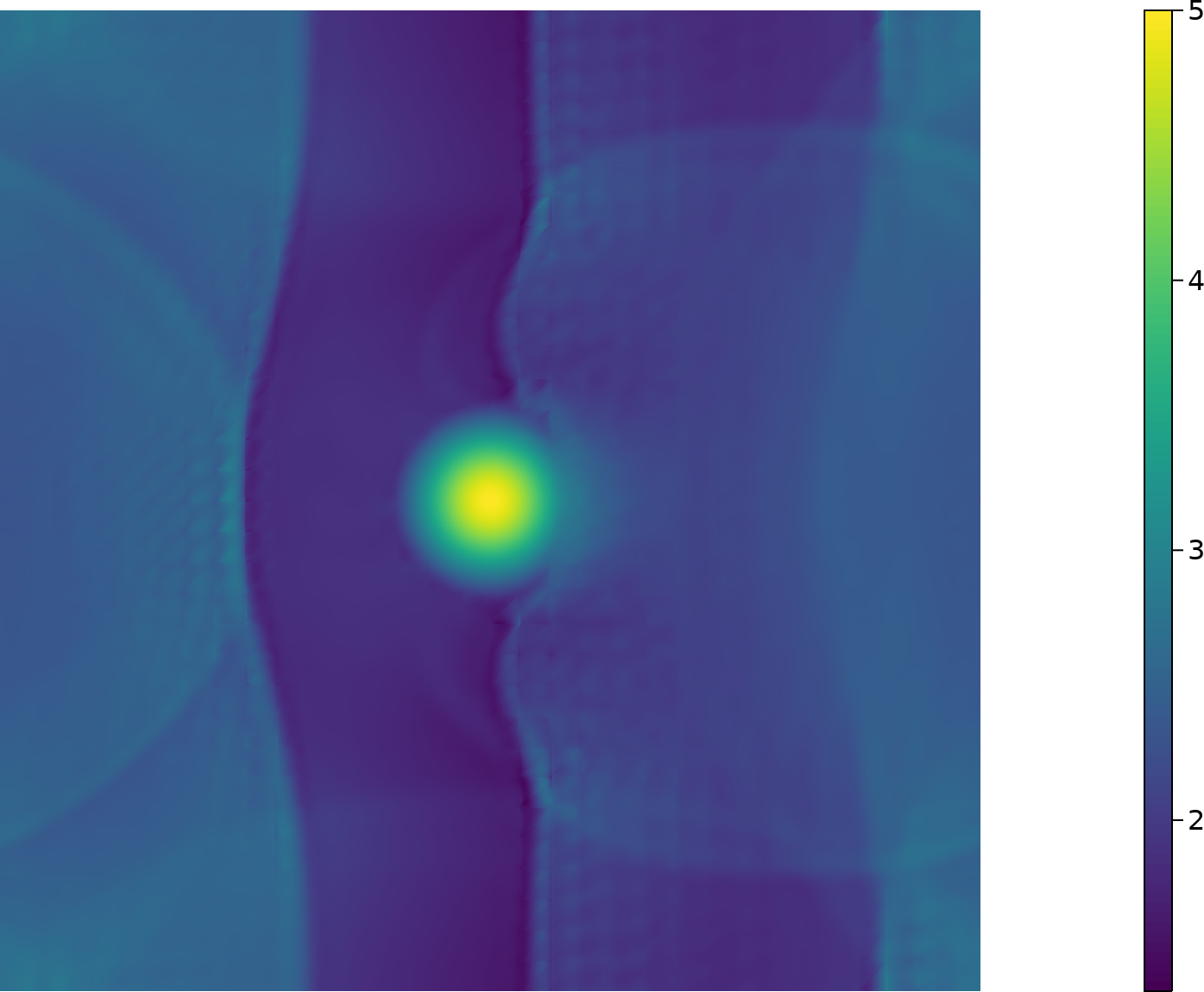}
\includegraphics[width=0.32\textwidth, height = 0.2\textheight, 
trim={0 0 6cm 0},clip]{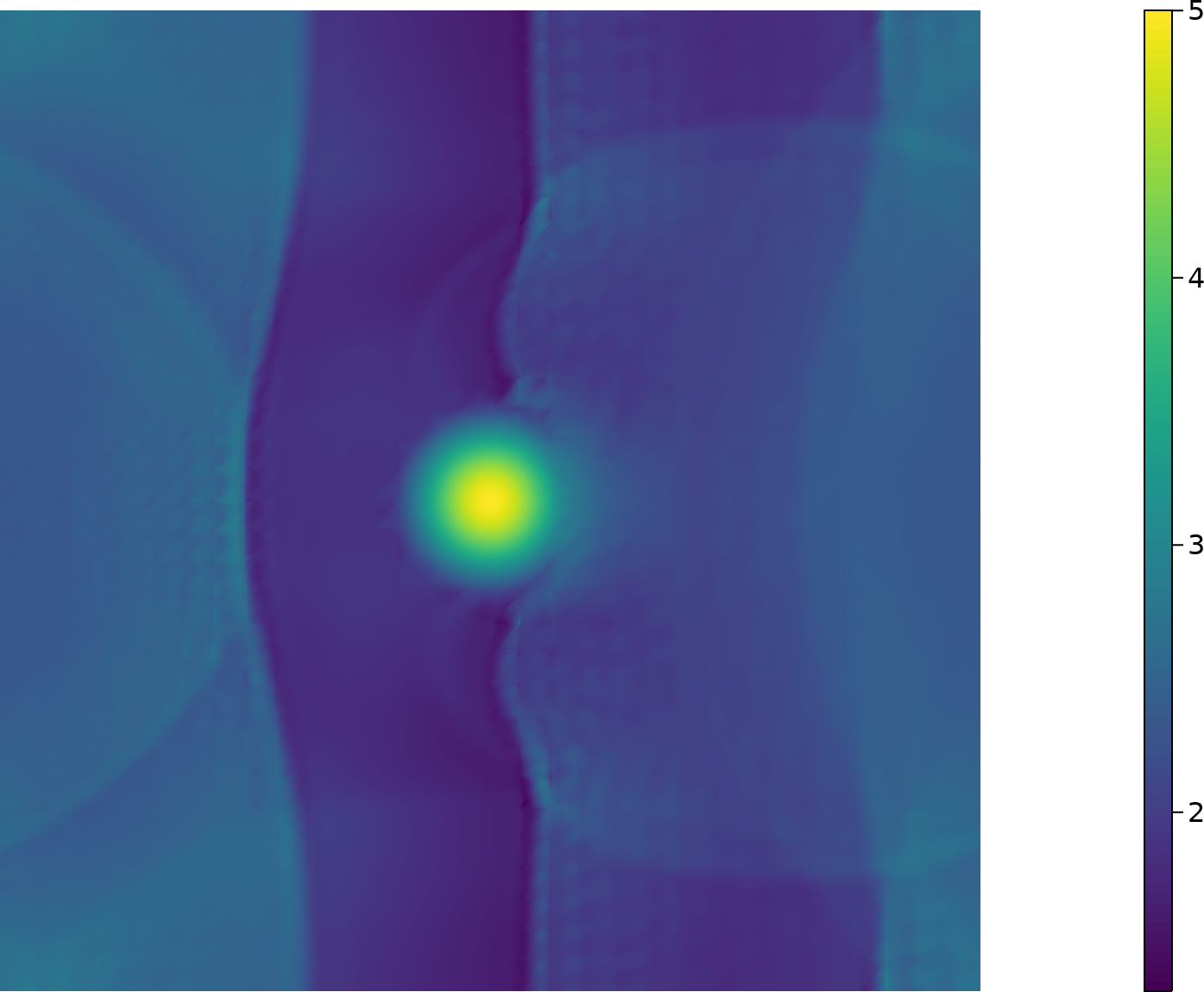}
\includegraphics[width=0.32\textwidth, height = 0.2\textheight, 
trim={0 0 6cm 0 cm},clip]{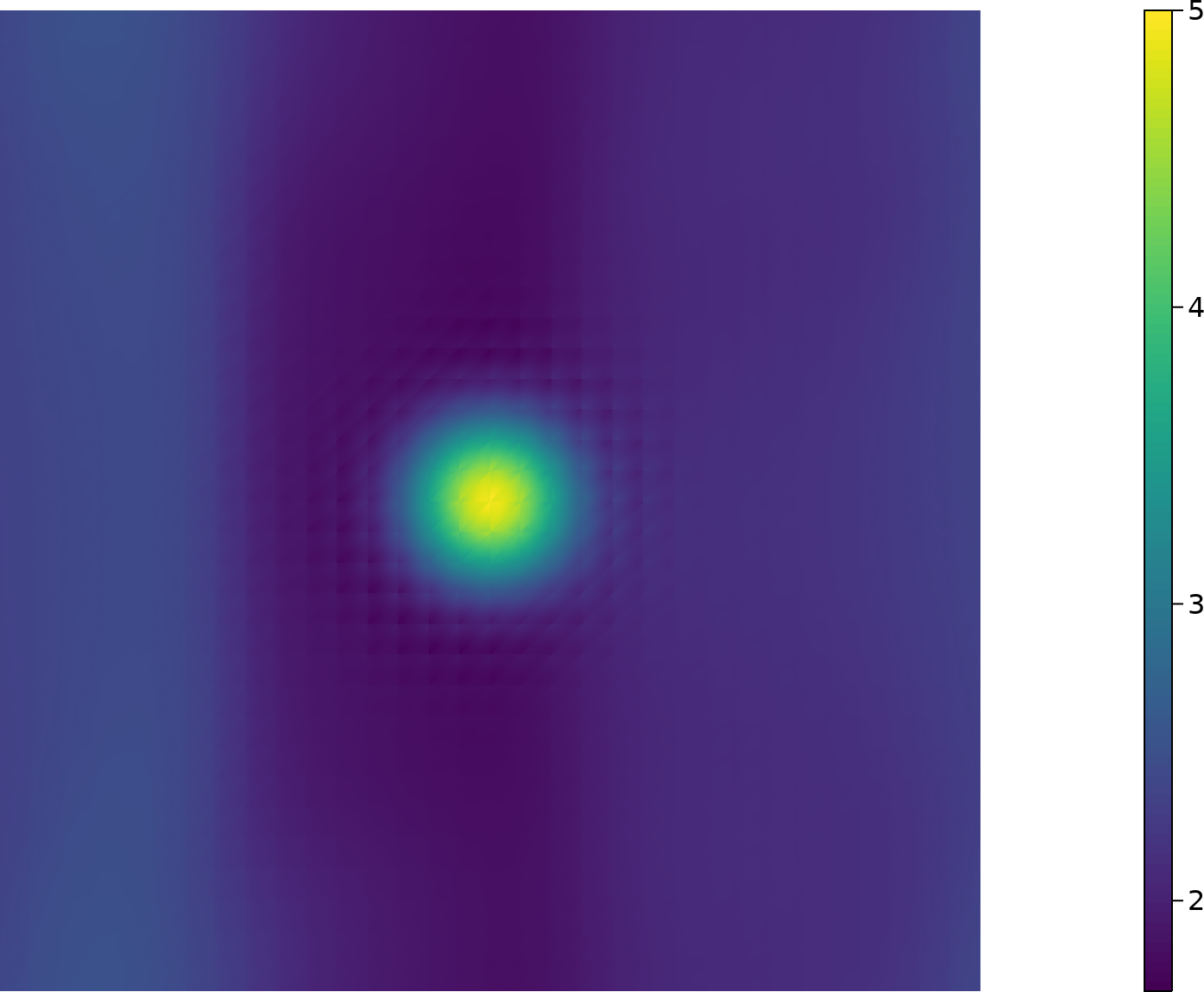}\\

\includegraphics[width=0.32\textwidth, height = 0.19\textheight]{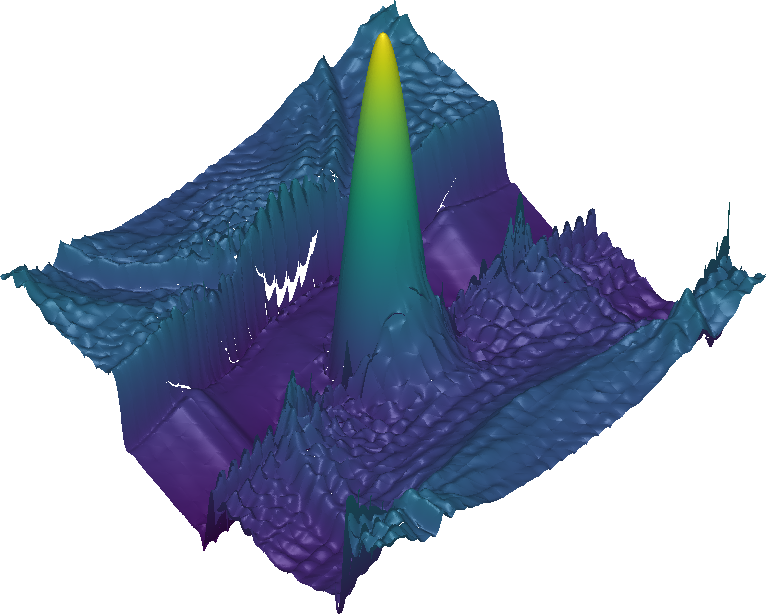}
\includegraphics[width=0.32\textwidth, height = 0.19\textheight]{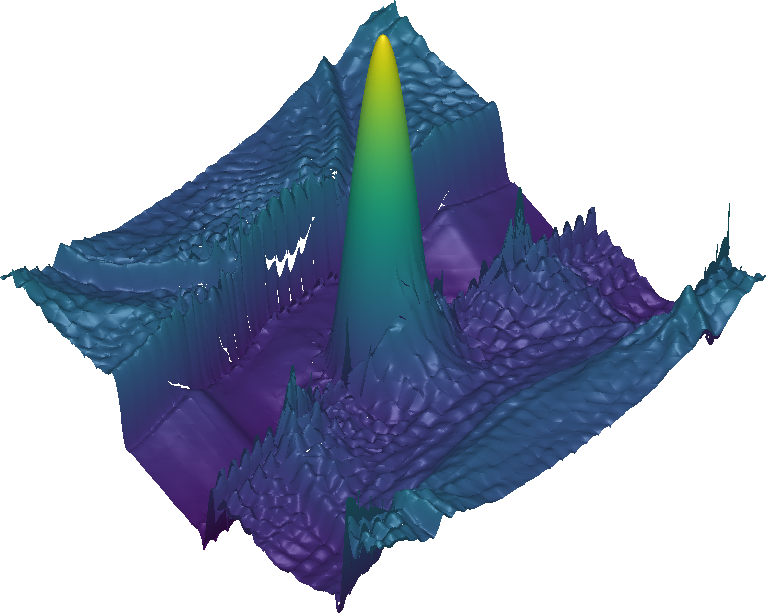}
\includegraphics[width=0.32\textwidth, height = 0.19\textheight]{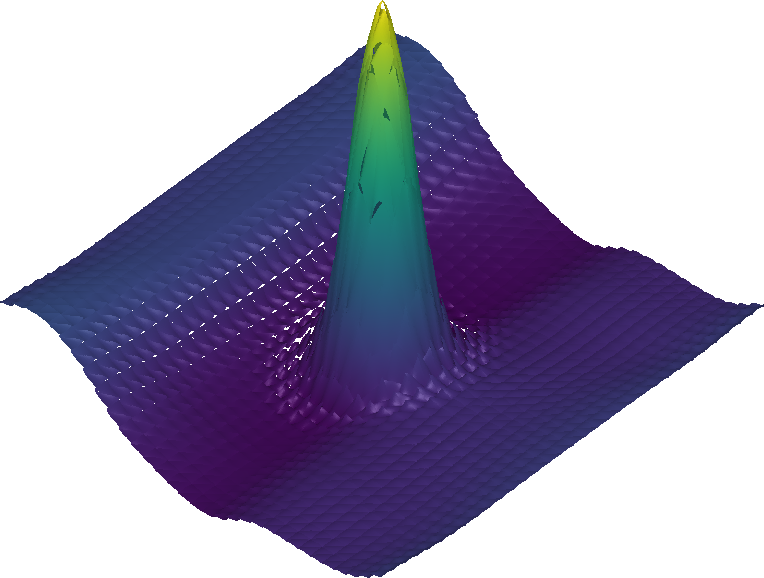}
\caption{Top view (top row) and side view (bottom row) of water surface and dry area for wave over a bump at T = 1. Left: node-wise limiting. Middle: element-wise limiting. Right: low order scheme}
\label{fig:CL_wave_4}
\end{center}
\end{figure}




\end{subequations}

\section{Conclusions}
\label{sec:6}
In this paper, we have presented a high order entropy stable positivity-preserving numerical scheme based on DG methods for the shallow water equations in 1D and on 2D with triangular meshes. This numerical method can preserve the positivity of water heights under an appropriate time step restriction. We apply the convex limiting technique on a low order positivity preserving method and a high order nodal SBP method to obtain this final scheme. This method can achieve high order of accuracy away from dry areas while remaining stable at the wet-dry front. We provide both 1D and 2D numerical experiments to demonstrate our theoretical findings.
\hspace{1em}

\section*{Acknowledgement}
Authors Xinhui Wu and Jesse Chan gratefully acknowledge support from the National Science Foundation under awards DMS-1719818, DMS-1712639, and DMS-CAREER-1943186.

Author Nathaniel Trask acknowledges funding under the Collaboratory on Mathematics and Physics-Informed Learning Machines for Multiscale and Multiphysics Problems (PhILMs) project funded by DOE Office of Science (Grant number DE-SC001924) and the DOE Early Career program. Sandia National Laboratories is a multi-mission laboratory managed and operated by National Technology and Engineering Solutions of Sandia, LLC., a wholly owned subsidiary of Honeywell International, Inc., for the U.S. Department of Energy’s National Nuclear Security Administration under contract DE-NA0003525. This paper describes objective technical results and analysis.  Any subjective views or opinions that might be expressed in the paper do not necessarily represent the views of the U.S. Department of Energy or the United States Government.

\appendix
\section*{Appendix A: Construction of low order operators}
\label{sec:AppendixA}
\renewcommand{\thesubsection}{A.\Roman {subsection}}
\renewcommand{\theequation}{A.\arabic{equation}}
In this appendix, we describe the procedures to compute lower order operators. We use a meshfree approach on the SBP nodes. We start by building a connectivity graph on the SBP nodes. Then, we solve a graph-Laplacian problem based on the connectivity graph to form the low order operators. 
\subsection{Construction of sparse connectivity graph}
\label{subsec:sparse_con_graph}
\begin{figure}
\begin{center}
\includegraphics[width=.45\textwidth]{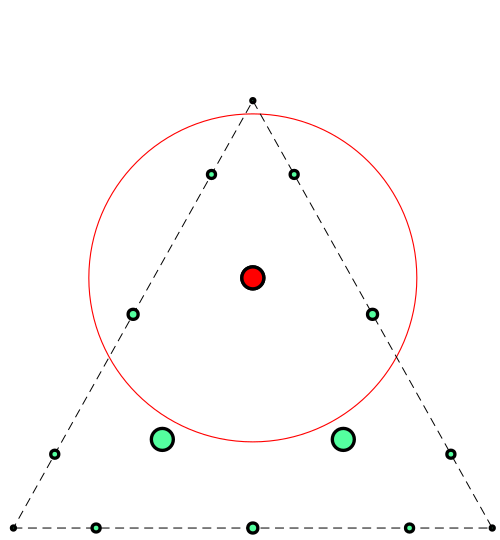}
\hspace{-0.2cm}
\includegraphics[width=.45\textwidth]{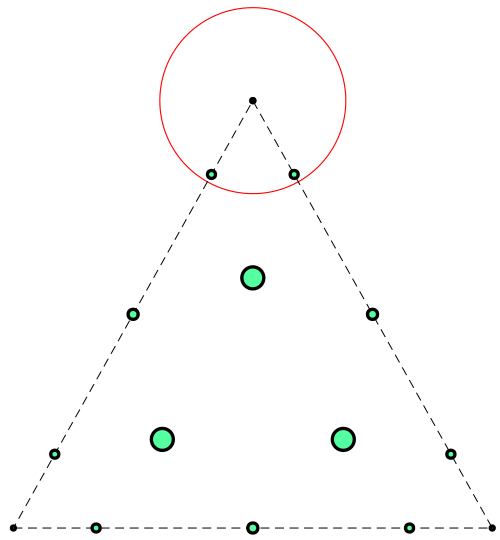}\\
\caption{Left: neighborhood of a interior nodes with larger weights. Right: neighborhood of a boundary node with smaller weights. The sizes of the nodes are proportional to its quadrature weights. The central nodes of the neighborhoods are marked in red.}
\label{fig:rescaled_radius}
\end{center} 
\end{figure}

To construct the low order operator $\fnt{Q}^{L}_{i}$, we need to introduce a connectivity graph. This connectivity defines the adjacency matrix $A$ from the previous section. We first use uniform radius $r = \alpha$ to define to determine the neighbors of each node, such that node $j \in I(i)$ if the distance between node $j$ and $i$ is less or equal to the radius $r$. However, this approach is less well suited to non-uniform node distributions, which is the case for the Gauss-Legendre quadrature node and Gauss-Lobatto quadrature nodes shown in Fig. \ref{fig:SBP_GLE_nodes} and \ref{fig:SBP_GLO_nodes}, respectively. Then, we investigate the scaling of the radii calculated from coefficient $\alpha$ and the quadrature weight of the node raised to some power. So now the radius $r_i$ for node $i$ can be written as
\begin{align}
    r_i = \alpha (w_i)^p,
\end{align}
where $w_i$ is the quadrature weight of the node $i$ and $p$ is some power. With our experimental approach, we test many different combinations of $\alpha$ and $p$ and discover that $\alpha = 1.0$ and $p = \frac{1}{4}$ would produce sparsest connectivity graph. For illustration, we plot the circles with radii $r_i$ on nodes with different quadrature weights in Fig. \ref{fig:rescaled_radius}. The red circles define the neighborhoods of nodes. We observe that nodes with greater weights form neighborhoods with larger radii.

\begin{figure}
\begin{center}
\subfloat[]{\includegraphics[width=.45\textwidth]{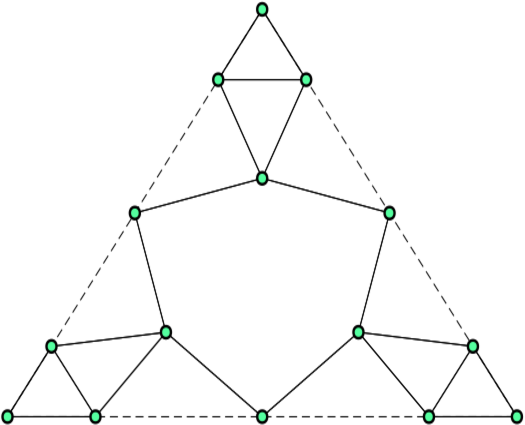}}
\hspace{0.5cm}
\subfloat[]{\includegraphics[width=.4\textwidth]{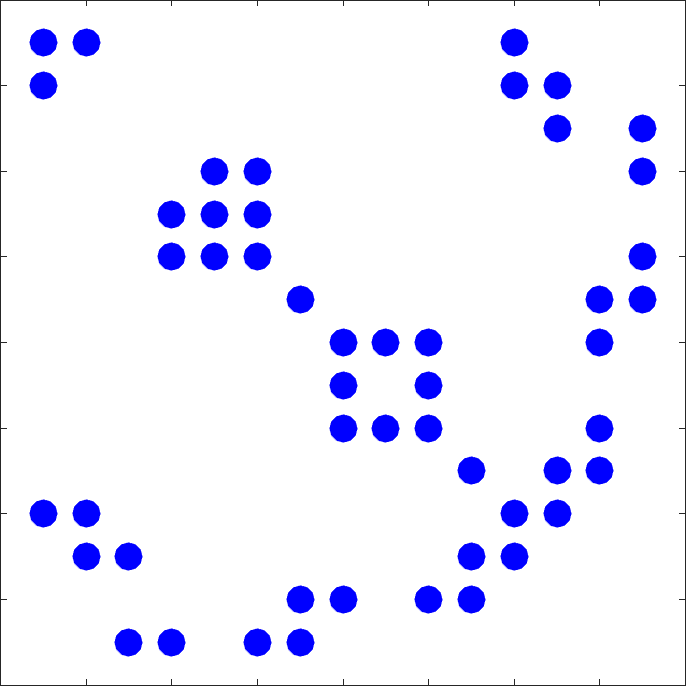}}\\
\subfloat[]{\includegraphics[width=.45\textwidth]{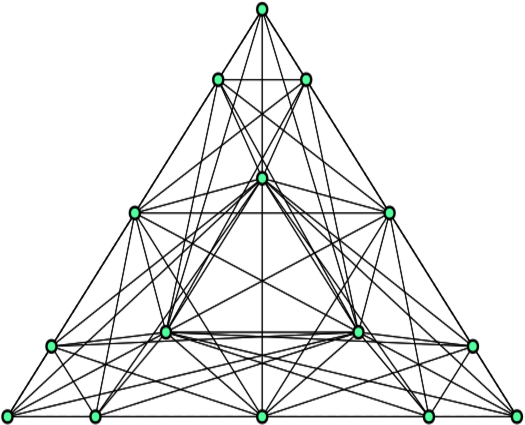}}
\hspace{0.5cm}
\subfloat[]{\includegraphics[width=.4\textwidth]{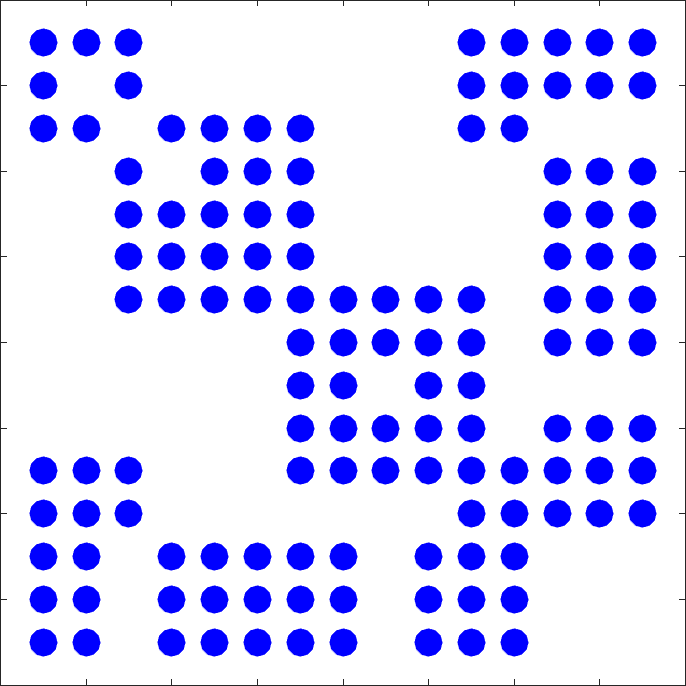}}
\caption{(a) Connectivity map of quadrature nodes with $\alpha = 1.0$. (b) Spy diagram for low order operator with $\alpha = 1.0$. (c) Connectivity map of quadrature nodes with $\alpha = 2.5$. (d) Spy diagram for low order operator with $\alpha = 2.5$.}
\label{fig:low_opt_sparse}
\end{center}
\end{figure}

Notice that we still need the low order operator constructed from this connectivity graph to satisfy the requirements in (\ref{eq:low_order_req}). Therefore, $\alpha$ can not be arbitrarily small. In Fig. \ref{fig:low_opt_sparse}, we notice that the smaller radius can produce a much sparser low order operator by reducing the number of connected nodes. Sparsity can increase the accuracy of the algorithm. It also improves the computing speed because of the fewer number of non-zero terms in the matrix.

\subsection{Meshfree construction of low order summation-by-parts matrices}
\label{subsec:meshfree_con}
This setup defines an undirected finite graph $(V, E)$, where the vertices consist of all the quadrature points.  $(j,k)$ is an edge of the graph, i.e., $(j,k) \in E$, iff $j \in I(k)$ and $k \in I(j)$. We use (V, E) to refer to the connectivity graph. Let $r$ denote the radius of a ball. Then we define a adjacency matrix $\fnt{A}$ where $\fnt{A}_{jk} = 1$, for $j\neq k$ if the distance between node $j$ and $k$ is less or equal to $r$ and $\fnt{A}_{jk} = 0$ otherwise. $j \in I(k)$ iff $\fnt{A}_{jk} = 1$, We let  $\fnt{A}_{jj} = |I(j)|$, the number of nodes adjacent to node $j$.

We first build a mesh-free low order differentiation operator $\fnt{Q}^{L}_{i}$ in the $i$th coordinate direction with the summation-by-part property. We require that
\begin{align}
    \fnt{Q}^{L}_{i} +  (\fnt{Q}^{L}_{i})^T &= \fnt{B}_i = {\rm diag}(\bm{n}_i\circ w^f_i),
    \label{eq:mesh_free_op_1_1}
\end{align}
\begin{align}
    \fnt{Q}^{L}_{i}\bm{1} &= \bm{0},
    \label{eq:mesh_free_op_1_2}
\end{align}
We can rewrite Eq. (\ref{eq:mesh_free_op_1_1}) as
\begin{align}
    \fnt{Q}^{L}_{i} = \bm{S} + \frac{1}{2}\fnt{B}_i, 
    \label{eq:mesh_free_op_2}
\end{align}
where $\bm{S}$ is skew-symmetric. We assume that $\bm{S}_{jk} = \phi_j-\phi_k$ for entries $\phi_j$ and $\phi_k$ of some ``potential" vector $\Phi$. Enforcing $\fnt{Q}^{L}_{i}\bm{1} = \bm{0}$ then implies that
\begin{align}
    -\bm{S}\bm{1} = \frac{1}{2}\fnt{B}_{i}\bm{1}
    \implies \underbrace{\sum_{k \in I(j)} (\phi_k-\phi_j)}_{-\bm{S}\bm{1}} = \frac{1}{2}(\fnt{B}_{i})_{jj}
    \label{eq:mesh_free_op_3}
\end{align}
This is equivalent to a Neumann-type graph-Laplacian problem in \cite{trask2020conservative}. We can define the graph Laplacian matrix $\fnt{L} = \fnt{D} - \fnt{A}$, where  $\fnt{D}$ is the diagonal degree matrix whose entries are the degrees of each vertex. Then (\ref{eq:mesh_free_op_3}) can be written as  $\fnt{L}\bm{\Phi} = \frac{1}{2}\fnt{B}\bm{1}$, and we solve the graph Laplacian problem subject to $\bm{\Phi}^T \bm{1} = \bm{0} $, where $\bm{\Phi}_j = \phi_j$. We solve the system with Lagrange multiplier and use the entries in vector $\bm{\Phi}$ to define the entries in $\bm{S}$. Last, we use Eq. (\ref{eq:mesh_free_op_2}) to obtain the differential operator $\fnt{Q}^{L}_{i}$.
\appendix
\section*{Appendix B: Effect of sparsity of the low order operators on numerical dissipation}
\label{sec:AppendixB}
\begin{figure}[H]
\begin{center}
\includegraphics[width=0.48\textwidth, trim={0 0 6cm 0},clip]{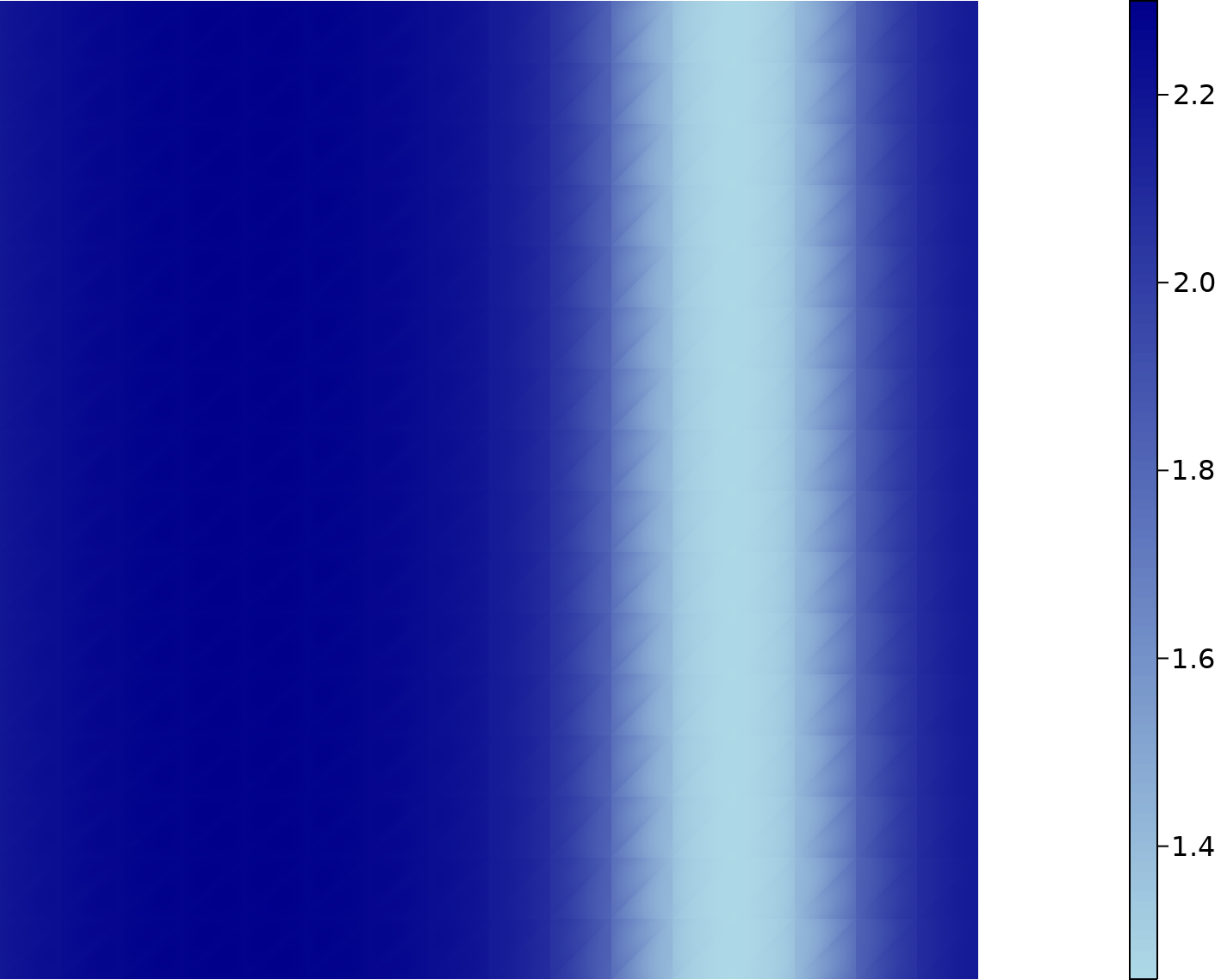}
\includegraphics[width=0.48\textwidth, trim={0 0 6cm 0},clip] {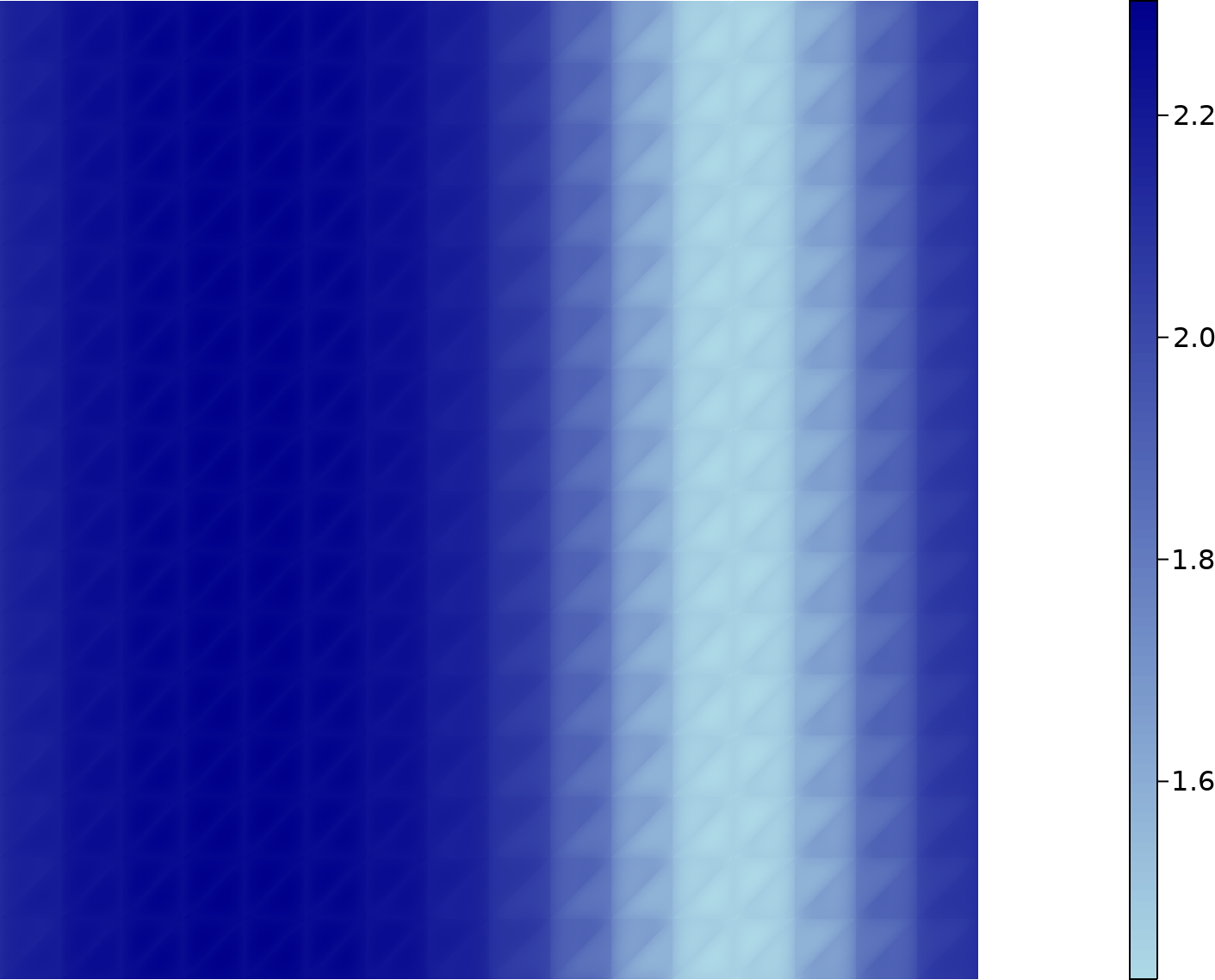}\\
\vspace{0.5cm}
\includegraphics[width=0.48\textwidth, height=0.3\textheight]{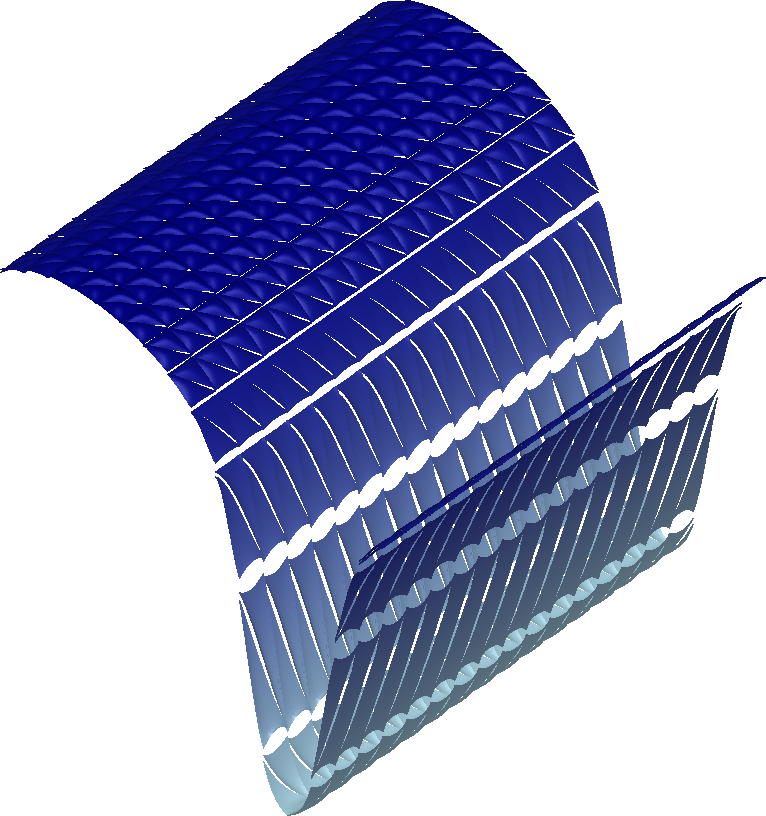}
\includegraphics[width=0.48\textwidth, height=0.285\textheight]{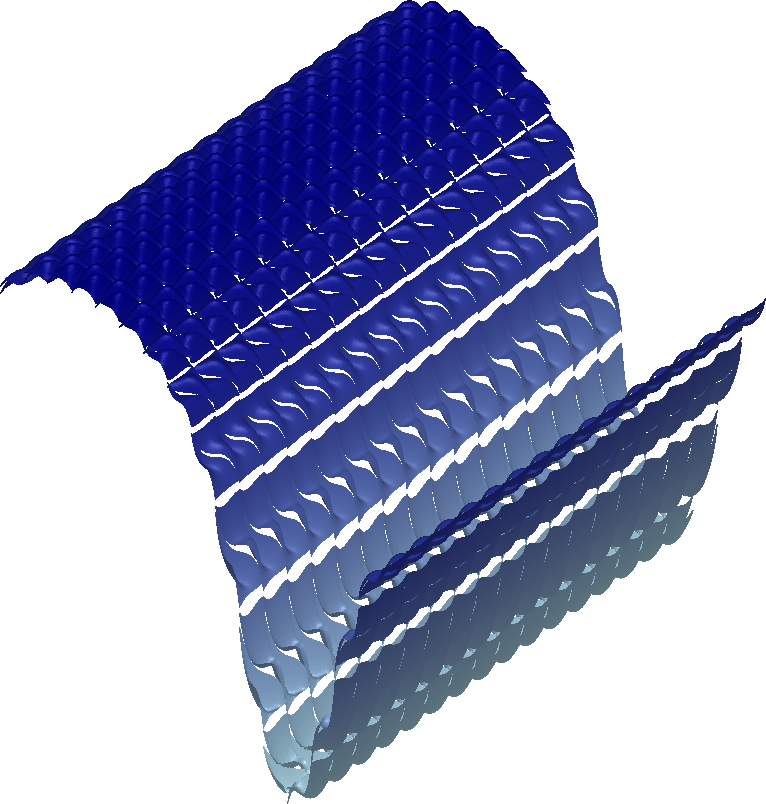}

\caption{Sine wave on low order scheme with top view (top row) and side view (bottom row) at $T=1$. Left: $\alpha = 1.0$. Right: $\alpha = 2.25$ }
\label{fig:alpha_low_order}
\end{center}
\end{figure}

\renewcommand{\thesubsection}{B.\Roman {subsection}}
\renewcommand{\theequation}{B.\arabic{equation}}
In this appendix, we demonstrate the effect of sparsity of low order operators on the quality of the low order positivity-preserving numerical solution. We test the low order scheme with low order operators constructed using different connectivity graphs. We utilize a sine wave initial condition on mesh defined on $[-1,1]\times[-1,1]$ with a $16\times16$ grid with periodic boundary conditions and a flat bottom topography. The initial conditions for this experiment are
\begin{align}
    &h(x,y) = \sin{(\pi x)}, \qquad hu(x,y) =  hv(x,y) = 0.
\end{align}
We run each low order scheme up to $T = 0.1$ with different connection radii $\alpha$ and present the results in Fig. \ref{fig:alpha_low_order}. We observe that the solutions become more dissipative as the radii $\alpha$ increases. The larger $\alpha$ becomes, the more connected the quadrature nodes are, resulting in graph viscosity being applied between more pairs of nodes. Choosing smaller values for $\alpha$ results in sparser low order operators, which sharpens the numerical solution and avoids introducing extraneous artificial dissipation.

\bibliographystyle{unsrt}
\bibliography{main.bib}

\end{document}